%% file: 513.tex


\input amstex
\input mathdefs \input citeadd

\sectno=-1   
\localtags
\ifx\shlhetal\undefinedcontrolsequence\let\shlhetal\relax\fi
\NoBlackBoxes
\define\mr{\medskip\roster}
\define\sn{\smallskip\noindent}
\define\mn{\medskip\noindent}
\define\bn{\bigskip\noindent}
\define\ub{\underbar}
\define\wilog{\text{without loss of generality}}
\define\ermn{\endroster\medskip\noindent}
\define\dbca{\dsize\bigcap}
\define\dbcu{\dsize\bigcup}
\define\nl{\newline}
\newbox\stickbox
\setbox 1 = \hbox{\raise1ex\hbox{$\bullet$}}
\setbox0\hbox to \wd1{\hss\vrule width 0.4pt\strut\hss}
\wd0=0cm
\setbox\stickbox\hbox{\box0\box 1}
\def\stick{\leavevmode\copy\stickbox}
\documentstyle {amsppt}
\topmatter
\title{pcf and infinite free subsets in an algebra \\
 Sh513} \endtitle
\author {Saharon Shelah \thanks {\null\newline I would like to thank 
Alice Leonhardt for the beautiful typing. \null\newline
 Research partially supported by the Israel Science Foundation \null\newline 
 First Typed at Rutgers - 98/July/16 \null\newline
 Done - 6-7/93 \null\newline
 Previous version - 98/Mar/19} \endthanks} \endauthor 
\affil{Institute of Mathematics\\
 The Hebrew University\\
 Jerusalem, Israel
 \medskip
 Rutgers University\\
 Mathematics Department\\
 New Brunswick, NJ  USA} \endaffil
\endtopmatter
\document  
\input alice2jlem

\newpage

\head {Anotated Content} \endhead  \resetall
\bn
\S1 $\quad$ Other variants of ``G.C.H. holds almost always"
\mr
\item "{{}}"  [We give another way to prove that for every $\lambda \ge
\beth_\omega$ for every large enough regular $\kappa < \beth_\omega$ we have
$\lambda^{[\kappa]} = \lambda$, dealing with sufficient conditions for
replacing $\beth_\omega$ by $\aleph_\omega$.]
\endroster
\bn
\S2 $\quad$ Large pcf$({\frak a})$ implies the existence of free sets
\mr
\item "{{}}"  [A nice example of the implication stated in the title is that
if pp$(\aleph_\omega) > \aleph_{\omega_1}$ then for every algebra $M$ of
cardinality $\aleph_\omega$ with countably many functions (or just
$< \aleph_\omega$ many), for some $a_n \in M$ (for $n < \omega$) we have
$a_n \notin c \ell_M(\{a_\ell:\ell \ne n,\ell < \omega\})$.  Generally if
pcf$({\frak a})$ is not just of cardinality $> |{\frak a}|$, but $\langle
J_{< \theta}[{\frak a}]:\theta \in \text{ pcf}({\frak a}) \rangle$ has large
rank (as defined below) than a relevant instance of IND connected to 
sup$({\frak a})$ holds.]
\endroster
\bn
\S3 $\quad$ Existence of free subsets implies restrictions on pcf
\mr
\item "{{}}"  [We have results of forms complementary to those of \S2 (though
not close enough).  So if IND$(\mu,\sigma)$ (in every algebra with universe
$\lambda$ and $\le \sigma$ functions there is an infinite independent subset)
then for no distinct regular $\lambda_i \in \text{ Reg} \backslash \mu^+$
(for $i < \kappa$) does $\dsize \prod_{i < \kappa} \lambda_i/[\kappa]
^{\le \sigma}$ have true cofinality.  We also look at IND$(\langle 
J^{\text{bd}}_{\kappa_n}:n < \omega \rangle),J_n$ an ideal on $\kappa_n$
(we ask for $\bar \alpha \in \dsize \prod_{n < \omega} \kappa_n$ such that
$n < \omega \Rightarrow \alpha_n > \text{ sup } c \ell_\mu(\{\alpha_\ell:
\ell \in (n,\omega)\}$) and more general version, and from assumptions as in
\S2 get results even for the non stationary ideal.]
\endroster
\bn
\S4 $\quad$ Sticks and Boolean Algebras
\mr
\item "{{}}"  [We deal with some other measurements of $[\lambda]
^{\ge \theta}$ and give an application by a construction of a Boolean
Algebra.]
\endroster
\bn
\S5 $\quad$  More on independence
\bn
\S6 $\quad$  Odds and ends
\mr
\item "{{}}"  [In 6.1 we deal with a replacement for $\Delta$-system lemma.
We have $> 2^\kappa$ sequences of length $\kappa$.  In \scite{6.3} we look at
how we can divide $F \subseteq \Pi_\alpha$ to few bounded sets.  In 
\scite{6.4} we relook at the characterization of a property of \cite{GHS},
generalizing the questions somewhat.  We then deal with freeness properties
for $F \subseteq {}^\delta\text{Ord}$ (modulo an ideal) and we give a correct
version of \cite[Ch.IX,3.5]{Sh:g} on characterizing cov$(\lambda,\lambda,
\theta,\sigma)$ when $\sigma > \aleph_0$ concerning the obtainment of the pp
version.  We shall continue in \cite{Sh:589}.]
\endroster
\newpage

\head {\S1 Other variants of ``G.C.H. holds almost always"} \endhead\resetall
\bigskip

We essentially redo the proof of \cite{Sh:460}, \S2 in another more general
way.
\demo{\stag{1.1} Notation}  1) ${\frak F}_\kappa(A)$ is the family of
$\kappa$-complete filters ${\frak D}$ on ${\Cal P}(A)$ (so ${\frak D} 
\subseteq {\Cal P}({\Cal P}(A))$ (so the points are subsets of $A$, and the
members of ${\frak D}$ (which are $\subseteq {\Cal P}(A)$) which we shall
be most interested in are ideal and their compliments. \nl
2) We say ${\frak D} \in {\frak F}_\kappa(A)$ has $\sigma$-complete
character if for any $Y \subseteq {\Cal P}(A)$ we have: $Y \in {\frak D}$ iff
fil$_\sigma(Y) \in {\frak D}$ where fil$_\sigma(Y)$ is the $\sigma$-complete
ideal on $A$ generated by $Y$. \nl
For an ideal $I$ on $X,I^+ = {\Cal P}(X) \backslash I$, similarly for a
filter.
\enddemo
\bigskip

\definition{\stag{1.2} Definition}  1) For ${\frak D} \in {\frak F}_\kappa(A)$,
cardinals $\mu < \lambda$ and $\sigma$ such that $|A| < \mu < \lambda$, we
say that $\lambda$ is $({\frak D},\mu,\sigma)$-inaccessible when: if
${\frak a}_t \subseteq (\mu,\lambda) \cap \text{ Reg}$ for $t \in A,
|{\frak a}_t| < \sigma$ \ub{then} $\{B \subseteq A:
\text{pcf}_{\sigma\text{-complete}}(\cup\{{\frak a}_t:t \in B\}) \subseteq
\lambda\} \in {\frak D}$. \nl
2) If we omit $\mu$ we mean: for $\mu = (|A| + \sigma)^+$.
\enddefinition
\bigskip

\proclaim{\stag{1.3} Theorem}  Suppose $\langle \kappa_n:n < \omega
\rangle$ is a strictly increasing sequence of regular cardinals $> \aleph_2$.
Stipulate $\kappa_{-1} = \aleph_1$ and assume ${\frak D}_n \in
{\frak F}_{\kappa_{n-1}}(\kappa_n)$ for $n < \omega$ and $\kappa =
\dsize \sum_{n < \omega} \kappa_n$ satisfies:
\mr
\item "{$\otimes$}"  if $n < \omega,\aleph_0 < \theta = \text{ cf}(\theta) <
\kappa_n,h:\kappa_{n=1} \rightarrow \theta \text{ and}$ \nl
$Y \in {\frak D}^+_{n+1}$ (so $Y \subseteq {\Cal P}(\kappa_{n+1}))$
\ub{then} for some $\zeta < \theta$ we have
{\roster
\itemitem{ $(*)_{y,\zeta}$ }  $\{B \in Y:\text{sup Rang}(h \restriction B)
< \zeta\} \in {\frak D}^+_{n+1}$.
\endroster}
\ermn
If $\lambda > \kappa$ \ub{then} for every $n < \omega$ large enough, $\lambda$
is $({\frak D}_n,\kappa,\aleph_1)$-inaccessible.
\endproclaim
\bigskip

\remark{\stag{1.4} Remark}  1) We can replace $\omega,\aleph_1$ by $\theta,
\theta^+$ or $< \theta,\theta$ (when $\theta$ is regular uncountable)
respectively (so $\kappa = \dsize \sum{i < \theta} \kappa_i$, etc.) (why?
repeat the proof or force by Levy$(\aleph_0,< \sigma))$.  Of course, we can
replace ${\frak F}(\kappa_n)$ by ${\frak F}(A)$ if $|A| = \kappa_n$. \nl
2) Note that the set defined in \scite{1.2}(2) is always an ideal on $A$.
\endremark
\bigskip

\demo{Proof}  We prove this by induction on $\lambda$. If $\lambda = \kappa^+$
this is empty (as $(\kappa,\lambda) = \emptyset$).  Also if $\lambda <
\kappa^{+ \omega_1}$ this is trivial, as $\dbcu_{t \in A} {\frak a}_t$ is
countable $(\subseteq \{\kappa^{+(\alpha+1)}:\alpha < \omega_1\} \cap
\lambda)$ hence pcf$_{\aleph_1\text{-complete}} (\dbcu_{t \in A} {\frak a}_t)
= \dbcu_{t \in A} {\frak a}_t \subseteq \{\kappa^{+(\alpha+1)}:\alpha <
\omega_1\} \cap \lambda$.  Also if this holds for $\lambda$ it holds for
$\lambda^+$ because pcf$({\frak a} \cup \{\lambda\}) \subseteq
\text{ pcf}({\frak a}) \cup \{\lambda\}$.  So we can assume that $\lambda$ is
a limit cardinal.  If the conclusion fails then for some infinite $W
\subseteq \omega$, for $n \in W$ we have a sequence $\langle {\frak a}^n
_\alpha:\alpha < \kappa_n \rangle$ (where ${\frak a}^n_\alpha \subseteq
(\kappa,\lambda) \cap \text{ Reg})$ which is a counterexample, i.e.
$Y_n =: \{B \subseteq \kappa_n:\lambda \nsupseteq 
\text{ pcf}_{\aleph_1\text{-complete}}(\dbcu_{\alpha \in B} {\frak a}^n
_\alpha)\} \in {\frak D}^+_n$.  If cf$(\lambda) > \kappa$ then $\cup\{
{\frak a}^n_\alpha:n < \omega,\alpha < \kappa_n\}$ is a subset of $\lambda$
of cardinality $\le \kappa$, hence is bounded by some $\lambda' < \lambda$,
so apply the induction hypothesis on $\lambda'$.  If $\aleph_0 <
\text{ cf}(\lambda) < \kappa$ let $\lambda = \sum\{\lambda_\zeta:\zeta <
\text{ cf}(\lambda)\},\dsize \bigwedge_{\zeta < \xi} \lambda_\zeta <
\lambda_\xi < \lambda$. Now as cf$(\kappa) = \aleph_0$ and $\kappa =
\dsize \sum_{n < \omega} \kappa_n$, for some $n(*) < \omega$ we have
cf$(\lambda) < \kappa_{n(*)}$. For every $n \in W \backslash (n(*)+2)$, we
define a function $h_n:\kappa_n \rightarrow \text{ cf}(\lambda)$ by $h_n
(\alpha) = \text{ Min}\{\zeta < \kappa_{n(*)}:{\frak a}_\alpha \subseteq
\lambda_\zeta\}$. Hence by the assumption $\otimes$, as cf$(\lambda) <
\kappa_{n(*)} < \kappa_{n-1}$, for some $\zeta_n < \text{ cf}(\lambda)$ we
have $\{B \subseteq \kappa_n:\dsize \bigwedge_{\alpha \in B} {\frak a}_\alpha
\subseteq \lambda_{\zeta_n}\} \in {\frak D}^+_n$.  Now we can contradict the
induction hypothesis for $\lambda' = \sup\{\lambda_{\zeta_n}:n \in W
\backslash (n(*)+2)\}$.  We are left with the case cf$(\lambda) = \aleph_0$
so let $\lambda = \dsize \sum_{n < \omega} \lambda_n,\lambda_0 = \kappa^+,
\lambda_n < \lambda_{n+1}$.  Note that as our ${\frak D}_m$ is
$\kappa_{m-1}$-complete (see \scite{1.1}) hence is $\aleph_1$-complete,
\wilog \, for some $\langle \lambda_k:k < \omega \rangle$ strictly
increasing and $\lambda = \dsize \sum_{k < \omega} \lambda_k,\lambda_0 =
\kappa^+$ define for each $n,k < \omega,Y^k_n = \{B \subseteq \kappa_n:\lambda
\nsupseteq \text{ pcf}_{\aleph_1\text{-complete}}(\dbcu_{\alpha \in B}
{\frak a}^n_i \cap [\lambda_k,\lambda_{k+1})\}$.  So $Y_n = \dbcu_{k < \omega}
Y^k_n$, but $Y_n \in {\frak D}^+_n$, and ${\frak D}^+_n$ is $\aleph_1$-
complete hence for some $k_n < \omega,Y^{k_n}_n \in {\frak D}^+_n$, so
possibly shrinking $W$ we get $\langle k_n:n \in W \rangle$ is constant or
strictly increasing, the former contradicts the induction hypothesis on
$\lambda_{k_{\text{Min}(W)^{+1}}}$.  Now renaming the $\lambda_k$'s we get
$k_n = n$ and we can replace ${\frak a}^n_\alpha$ by $a^n_\alpha \cap
[\lambda_n,\lambda_{n+1})$.  So \wilog \, Min$(W) > 4$ and for $n \in W$ we
have $\dbcu_\alpha {\frak a}^n_\alpha \subseteq [\lambda_n,\lambda_{n+1})$
and $\lambda_n < \lambda$, of course.

Let $n(*) = \text{ min}(W)$.  We try to define by induction on $k < \omega,
\langle \theta_t:t \in \omega_k \rangle,w_k = \dbcu_{i < \kappa_{n(*)}}
w_{k,i},J_k$ and $h_{k-1}$ if $k > 0$ such that:
\mr
\item "{$(a)$}"  $\theta_t \in \text{ Reg} \cap \lambda \backslash \kappa$
for $t \in w_k$
\sn
\item "{$(b)$}"  $w_k = \dbcu_{i < \kappa_{n(*)}} w_{k,i}$ is disjoint to
$\dbcu_{\ell < k} w_\ell$
\sn
\item "{$(c)$}"   $\langle w_{k,i}:i < \kappa_{n(*)} \rangle$ is a sequence
of pairwise disjoint, countable sets
sn
\item "{$(d)$}"  $w_{0,i} = \{i\} \times {\frak a}^{n(*)}_i$ and
$\theta_{(i,\tau)} = \tau$ for $\tau \in {\frak a}^{n(*)}_i$
\sn
\item "{$(e)$}"  $h_k$ is a function from $w_{k+1}$ to $w_k$ mapping
$w_{k+1,i}$ into $w_{k,i}$
\sn
\item "{$(f)$}"  $J_k = \{w \subseteq w_k:\lambda \supseteq
\text{ pcf}_{\aleph_2\text{-complete}}(\{\theta_t:t \in w\})\}$ is a proper
ideal
\sn
\item "{$(g)$}"  if $w \in J^+_k$ then $\{t \in w_{k+1}:h_k(t) \in w\} \in
J^+_{k+1}$
\sn
\item "{$(h)$}"  $t \in w_{k+1} \Rightarrow \theta_t < \theta_{h_k(t)}$.
\ermn
During the induction, $h_k$ is defined in the $k$-th step.

If we succeed, we shortly get a contradiction (by observation \relax
\cite[2.2]{Sh:460}).  For $k = 0$ define $w_{0,i},\theta_\tau$ for $\tau \in w_0 =
\dbcu_{i < \kappa_{n(*)}} w_{0,i}$ by clause (d) and the clause (f) holds as
otherwise $\{\theta_t:t \in w_0\}$ can be represented as $\dbcu_{\varepsilon
< \omega_1} {\frak b}_\varepsilon$ with max pcf$({\frak b}_\varepsilon) <
\lambda$, let $h:\kappa_{n(*)} \rightarrow \omega_1$ be $h(i) =
\text{ sup}\{\text{min}\{\varepsilon:\tau \in {\frak b}_\varepsilon\}:\tau
\in {\frak a}^{n(*)}_i\}$ and apply $\otimes$ from the hypothesis to get
$Y \subseteq Y_{n(*)}$ and $\zeta < \omega_1$ such that $Y \in {\frak D}^+
_{n(*)}$ and $\{B \in Y:\sup \text{ Rang}(h \restriction B) < \zeta\} \in
{\frak D}^+_{n(*)}$, but $B \in Y$ implies $\lambda \nsupseteq 
\text{ pcf}_{\aleph_1\text{-complete}}(\dbcu_{\alpha \in B} {\frak a}^{n(*)}
_\alpha)$ because $B \in Y_{n(*)}$ and

$$
\text{pcf}_{\aleph_1\text{-complete}}(\dbcu_{\alpha \in B} {\frak a}^{n(*)}
_\alpha) \subseteq \text{ pcf}_{\aleph_1\text{-complete}}
(\dbcu_{\varepsilon < \zeta} {\frak b}_\varepsilon) \subseteq
\dbcu_{\varepsilon < \zeta} \text{ pcf}_{\aleph_1\text{-complete}}
({\frak b}_\varepsilon) \subseteq \lambda;
$$
\mn
contradiction.  So assume $w_k,\langle w_{k,i}:i < \kappa_{n(*)} \rangle,
J_k$ are as required and we shall define $w_{k+1},\langle w_{k+1,i}:i <
\kappa_{n(*)} \rangle,J_{k+1},h_k$. \nl
Now for any $t \in w_k$ by the induction hypothesis for some $g(t) < \omega$
we have
\mr
\item "{$(*)_1$}"  if $m \in [g(t),\omega)$ and ${\frak b}_i \subseteq
\text{ Reg } \cap \theta_t \backslash \kappa$ is countable for $i < \kappa_m$
then
$$
\{B \subseteq \kappa_m:\text{pcf}_{\aleph_1\text{-complete}}
(\cup\{{\frak b}_i:i \in B\}) \subseteq \theta_t\} \in {\frak D}_m
$$
\item "{$(*)_2$}"  $g(t) > n(*) +1$.
\ermn
Let $u_{k,m} = \{t \in w_k:g(t) = m \text{ and } \theta_t > \kappa^+\}$.
We shall prove
\mr
\item "{$\boxtimes_{k,m}$}"  if $u = u_{k,m} \notin J_k$, then we can find
$\langle {\frak c}_t:t \in u_{k,m} \rangle,{\frak c}_t \subseteq \text{ Reg }
\cap \theta_t \backslash k$ countable such that: \nl
$u \subseteq u_{k,m},u \in J^+_\kappa$ implies
pcf$_{\aleph_1\text{-complete}}(\dbcu_{t \in u} {\frak c}_t) \nsubseteq
\lambda$.
\ermn
As $J_k$ is $\aleph_1$-complete (by its definition; even more) this suffices
for carrying the induction. \nl
[Why?  Let $\langle {\frak c}^m_t:t \in u_{k,m} \rangle$ for $m < \omega$
such that $u_{k,m} \notin J_k$ be as above, let

$$
w_{k+1,i} = \bigcup\{\{t\} \times {\frak c}^m_t:\text{for some } m,
u_{k,m} \notin J_k \text{ and } t \in u_{k,m} \cap w_{k,i}\},
$$
\mn
$\theta_{(t,\tau)} = \tau,w_{k+1} = \dbcu_{i < k_{n(*)}} w_{k+1,i}$ and we
define the function $h_{k+1}:w_{k+1} \rightarrow w_k$ by $h_{k+1}((t,
\sigma)) = t$ (note: every $t$ belongs to at most one $u_{k,m} \, 
(m < \omega)$ and $w_k \backslash \cup\{u_{k,m}:u_{k,m} \notin J_k\} =
\emptyset \text{ mod } J_k$.]
\enddemo
\bigskip

\demo{Proof of $\boxtimes_{k,m}$}  For each $\tau \in {\frak a}_m$, we
apply \cite[Ch.I,1.6]{Sh:g} or \cite[Ch.IX,4.1]{Sh:g} on $\langle \theta_t:t
\in u_{k,m} \rangle,J = J_k \restriction u_{k,m}$ and $\tau$ (possible as
$|u_{k,m}| < \kappa < \text{ min}\{\theta_t:t \in u_{k,m}\}$), each
$\theta_t (t \in u_{k,m})$ regular and $\dsize \prod_{t \in u_{k,m}}
\theta_t/J$ is $\lambda^+$-directed, $\tau < \lambda^+$ (the cases $\tau <
\theta_t$ can be ignored for several reasons, e.g. $\dbcu_\alpha
{\frak a}^m_\alpha \subseteq [\lambda_m,\lambda_{m+1}))$.  So we can find
$\langle \theta_{t,\tau}:\tau \in {\frak a}_m \rangle$ such that:
\mr
\item "{$(\alpha)$}"  $\theta_{t,\tau}$ is regular $\kappa^+ \le
\theta_{t,\tau} < \theta_t$
\sn
\item "{$(\beta)$}"  $\dsize \prod_{t \in u_{k,m}} \theta_{t,\tau}/J_k$
has true cofinality $\tau$
\ermn
(note that $t \in u_{k,m} \Rightarrow \kappa^+ < \theta_t$, so we can
assume $\theta_{t,\tau} \ge \kappa^+$).

Now for each $t \in u_{k,m},\theta_{t,\tau} \in \text{ Reg } \cap \theta_t
\backslash \kappa$ for $\tau \in {\frak a}_m$, but $g(t) = m$ (as $t \in
u_{k,m}$), hence by the definition of $g$, letting ${\frak a}^{m,t}_i =
\{\theta_{t,\tau}:\tau \in {\frak a}^m_i\}$, for $i < \kappa_m$, we have

$$
\Gamma^m_t =: \{B \subseteq \kappa_m:\theta_t \supseteq
\text{ pcf}_{\aleph_1\text{-complete}}(\dbcu_{i \in B} {\frak a}^{m,t}_i)\}
\in {\frak D}_m.
$$
\mn
But ${\frak D}_m$ is $\kappa_{n(*)+1}$-complete (as $m = g(t) > n(*)+1$) and
$|u_{k,m}| \le \kappa_{n(*)} < \kappa_{n(*)+1}$ hence $\Gamma^* = 
\dbca_{t \in u_{k,m}} \Gamma^m_t \in {\frak D}_m$.  On the other hand (by
the choice of ${\frak a}_m$)

$$
\Gamma_m =: \{B \subseteq \kappa_m:\lambda \supseteq
\text{ pcf}_{\aleph_1\text{-complete}}(\dbcu_{i \in B} {\frak a}^m_t)\}
\notin {\frak D}_m.
$$
\mn
So there is $B \in \Gamma_m \cap \Gamma^* = \Gamma_m \cap
\dbca_{t \in u_{k,m}} \Gamma^m_t$.

Let

$$
{\frak a}^* = \{\theta_t:t \in u_{k,m}\} \cup \{\theta_{t,\tau}:t \in
u_{k,m},\tau \in {\frak a}_m\} \cup {\frak a}_m,
$$
\mn
and for simplicity assume $\lambda \cap \text{ pcf}({\frak a}^*)| <
\text{ min}({\frak a}^*)$ hence there is a smooth close generating sequence
$\langle {\frak b}_\sigma[{\frak a}^*]:\sigma \in \lambda \cap
\text{ pcf}({\frak a}^*) \rangle$ (see e.g. \cite[6.7]{Sh:430}, if not use
\cite{Sh:430}, \scite{6.7} - \sciteu{6.7F}).  Clearly $B \ne \emptyset$.  For
each $t \in u_{k,m}$ we know $B \in \Gamma^m_t$ hence $\theta_t \supseteq
\text{ pcf}_{\aleph_1\text{-complete}}(\dbcu_{i \in B} {\frak a}^{m,t}_i)$.
So we can find countable ${\frak c}_t \subseteq \theta_t \cap
\text{ pcf}(\dbcu_{i \in B} {\frak a}^{m,t}_i)$ such that

$$
\dbcu_{i \in B} {\frak a}^{m,t}_i \subseteq \dbcu_{\sigma \in {\frak c}_t}
{\frak b}_\sigma[{\frak a}^*].
$$
\mn
The pcf calculus verifies clause (g) (as in \cite[\S2]{Sh:460}).
\hfill$\square_{\scite{1.3}}$
\enddemo
\bn
It is now natural to look for suitable filters ${\frak D}$, the simplest
ones are.
\definition{\stag{1.5} Definition}  For $\sigma < \theta < \kappa$ and $\mu
\le \kappa$ (always $\sigma,\theta$ regular) let ${\frak D} =
{\frak D}_{\sigma,\theta,\kappa,\mu}$ be the following filter on ${\Cal P}
(\kappa):Y \in {\frak D}$ \ub{iff} there are functions $f_\alpha:\kappa
\rightarrow \theta_\alpha$ for $\alpha < \chi(*) < \theta$ where
$\theta_\alpha \in [\sigma,\theta) \cap \text{ Reg}$ such that $Y \supseteq
\{a \subseteq \kappa:|a| \ge \mu$ and for every $\alpha < \theta$ for some
$\zeta < \theta_\alpha$ we have Rang$(f_\alpha \restriction a) \subseteq
\zeta\}$.  If $\mu = \theta$ we may omit it.
\enddefinition
\bigskip

\demo{\stag{1.6} Observation}  1) If $\sigma < \theta < \kappa_1 \le
\kappa_2$ and $\emptyset \in {\frak D}_{\sigma,\theta,\kappa_1,\mu}$ then
$\emptyset \notin {\frak D}_{\sigma,\theta,\kappa_2,\mu}$. \nl
2) ${\frak D}_{\sigma,\theta,\kappa}$ is a $\theta$-complete filter. \nl
3) If $2^{< \theta} < \kappa$ then $\emptyset \notin 
{\frak D}_{\sigma,\theta,\kappa}$ and this is preserved by $\sigma$-c.c.
forcing.
\enddemo
\bigskip

\demo{Proof}  Straight.
\enddemo
\bigskip

\demo{\stag{1.7} Conclusion}  Let $\mu$ be a limit singular cardinal of
cofinality $< \sigma < \mu$ and:
\mr
\item "{$\otimes$}"  for every $\theta \in (\sigma,\mu)$ for some $\kappa
\in (\theta,\mu)$ we have: $\emptyset \notin {\frak D}_{\sigma,\theta,
\kappa}$.
\ermn
\ub{Then} for every $\lambda > \mu$, for some $\theta = \theta_\mu \in
(\sigma,\mu)$ for every $\kappa \in (\theta,\mu)$ we have
\mr
\item "{$(*)$}"  if $\lambda_i \in (\mu,\lambda) \cap \text{ Reg}$ for $i <
\kappa$ then
$$
\{a \subseteq \kappa:\text{pcf}_{\sigma\text{-complete}}\{\lambda_i:i \in
a\} \subseteq \lambda\} \in {\frak D}_{\sigma,\theta,\kappa}.
$$
\endroster
\enddemo
\bigskip

\demo{Proof}  Assume $\lambda$ is a counterexample.  Without loss of generality$\sigma$ is regular, choose by induction on $\zeta < \sigma,\kappa_\zeta \in
(\sigma,\mu) \cap \text{ Reg}$ as follows
\mr
\item "{{}}"  $\kappa_0 \in (\sigma,\mu) \cap$ Reg arbitrary;
\sn
\item "{{}}"  $\kappa_\zeta \in (\dbcu_{\epsilon < \zeta} \kappa^+_\epsilon,
\mu) \cap$ Reg is minimal $\kappa$ which is a counterexample to $\otimes$ for
\sn
\item "{{}}"  $\theta = \dbcu_{\epsilon < \zeta} \kappa^+_\epsilon$ (in
particular $\emptyset \notin {\frak D}_{\sigma,\cup\{\kappa^+_\epsilon:
\epsilon < \zeta\},\kappa_\zeta}$, so $\kappa_\zeta < \mu$).
\ermn
Let $\kappa = \dbcu_{\zeta < \sigma} \kappa_\zeta$ and apply \scite{1.3} for
${\frak D}_{\sigma,\dbcu_{\epsilon < \zeta} \kappa^+_\epsilon,\kappa_\zeta}$
from Definition \scite{1.5}, more exactly the variant with $\sigma$ instead
$\aleph_1$ (see \scite{1.4}).  \hfill$\square_{\scite{1.7}}$
\enddemo
\bigskip

\remark{\stag{1.8} Remark}  1) In the proof of \scite{1.3} we can change the
universe during the proof, so weaken the demand $\otimes$. \nl
2) The problematic example is: $T \subseteq {}^{\omega_1}\omega$ a Kurepa
tree, say $T \cap {}^\alpha 2 = \{\gamma_\alpha(n):m < \omega\},\eta_j \in
\text{ lim}_{\omega_1}(T)$ for $j < j^*$ and $\{a \subseteq \omega_1 \cup
j^*:\delta =: a \cap \omega_1$ a limit ordinal and for every $n < \omega$
for some $j \in j^* \cap a$ we have $\gamma_\alpha(n) = \eta_j(\alpha)\} \in
{\Cal D}_{< \aleph_1}(\omega_1 \cup j^*)$. \nl
3) We can replace in \scite{1.7} above cf$(\mu) < \sigma$ by cf$(\mu) \ne
\sigma$.  
We can replace $\emptyset \notin {\frak D}_{\sigma,\theta,\kappa}$ by
$\emptyset \in {\frak D}_{\sigma,\theta,\kappa;\Upsilon}$ for any fix
$\Upsilon \in [\sigma,\mu)$.  Note that the case $\Upsilon < \sigma$ is not
interesting. \nl
4) Note that the meaning of $\emptyset \in 
{\frak D}_{\sigma,\theta,\kappa;\Upsilon}$ is that there are $\alpha(*) <
\theta$ and functions $f_\alpha:\kappa \rightarrow \theta_\alpha$ where
$\theta_\alpha \in \text{ Reg } \cap [\sigma,\theta)$ such that for no $u\in
[\kappa]^\Upsilon$ do we have $\alpha < \alpha(*) \Rightarrow
\text{ sup Rang}(f_\alpha \restriction u) < \theta_\alpha$.  Recall if
$\Upsilon$ is omitted it means $\Upsilon = \theta$.
\endremark
\bigskip

\definition{\stag{1.9} Definition}  Assume $\bold J \subseteq \text{ Id}
(\kappa)$ (= the family of ideals on $\kappa$). \nl
1) We say $(\lambda,\mu)$ is $\bold J$-inaccessible if $\kappa \le \mu <
\lambda$ and there are no $\lambda_i \in \text{ Reg } \cap (\mu,\lambda)$
for $i < \kappa$ and $J \in \bold J$ such that $\dsize \prod_{i < \kappa}
\lambda_i/\bold J$ is $\lambda$-directed (equivalently for some such
$\lambda_i$'s, $\dsize \prod_{i < \kappa} \lambda_i/J$ has true cofinality
and it is $\ge \lambda$). \nl
2) $(\lambda,*)$ means $(\lambda,\mu)$ for some $\mu \in [\kappa,\lambda),
\lambda$ means $(\lambda,\kappa)$. \nl
3) $\bold J$ is $\sigma$-indecomposable when: if $J \in \bold J$ and $h:
\text{Dom}(J) \rightarrow \sigma$ then for some $\zeta < \sigma$ and $I \in
\bold J$ we have $J \restriction h^{-1}\{\varepsilon:\varepsilon < \zeta\}
\le^* I$ (see below). \nl
4) For ideals $I_\ell$, on $A_\ell \, (\ell = 0,1) I_0 \le^* I_1$ if there
is $B_0 \in I^+_0$ and $B_1 \in I^+_1$ and one-to-one function $g$ from
$B_0$ into $B_1$ such that

$$
Y \cap B_0 \in I_0 \Rightarrow g''(Y \cap B_0) \in J_1
$$
\mn
5) $\bold J$ is $[\sigma,\kappa)$-indecomposable if it is 
$\theta$-indecomposable for every $\theta \in [\sigma,\kappa) \cap$ Reg.
\enddefinition
\newpage

\head {\S2 large pcf$({\frak a})$ implies the existence of free sets}
\endhead  \resetall  
\bigskip

\definition{\stag{2.1} Definition}  1) Let $\bar A = \langle A_\alpha:\alpha
< \alpha^* \rangle$ be a sequence of subsets of $\kappa$, no
$A_\alpha$ in the ideal generated by $\{A_\beta:\beta < \alpha\}$.  We
define functions $\text{rk } = \text{ rk}_{\bar A}$, rk$' = \text{
rk}'_{\bar A}$ from ${\Cal P}(\kappa)$ to the ordinals by:
\mr
\item "{{}}"  $\text{rk}(A) \ge \zeta$ \ub{iff} for every $\xi < \zeta$ for
some $\alpha,A \ne A \cap A_\alpha$ and $\text{rk}(A \cap A_\alpha) \ge
\zeta$
\sn
\item "{{}}"  $\text{rk}'(A) \ge \zeta$ \ub{iff} for every $\xi < \zeta$ for
some $\alpha$, we have rk$'(A \cap A_\alpha) \ge \xi$ and $A \backslash
A_\alpha,A \cap A_\alpha$ are not in id$_{\bar A \restriction \alpha} =$ the
ideal generated by $\{A_\beta:\beta < \alpha\}$.
\ermn
2) Let $\bar J = \langle (J_\alpha,J'_\alpha):\alpha < \alpha^* \rangle$ be
a sequence of pair of ideals on $\kappa$ such that $[\alpha < \beta
\Rightarrow J_\alpha \subseteq J'_\alpha \subseteq J_\beta \subseteq
J'_\beta]$ and for some $A_\alpha \in J^+_\alpha,J'_\alpha = J_\alpha +
A_\alpha$ we define rk$'_{\bar J}(A)$ for $A \subseteq \kappa$ by:
\mr
\item "{{}}"   rk$'_{\bar J}(A) \ge \zeta$ \ub{iff} for every $\xi < \zeta$
for some $\alpha < \alpha^*$ we have: rk$'_{\bar J}(A \cap A_\alpha) \ge \xi$
and $A \backslash A_\alpha,A \cap A_\alpha$ are not in $J_\alpha$.
\ermn
3) We identify $\bar A$ with $\langle (\text{id}_{\bar A \restriction
\alpha},\text{id}_{\bar A \restriction (\alpha +1)}):\alpha < \alpha^* 
\rangle$ (see \scite{2.2}(3) below).  If $\bar J = \langle (J_\alpha,
J_{\alpha +1}):\alpha < \alpha^* \rangle $ is as required in (2) we may write
$\langle J_\alpha:\alpha < \alpha^* \rangle$ instead $\bar J$.  We can replace
$\kappa$ by any other set.  We may write rk$^{(\prime)}(A,\bar A)$ or rk$'(A,
\bar J)$ instead rk$^{(\prime)}_{\bar A}(A)$ or rk$'_{\bar J}(A)$.
\enddefinition
\bigskip

\proclaim{\stag{2.2} Claim}  1) rk$_{\bar A}$,rk$'_{\bar A}$ are well defined
(values: ordinals or $\infty$) and nondecreasing in $A$ (under $\subseteq$).
\nl
2) rk$_{\bar A}(A) \ge rk'_{\bar A}(A)$. \nl
3) rk$'_{\bar A}$ depend just on $\langle \text{id}_{\bar A \restriction
\alpha}:\alpha \le \alpha^* \rangle$ and for $A \subseteq \kappa$, we have
rk$'_{\bar A}(A) = \text{ rk}'_{\bar J}(A)$ where $\bar J = \langle 
(\text{id}_{\bar A \restriction \alpha},\text{id}_{\bar A \restriction
(\alpha +1)}):\alpha \le \alpha^* \rangle$ (so we may write
rk$'_{\langle J_\alpha:\alpha \le \alpha^* \rangle}(A)$ with $J_\alpha =
\text{ id}_{\bar A \restriction \alpha})$. \nl
4) If rk$'_{\bar A}(\kappa) = \zeta$, \ub{then} we can find $\bar Y =
\langle Y_\varepsilon:\varepsilon < \zeta \rangle$, an increasing sequence of
subsets of $\alpha^*,\dbcu_{\varepsilon < \zeta} Y_\varepsilon = \alpha^*$
and $\{A_\alpha:\alpha \in Y_\varepsilon \backslash \dbcu_{\xi < \varepsilon}
Y_\xi\}$ are almost disjoint modulo the ideal generated by $\{A_\alpha:
\alpha \in \dbcu_{\xi < \varepsilon} Y_\xi\}$. \nl
5) If $\bar J = \langle(J_\alpha,J'_\alpha):\alpha < \alpha^* \rangle$ is as
in \scite{2.1}(2), where $J_\alpha,J'_\alpha$ are ideals on $\kappa,
J'_\alpha = J_\alpha + A_\alpha$ and $\bar A = \langle A_\alpha:\alpha <
\alpha^* \rangle$ \ub{then} for every $B \subseteq \kappa$ we have
rk$_{\bar A}(B) \ge \text{ rk}'_{\bar A}(B) \ge \text{rk}'_{\bar J}(B)$.
\endproclaim
\bigskip

\demo{Proof}  Straight: e.g. for the fourth part use $Y_\varepsilon =:
\{\alpha:\text{rk}'_{\bar A}(A_\alpha) \le \varepsilon\}$.
\hfill$\square_{\scite{2.2}}$
\enddemo
\bigskip

\proclaim{\stag{2.3} Claim}  1) If rk$_{\bar A}(\kappa) \ge \kappa^+$
\ub{then} for some $\bar \alpha = \langle \alpha_n:n < \omega \rangle$ we
have $\alpha_n < \alpha_{n+1} < \kappa$ and for every $\ell < k < \omega$
for some $\alpha < \alpha^*$ we have $A_\alpha \cap \{\alpha_\ell,
\alpha_{\ell +1},\dotsc,\alpha_k\} = \{\alpha_{\ell +1},\dotsc,\alpha_k\}$.
\nl
2) If rk$'_{\bar J}(\kappa) \ge \beta$ and $\bar J = \langle(J_\alpha,
J'_\alpha):\alpha < \alpha^* \rangle$ as in \scite{2.1}(2) \ub{then} for
some $\Gamma \subseteq \alpha^*,|\Gamma| \le |\beta|$ we have rk$'_{\bar J
\restriction \Gamma}(\kappa) \ge \beta$. \nl
3) If rk$_{\bar J}(B) \ge \beta,\bar J = \langle(J_\alpha,J'_\alpha):
\alpha < \alpha^* \rangle$ as in \scite{2.1}(2) and $J'_\alpha = J_\alpha
+ A_\alpha$ \ub{then} we can find $\Gamma \subseteq \alpha^*$ such that
\mr
\item "{$(*)$}"  $|\Gamma| \le |\beta| + \aleph_0$ (even $|\Gamma| <
|\beta|^+ + \aleph_0$) and if $A_\alpha \subseteq A'_\alpha \in J'_\alpha$
then rk$_{\langle A'_\alpha:\alpha \in J \rangle}(B) \ge \beta$.
\endroster
\endproclaim
\bigskip

\demo{Proof}  1) Let $\text{rk } = \text{ rk}_{\bar A}$; choose by induction
on $n$ an ordinal $\alpha_n < \kappa$ and for every $\zeta < \kappa^+$ a
decreasing sequence $\langle B^n_{\zeta,0},\dotsc,B^n_{\zeta,n} \rangle$
of sets such that
\mr
\item "{$(\alpha)$}"  $(\forall \ell \le n)(\forall m \le n)[\alpha_\ell
\in B^n_{\zeta,m} \Leftrightarrow \ell > m]$,
\sn
\item "{$(\beta)$}"  rk$(B^n_{\zeta,n}) \ge \zeta$
\sn
\item "{$(\gamma)$}"  each $B^m_{\zeta,n}$ is the intersection of finitely
many $A_\alpha$'s
\ermn
For $n=0$, for every $\zeta < \kappa^+$ there is $\alpha_\zeta < \alpha^*$
such that $\kappa \cap A_{\alpha_\zeta} \ne \kappa$, rk$(\kappa \cap
A_{\alpha_\zeta}) \ge \zeta$, and choose $\alpha^0_\zeta \in \kappa
\backslash A_{\alpha_\zeta}$.  So for some $\alpha_0 < \kappa,\kappa^+ =
\sup\{\zeta < \kappa^+:\alpha^0_\zeta = \alpha_0\}$ and let $B^0_{\zeta,0} =
A_{\alpha_{\xi(\zeta)}}$ where $\xi(\zeta) < \kappa^+$ is the minimal
$\xi > \zeta$ such that $\alpha^0_\xi = \alpha_0$, as in demand $(\beta)$ we
ask ``$\ge \zeta$" not ``$= \zeta$", we succeed.  If we have defined for $n$,
for each $\zeta < \kappa^+$, as rk$(B^n_{\zeta +1,n}) \ge \zeta +1$, there is
$\beta(\zeta,n) < \ell g(\bar A)$ such that $\neg[B^n_{\zeta +1,n} \subseteq
A_{\beta(\zeta,n)}]$ but rk$_{\bar A}(B^n_{\zeta +1,n} \cap
A_{\beta(\zeta,n)}) \ge \zeta$, and choose $\gamma(\zeta,n) \in 
B^n_{\zeta +1,n} \backslash A_{\beta(\zeta,n)}$ so for some $\alpha_{n+1} <
\kappa$, the set $S_n = \{\zeta < \kappa^+:\gamma(\zeta,n) = \alpha_{n+1}\}$
is unbounded in $\kappa^+$.  For every $\zeta < \kappa^+$ let $\xi(\zeta,n) =
\text{ min}\{\xi:\xi \in S_n,\xi > \zeta\}$, let $B^{n+1}_{\zeta,\ell}$ be
$B^n_{\xi(\zeta,n),\ell}$ if $\ell \le n$ and $B^n_{\xi(\zeta,n)+1,n} \cap
A_{\beta(\zeta,n)}$ if $\ell = n+1$.  In the end we know that for $\ell < k
< \omega$, for every $\zeta < \kappa^+$ we have $B^k_{\zeta,\ell} \cap
\{\alpha_\ell,\alpha_{\ell +1},\dotsc,\alpha_k\} = \{\alpha_{\ell +1},
\dotsc,\alpha_k\}$; also $B^k_{\zeta,\ell}$ has the form $\dbca_{m < m(*)}
A_{\alpha(m)}$ for some $\alpha(m) < \ell g(\bar A)$, so for some $m$ we
have $\alpha_\ell \notin A_{\alpha)(m)}$, but $\{\alpha_{\ell +1},\dotsc,
\alpha_k\} \subseteq B^k_{\zeta,\ell} \subseteq A_{\alpha(m)}$ so $\alpha(m)$
is as required in \scite{2.3}(1) for our $\ell < k < \omega$.  Lastly
$\dsize \bigwedge_{n<m} \alpha_n \ne \alpha_m$ hence by Ramsey theorem
\wilog \, $\alpha_n < \alpha_{n+1}$ and we are done. \nl
2) By induction on $\beta$ or by part (3). \nl
3) We can find $\langle (B_\eta,j_\eta):\eta \in ds(\beta) \rangle$ where
$ds(\beta) = \{\eta:\eta$ is a (strictly) decreasing sequence of cardinals
$< \beta\},B_{<>} = B$,rk$'_{\bar J}(B_\eta) \ge \text{ min}(\{\beta\} \cup
\{\eta(\ell):\ell < \ell g(\eta)\}$ and $j_\eta < \alpha^*$ and if $\nu =
\eta \char 94 \langle \gamma \rangle \in ds(\beta)$ then $B_\nu \ne 
\beta_\eta,B_\nu = B_\eta \cap A_{j_\eta} \notin J_{j_\eta}$ and $B_\eta
\backslash B_\nu \notin J_{j_\eta}$.  Let $\Gamma = \{j_\eta:\eta \in ds
(\eta)\}$, now if $A_\alpha \subseteq A'_\alpha \in J_\alpha$ for $\alpha
\in \Gamma$ then we can prove by induction on $\gamma < \beta$ that:
$\eta \in ds(\beta) \Rightarrow \text{ rk}'_{\langle A'_\alpha:\alpha \in
\Gamma \rangle}(B_\eta) \ge \text{max}(\{\beta\} \cup \{\eta(\ell):\ell <
\ell g(\eta)\}$.  \hfill$\square_{\scite{2.3}}$
\enddemo
\bigskip

\definition{\stag{2.4} Definition}  1) For $\bar \lambda = \langle \lambda_n:
n < \omega \rangle$ strictly increasing, let IND$(\bar \lambda)$ mean:
\mr
\item "{$(*)_{\bar \lambda}$}"  for every algebra $M$ with universe
$\dbcu_{n < \omega} \lambda_n$ and $\aleph_0$ functions (all finitary) there
is $\bar \alpha = \langle \alpha_n:n < \omega \rangle$ such that:
{\roster
\itemitem{ $(a)$ }  $\alpha_n < \lambda_n$
\sn
\itemitem{ $(b)$ }  $\alpha_n$ is not in the $M$-closure of

$$
\{\alpha_\ell:\ell \in (n,\omega)\} \cup \{i:\dsize \bigvee_{m<n} i <
\lambda_m\}.
$$
\endroster}
\ermn
2) IND$(\lambda)$ means that $\lambda > \text{ cf}(\lambda) = \aleph_0$ and
for every (equivalently some, see below) $\bar \lambda = \langle \lambda_n:
n < \omega \rangle$ strictly increasing with limit $\lambda$ we have
\mr
\item "{$(*)'_{\bar \lambda}$}"   for every algebra $M$ with universe
$\lambda$ and countably many functions there is $\langle \alpha_n:n \in
\omega \rangle$ such that:
{\roster
\itemitem{ $(a)$ }  $w \subseteq w$ is infinite
\sn
\itemitem{ $(b)$ }  $\alpha_n < \lambda_n$ for $n \in w$
\sn
\itemitem{ $(c)$ }  for $n \in w,\alpha_n$ is not in the $M$-closure of
$$
\{\alpha_\ell:\ell \in w,\ell > n\} \cup \{i:\dsize \bigvee_{\ell \in n
\cap w} i < \lambda_\ell\}
$$
\endroster}
\ermn
3) IND$(\lambda,\kappa) = \text{ IND}^0(\lambda,\kappa)$ means: if $M$ is a
model with universe $\lambda$ and $\kappa$ functions we can find
$\bar \alpha = \langle \alpha_n:n < \omega \rangle$ such that

$$
\alpha_n < \lambda,\alpha_n \notin c \ell_M\{\alpha_\ell:\ell < \omega,
\ell > n\}.
$$
\mn
4) IND$^1(\lambda,\kappa)$ is defined similarly but demanding

$$
\alpha_n \notin c \ell_M\{\alpha_\ell:\all < \omega,\ell \ne n\}.
$$
\enddefinition
\bigskip

\demo{\stag{2.5} Observation}  1) In \scite{2.4}(2), if $(*)'_{\bar \lambda}$
holds for one $\bar \lambda$ with limit $\lambda$, then it holds for every
$\bar \lambda' = \langle \lambda'_n:n < \omega \rangle$ with limit $\lambda$.
\nl
2) If $\lambda$ is uncountable with cofinality $\aleph_0,P$ a forcing notion
of cardinality $\le \mu < \lambda$ or satisfying the $\mu^+$-c.c. for some
$\mu < \lambda$, or $\lambda$-complete \ub{then}: IND$(\lambda) 
\Leftrightarrow \Vdash_P ``\text{IND}(\lambda)"$ and if $\kappa \in [\mu,
\lambda),\mu$ as above \ub{then} IND$(\lambda,\kappa) \Leftrightarrow
\Vdash_P ``\text{IND}(\lambda,\kappa)"$ and if in addition $\mu < \lambda_n
< \lambda_{n+1}$, \ub{then} IND$(\langle \lambda_n:n < \omega \rangle)
\Leftrightarrow \Vdash_P ``\text{ IND}(\langle \lambda_n:n < \omega 
\rangle)"$. \nl
3) IND$(\bar \lambda) \Rightarrow \text{ IND}(\dsize \sum_{n < \omega}
\lambda_n) \Rightarrow \text{ IND}^1(\lambda,\kappa) \Rightarrow
\text{ IND}^0(\lambda,\kappa)$ if $\lambda = \dbcu_{n < \omega} \lambda_n$
and $\lambda_n < \lambda_{n+1}$ and $\lambda_0 > \kappa$.  If $\kappa <
\lambda \le \lambda'$ then

$$
\text{IND}^i(\lambda,\kappa) \Rightarrow \text{ IND}^i(\lambda',\kappa).
$$
\mn
4) If ($i \in \{0,1\}$ and) IND$^i(\lambda,\kappa),\lambda$ minimal for
this $\kappa$ \ub{then}
\mr
\item "{$(a)$}"  $\kappa \le \kappa_1 < \lambda \Rightarrow \text{ IND}^i
(\lambda,\kappa_1)$
\sn
\item "{$(b)$}"  cf$(\lambda) = \aleph_0$ and IND$^1(\lambda,\kappa)$
\ub{or} $\lambda$ is inaccessible
\sn
\item "{$(c)$}"   if $\lambda = \dsize \sum_{n < \omega} \lambda_n$ and
$\lambda_n < \lambda_{n+1}$ \ub{then} not only $(*)'_{\bar \lambda}$ (where
$\bar \lambda = \langle \lambda_n:n < \omega \rangle$) but if $P$ is a c.c.c.
forcing adding a dominating real \ub{then} in $V^P$ for some infinite $w
\subseteq \omega$ we have $(*)_{\bar \lambda \restriction w}$ holds.
\ermn
5) IND$^1(\lambda,\kappa)$ is equivalent to IND$^0(\lambda,\kappa)$.
\enddemo
\bigskip

\demo{Proof}  1, 2), 3)  Check. \nl
4) \ub{Clause (a) of (a)}:  Assume not and first let $i=0$.  Let $\chi = 
\beth_3(\lambda)^+$ and let $M$ be the model with universe $\lambda$ and the
functions ($n$-place from $\lambda$ to $\lambda$ for some $n$) definable in
$({\Cal H}(\chi),\in,<^*_\chi)$ with the parameters $\lambda,\kappa,\kappa_1$.
Clearly $(M,\alpha)_{\alpha < \kappa_1}$ exemplifies $\neg \text{IND}^i(
\lambda,\kappa_1)$.  Let $F^-_n,F^+_n$ be such that for $\bar \beta = \langle
\beta_\ell:\ell < n \rangle,\beta_\ell < \lambda$ we have:
$F^+_n(-,\bar \beta)$ is a one-to-one function from $c \ell_M(\{\beta_0,\dotsc,\beta_{n-1}\} \cup \kappa_1)$ onto $\kappa_1$ and $F^-_n(-,\beta)$ is its
inverse.  We can apply the assumption IND$^i(\lambda,\kappa)$ to the model
$(M,F^+_n,F^-_n,\beta)_{n < \omega,\beta < \kappa}$, so there are $\alpha_n
(n < \omega)$ as in \scite{2.4}(3).  By the assumption for no infinite $w
\subseteq \omega$ is $\{\alpha_n:n \in w\}$ as required in \scite{2.4}(3) for
$(M,\beta)_{\beta < \kappa_1}$ hence for some infinite $w \subseteq \omega,
\dsize \bigwedge_{n \in w} \alpha_n \in c \ell_M(\{ \alpha_\ell:n < \ell \in
w\} \cup \kappa_1)$, just choose by induction on $\ell < \omega,u_\ell,
n_\ell$ such that: $\dsize \bigwedge_{m < \ell} n_m < n_\ell < \omega,u_\ell
\subseteq (n_\ell,\omega)$ is infinite, $u_{\ell +1} \subseteq u_\ell
\subseteq w$ and $\alpha_{n_\ell} \notin c \ell_M(\{\alpha_n:n \in u_\ell\})$,
we cannot succeed so $w$ or some $u_\ell$ is as required.  By renaming $w =
\omega$, so for every $n$ for some $k_n \in (n+1,\omega)$ we have $\alpha_n
\in c \ell)M(\{\alpha_{n+1},\dotsc,\alpha_{k_n}\} \cup \kappa_1)$; as we can
increase $k_n$, \wilog \, $k_n < k_{n+1}$ hence $m \le n \Rightarrow \alpha_m
\in c \ell_M(\{\alpha_{n+1},\dotsc,\alpha_{k_n}\} \cup \kappa_1)$ (just prove
this by induction on $n$), so $\gamma_{n,m} =: F^+_{k_n-n}(\alpha_m,
\alpha_{n+1},\dotsc,\alpha_{k_n}) < \kappa_1$ and for each $n$ we have:
$\langle \gamma_{n,m}:m \le n \rangle$ is with no repetitions.  Choose by
induction on $\ell,m_\ell \in [\ell^2,(\ell +1)^2)$ such that 
$\gamma_{(\ell +1)^2,m_\ell} \notin \{\gamma_{(q+1)^2,m_q}:q < \ell\}$.  But
as $\neg \text{ IND}^i(\kappa_1,\kappa)$ (because $\lambda > \kappa_1$ was
minimal such that ...) for some $\ell < p < \omega$ we have 
$\gamma_{(\ell+1)^2,m_\ell} \in c \ell_M(\{\gamma_{(q+1)^2,m_q}:\ell < q <
p\} \cup \kappa)$ and using some $F^-_{p^2-(\ell +1)^2}$ we have
$\alpha_{m_\ell} \in c \ell_M(\{\alpha_q:q \text{ is } \ge(\ell +1)^2
\text{ but } \le k_{(p^2)}\} \cup \kappa)$. \nl
[Why?  First note that $\gamma_{(q+1)^2,m_q}$ belong to this model for
$q=\ell +1,\dotsc,p-1$, (using $F_{k_{(q+1)^2}}$) hence also 
$\gamma_{(\ell +1)^2,m_\ell}$ belongs to this model by the choice of $\ell$
and $p$; a contradiction.]

If $i=1$ the proof is similar: choose, by induction on $\ell,k_\ell,m_\ell,
m_\ell$ such that $k_\ell < m_\ell < k_{\ell +1}$ and $\alpha_{m_\ell} \in c
\ell_M(\{\alpha_n:n \in [k_\ell,m_\ell)$ or $n \in (m_\ell,k_{\ell +1})\} \cup
\kappa_1)$, this is possible as otherwise $\{\alpha_n:n \in [k_\ell,\omega)\}$
contradict ``$M$ exemplifies $\neg \text{IND}^1(\lambda,\kappa_1)$". Let

$$
\gamma_\ell = F^+_{k_{\ell +1}-k_\ell-1}(\alpha_{m_\ell};\alpha_{k_\ell},
\alpha_{k_\ell +1},\dotsc,\alpha_{m_\ell-1},\alpha_{m_\ell +1},\dotsc,
\alpha_{k_{\ell +1}-1}) < \kappa_1.
$$
\mn
For some $\ell < \ell(*) < \omega$ we have $\gamma_\ell \in c \ell_M(\{
\gamma_0,\dotsc,\gamma_{\ell-1},\gamma_{\ell +1},\dotsc,\gamma_{\ell(*)-1}\}
\cup \kappa)$ (because $\neg \text{ IND}^1(\kappa_1,\kappa)$ as $\lambda$ is
first and the choice of $M$), hence $\alpha_\ell \in c \ell_M(\{\alpha_n:n <
\omega,n \le \ell\} \cup \kappa)$; a contradiction.
\mn
\ub{Clause (b) of (4)}:  By the definition easily $\aleph_0 < \text{ cf}
(\lambda) \le \kappa$ is impossible. \nl
[Why?  Let $\lambda = \dsize \sum_{i < \kappa} \lambda_i$, with $\lambda_i
< \lambda$, and by the minimality of $\lambda$ let $M_i$ be a model with
universe $\lambda_i$ and $\le \kappa$ functions exemplifying $\neg \text{IND}
(\lambda_i,\kappa)$ and lastly let $M$ be the model with universe $\lambda$ 
and the functions of all the $M_i$; check that $M$ exemplifies $\neg
\text{ IND}(\lambda,\kappa)$.] \nl
By \scite{2.5}(4) clause (a) it follows that

$$
[\text{cf}(\lambda) > \aleph_0 \and \kappa_1 < \lambda \Rightarrow \kappa_1
< \text{ cf}(\lambda)].
$$
\mn
So if cf$(\lambda) > \aleph_0$ then $\lambda$ is regular, it is inaccessible
as it is not a successor as trivially

$$
\neg \text{IND}^i(\mu,\kappa) \Rightarrow \neg\text{IND}(\mu^+,\kappa).
$$
\mn
We still have to prove $\text{IND}^1(\lambda)$ when cf$(\lambda) = \aleph_0$;
if $i=1$ this is trivial so assume $i=0$.  So assume cf$(\lambda) =
\aleph_0,\bar \lambda = \langle \lambda_n:n < \omega \rangle,\lambda_n <
\lambda_{n+1}$ and $\lambda = \dsize \sum_{n < \omega} \lambda_n$.  We should
prove IND$^1(\lambda,\kappa)$, so let us be given $M$ as in \scite{2.4}(4).
Without loss of generality, expanding $M$ by examples to $\neg \text{IND}
(\lambda_n,\kappa)$, we get $M^+$. \nl
As IND$^0(\lambda,\kappa)$ there are $\alpha_n < \lambda$ (for $n < \omega$)
such that $\alpha_n \notin c \ell_{M^+}\{\alpha_\ell:\ell < \omega,\ell > 
n\}$.  Without loss of generality $\langle \alpha_\ell:\ell < \omega 
\rangle$ is
strictly increasing, and letting $m_n = 
\text{ min}\{m:\lambda_m > \alpha_n\}$,
\wilog \, $\langle m_n:n < \omega \rangle$ is constant or strictly increasing.
If $\dsize \bigwedge_n M-n = m(*)$ we get a contradiction to ``$M^+$ expand a
counterexample to IND$^0(\lambda_{m(*)},\kappa)$", so $\dsize \bigwedge_n
m_n < m_{n+1}$.  It suffices to prove that for every $n_1 < \omega$ there is
$n_2 \in (n_1,\omega)$ such that $\alpha_{n_2} \notin c \ell_M(\{\alpha_\ell:
\ell < \omega,\ell > n_2\} \cup \lambda_{n_1})$ as then we can produce an
example for $(*)'_{\bar \lambda}$.  But if for some $n_1$ there is no $n_2$
as above the proof of clause (a) gives a contradiction.  
\mn
\ub{Clause (c) of (4)}:  Left to the reader. \nl
5) By \scite{2.5}(3),

$$
\text{IND}^1(\lambda,\kappa) \Rightarrow \text{IND}^0(\lambda,\kappa)
\text{ and } \lambda \le \lambda' \and \text{ IND}^i(\lambda,\kappa)
\Rightarrow \text{ IND}^i(\lambda,\kappa).
$$
\mn
Hence it suffices to prov: if $\lambda$ is minimal such that $\text{IND}^0
(\lambda,\kappa)$ then $\text{IND}^1(\lambda,\kappa)$.  let $M$ be a model
with universe $\lambda$ and vocabulary of cardinality $\le \kappa$ and
we shall prove the conclusion of \scite{2.4}(4) (= the Definition of
IND$^1(\lambda,\kappa))$.  Let for $n < \omega,F^+_n,F^-_n$ be $(n+2)$-place
functions from $\lambda$ to $\lambda$ such that
\mr
\item "{$(*)$}"  if $\gamma < \lambda$ and $\bar \beta \in {}^n \lambda$ then
$F^+_n(-,\bar \beta,\gamma)$ is a one-to-one function from $c \ell_M(\{
\beta_\ell:\ell < \ell g(\bar \beta)\} \cup \{i:i \le \gamma \text{ or }
i < \kappa\})$ onto $|\kappa + \gamma|$ and $F^-_n(-,\bar \beta,\gamma)$ be
the inverse function.
\ermn
Let $M^* = (M,F^+_n,F^-_n)_{n < \omega},M^+ = (M^*,i)_{i < \kappa}$ but as
IND$(\lambda,\kappa)$ we can apply Definition \scite{2.4}(4) and get
$\langle \alpha_n:n < \omega \rangle,\alpha_n < \lambda,\alpha_n \notin
c \ell_{M^+}(\{\alpha_\ell:\ell \in (n,\omega)\})$.  Without loss of 
generality $\langle \alpha_n:n < \omega \rangle$ is strictly increasing;
$\alpha_n > \kappa$ (as without loss of generality each 
$i \le \kappa$ is an individual
constant of $M$), clearly it suffices to prove
\mr
\item "{$(**)$}"   for any $n < \omega$ for some $m \in (n,\omega)$, \nl
$\alpha_m \notin c \ell_M(\{\alpha_\ell:\ell < n \text{ or } \ell > m\})$
hence it suffices to prove
\sn
\item "{$(**)'$}"  for any $n < \omega$ for some $m \in (n,\omega)$,
$$
\alpha_m \notin c \ell_M(\{\alpha_\ell:\ell > m\} \cup \{i:i \le \alpha_m\}).
$$
\ermn
But here we can repeat the proof of clause (a) of \scite{2.4}(4) for the case
$i=0$.  \hfill$\square_{\scite{2.5}}$
\enddemo
\bigskip

\proclaim{\stag{2.6} Claim}  1) Assume $\lambda > \text{ cf}(\lambda),
|{\frak a}| < \min({\frak a})$ and sup$({\frak a}) = \lambda$ and \nl
rk$'_{\langle J_{< \theta}[{\frak a}]:\theta \in \text{pcf}({\frak a})
\rangle}({\frak a}) \ge |{\frak a}|^+$.  \ub{Then} IND$(\lambda,|{\frak a}|)$.
\nl
2) Moreover, for any model $M$ with universe $\lambda$ and $|{\frak a}|$
functions and ${\frak c}$ such that ${\frak a} \subseteq {\frak c} \subseteq
\text{ pcf}({\frak a}),|{\frak c}| < \min({\frak a})$ and $\langle
{\frak b}_\theta[{\frak a}]:\theta \in \text{ pcf}({\frak a}) \rangle$ a
generating sequence we can find $\bar \alpha = \langle \alpha_\theta:\theta
{\frak a} \rangle \in \Pi {\frak a}$ such that, defining for ${\frak b}
\subseteq {\frak a}$:

$$
c \ell_{M,\bar \alpha}({\frak b}) = {\frak b} \cup \{\theta \in {\frak a}:
\alpha_\lambda \in c \ell_M(\{\alpha_\mu:\mu \in {\frak b}\})\};
$$
\mn
we have
\mr
\item "{$\otimes_1$}"  $[\theta \in {\frak c} 
\Rightarrow c \ell_{M,\bar \alpha}({\frak b}_\theta[{\frak a}]) \in
J_{\le \theta}[{\frak a}]]$;
\sn
\item "{$\otimes_2$}"  $c \ell_{M,\bar \alpha}(-)$ is a closure operation
on ${\frak a}$, i.e.

$$
{\frak b}_1 \subseteq {\frak b}_2 \Rightarrow c \ell_{M,\bar \alpha}
({\frak b}_1) \subseteq c \ell_{M,\bar \alpha}({\frak b}_2),
$$

$$
{\frak b} \subseteq c \ell_{M,\alpha}({\frak b}),
$$

$$
c \ell_{M,\bar \alpha}(c \ell_{M,\bar \alpha}({\frak b})) = c \ell
_{M,\bar \alpha}({\frak b}).
$$
\endroster
\endproclaim
\bigskip

\remark{Remark}  See \scite{3.14} - \scite{3.17} for more.
\endremark
\bigskip

\demo{Proof}  1) Let us define $\bar J = \langle (J_{< \theta}[{\frak a}],
J_{\le \theta}[{\frak a}]):\theta \in \text{ pcf}({\frak a}) \rangle$.  
We prove part (1) assuming part (2).  Choose ${\frak c} \subseteq 
\text{ pcf}({\frak a}),|{\frak c}| = |{\frak a}|^+$ such 
that rk$'_{\bar J \restriction {\frak c}}({\frak a}) \ge |{\frak a}|^+$ 
(this is possible by \scite{2.3}(2) and \wilog
\, ${\frak a} \subseteq {\frak c}$ and let $\langle {\frak b}_\theta
[{\frak a}]:\theta \in \text{ pcf}({\frak a}) \rangle$ be a generating 
sequence for ${\frak a}$ (exists by \cite[2.6]{Sh:371}).  For proving 
IND$(\lambda,|{\frak a}|)$ let $M$ be a model with universe $\lambda$ 
and $\le |{\frak a}|$
functions, by part (2) there is a sequence $\bar \lambda = \langle
\alpha_\tau:\tau \in {\frak c} \rangle$ as there.  Let for $\theta \in
{\frak c},{\frak d}_\theta =: c \ell_{M,\bar \alpha}({\frak b}_\theta
[{\frak a}]) \subseteq {\frak a}$ (as defined in part (2)) so 
${\frak b}_\theta[{\frak a}] \subseteq {\frak d}_\theta \in 
J_{\le \theta}[{\frak a}]$ hence $J_{< \theta}[{\frak a}] + 
{\frak d}_\theta = J_{\le \theta}[{\frak a}]$.  So
by \scite{2.2}(5) we know rk$_{\langle {\frak d}_\theta:\theta \in {\frak c}
\rangle}({\frak a}) \ge \text{ rk}'_{\langle {\frak d}_\theta:\theta \in
{\frak c} \rangle}({\frak a}) \ge \text{ rk}'_{\bar J \restriction {\frak c}}
({\frak a}) \ge \kappa^+$ (by the choice of ${\frak c}$ above).  Now by
\scite{2.3}(1) we can find $\tau_n \in {\frak a}$ for $n < \omega$, pairwise
distinct and strictly increasing with $n$ such that for every $n < m < \omega$
for soe $\theta_{n,m} \in {\frak c}$ we have $\{\tau_n,\tau_{n+1},\dotsc,
\tau_m\} \cap {\frak d}_{\theta_{n,m}} = \{\tau_{n+1},\dotsc,\tau_m\}$,  note:
as $\tau_m \in {\frak d}_{\theta_{n,m}}$ necessarily $\theta_{n,m} \ge 
\tau_m$.
So by the choice of $\bar \alpha$, we have $\alpha_{\tau_n} \notin c \ell_M
(\{\alpha_{\tau_{n+1}},\alpha_{\tau_{n+2}},\dotsc,\alpha_{\tau_m}\})$.  So
$\langle \alpha_{\tau_n}:n < \omega \rangle$ are as required in the definition
of IND$(\lambda,|{\frak a}|$). \nl
2) Let $\bar{\frak b} = \langle {\frak b}_\theta[{\frak a}]:\theta \in
\text{ pcf}({\frak a}) \rangle$ be the generating sequence for ${\frak a}$;
\wilog \, max pcf$({\frak a}) \in {\frak c}$ and $\theta \in {\frak c} 
\backslash \{\text{max pcf}({\frak a})\} \Rightarrow \min(\text{pcf}
({\frak a}) \backslash \theta^+) \in {\frak c}$.  We know that
\mr
\item "{$(*)_{\frak a}$}"  there is $\bar F = \langle F_\theta:\theta \in
\text{ pcf}({\frak a}) \rangle$ such that:
{\roster
\itemitem{ $(a)$ }  $F_\theta \subseteq \Pi {\frak a}$ and $|F_\theta| \le
\theta$
\sn
\itemitem{ $(b)$ }  $\{f \restriction {\frak b}:f \in F_\theta\}$ is cofinal
in $\Pi {\frak b}$ for every ${\frak b} \in J_{\le \theta}[{\frak a}]$
\sn
\itemitem{ $(c)$ }  $F_\theta$ includes $\cup \{F_\tau:\tau \in \theta \cap
\text{ pcf}({\frak a})\}$ and is closed under some natural operations
\sn
\itemitem{ $(d)$ }  if $\tau \in \theta \cap {\frak c}$ then $f \in F_\theta
\Rightarrow (\exists g \in F_\tau)(f \restriction {\frak b}_\tau[{\frak a}]
= g \restriction {\frak b}_\tau[{\frak a}])$.
\endroster}
\ermn
[Why?  E.g. the proof of \cite{Sh:355},3.5.]

Let $M$ be a model with universe $\lambda$ and vocabulary of cardinality
$|{\frak a}|$.  For every $f \in \Pi {\frak a}$ (e.g., $f \in F_\theta,\theta
\in \text{ pcf}({\frak a}))$ 
we define $g_f \in \Pi {\frak a}$ by $g_f(\tau) =:
\sup[\tau \cap c \ell_M(\text{Rang }f)]$.  For every $\theta \in {\frak c}
\cap \text{ pcf}({\frak a}),\{g_f:f \in \dbcu_{\tau \in \theta \cap {\frak c}}
F_\tau\}$ is a subset of $\Pi {\frak a}$ of cardinality $< \theta$ (here
instead $|{\frak c}| < \min({\frak a})$, just $\theta \in {\frak c} 
\Rightarrow |{\frak c} \cap \theta| < \theta$ suffice), so there is
$g^\theta \in F_\theta$ such that:
\mr
\item "{$(*)_1$}"  $f \in \dbcu_{\tau \in \theta \cap {\frak c}} F_\tau
\Rightarrow g_f < g^\theta \text{ mod } J_{< \theta}[{\frak a}]$.
\ermn
Define $g^* \in \Pi {\frak a}$ by $g^*(\tau) = \sup\{g^\theta(\tau)+1:\theta
\in {\frak c}\}$ (remember $|{\frak c}| < \min({\frak a})$).  So there is
$h \in F_{\text{max pcf}({\frak a})}$ such that $g^* < h$ (see
$(*)_{\frak a}(b))$; we shall show that:
\mr
\item "{$(*)_2$}"  for any such $h$ the sequence $\langle h(\tau):\tau \in
{\frak a} \rangle$ is as required.
\endroster
\enddemo
\bigskip

\demo{Proof of $(*)_2$}  So let $\alpha_\tau = h(\tau)$.  Note that
$\otimes_2$ from \scite{2.6}(2) is trivial so we shall prove $\otimes_1$.
Assume $\theta \in {\frak c}$ and let ${\frak b} = {\frak b}_\theta
[{\frak a}] \in J_{\le \theta}[{\frak a}]$ so by clause (d) of $(*)_{\frak a}$
for some $f_1 \in F_\theta$ we have 
\mr
\item "{$\oplus_1$}"  $h \restriction {\frak b} = f_1 \restriction {\frak b}$;
\ermn
we can assume $\theta < \text{ max pcf}({\frak a})$ (otherwise conclusion
is trivial) and let $\sigma = \min(\text{pcf}({\frak a}) \backslash 
\theta^+)$, by an assumption made in the beginning of the proof $\sigma 
\in {\frak c}$ and so as $f_1 \in F_\theta$ by the choice of $g^\sigma$ we
have:

$$
g_{f_1} < g^\sigma \text{ mod } J_{< \sigma}[{\frak a}]
$$
\mn
but by the choice of $g^*$

$$
g^\sigma < g^*
$$
\mn
and by the demand of $h$

$$
g^* < h
$$
\mn
together

$$
g_{f_1} < h \text{ mod } J_{< \sigma}[{\frak a}]
$$
\mn
so for some ${\frak d} \in J_{< \sigma}[{\frak a}] = J_{\le \theta}
[{\frak a}]$ we have:
\mr
\item "{$\oplus_2$}"  $g_{f_1} \restriction ({\frak a} \backslash {\frak d})
< h \restriction ({\frak a} \backslash {\frak d})$.
\ermn
Now for any $\tau \in {\frak b}$
\mr
\item "{$\oplus_3$}"  $\tau \in c \ell_{M,\bar \alpha}({\frak b})
\Rightarrow \alpha_\tau \in c \ell_M[\{\alpha_\kappa:\kappa \in {\frak b}\}]
\Rightarrow \alpha_\tau \in c \ell_M[\text{Rang}(h \restriction {\frak b})]
\Rightarrow \alpha_\tau \in c \ell_M[(\text{Rang}(f_1 \restriction 
{\frak b})] \Rightarrow \alpha_\tau \in c \ell_M(\text{Rang}(f_1))
\Rightarrow \alpha_\tau \le g_{f_1}(\tau)$.
\ermn
By $\oplus_2 + \oplus_3$ the required conclusion follows.
\hfill$\square_{\scite{2.6}}$
\enddemo
\bn
E.g.
\proclaim{\stag{2.7} Claim}  If pp$(\aleph_\omega) > \aleph_{\omega_1}$
\ub{then} IND$(\aleph_\omega)$.
\endproclaim
\bigskip

\demo{Proof}  Let pp$(\aleph_\omega) = \aleph_{\alpha^*}$ (so $\alpha^* <
\omega_4,\alpha^* = \beta^* +1$, see \cite[Ch.IX,2.1]{Sh:g}) let ${\frak a} =
\{\aleph_{i+1}:5 \le i < \omega\}$ so we know pcf$({\frak a}) = \{
\aleph_{i+1}:5 \le i < \alpha^*\}$ and let $\bar {\frak b} = \langle 
{\frak b}_\lambda:\lambda \in \text{ pcf}({\frak a}) \rangle$ be a normal
generating sequence for pcf$({\frak a})$ (not ${\frak a}$!); \wilog \,
${\frak b}_{\aleph_{\alpha^*}} = \text{ pcf}({\frak a})$.  Let
$\bar {\frak c} = \langle {\frak b}_\lambda \cap {\frak a}:\lambda \in
\text{ pcf}({\frak a}) \rangle$.  Now by localization (\cite[Ch.VIII]{Sh:g},
2.6) we know that for some club $E$ of $\omega_1:\delta \in E \Rightarrow
\text{pp}(\aleph_\delta) > \aleph_{\omega_1}$ (so we can assume $\omega \in
E$).

Let $E = \{\beta_\zeta:\zeta < \omega_1\}$ (increasing in $\zeta$).  Hence by
\cite[Ch.II]{Sh:g}, 1.5A we have
\mr
\item "{$(*)_0$}"  for limit $\delta < \omega_1,\delta < \beta < \omega_1,
\delta \in E$, for some unbounded ${\frak d} \in \aleph_\delta \cap
\text{ Reg}$, we have $\aleph_{\beta +1} \in \text{ pcf}_{J^{\text{bd}}
_{\frak d}}({\frak d})$.
\ermn
Hence we can prove by induction on $\varepsilon < \omega_1$ that:
\mr
\item "{$(*)$}"  if $\beta_\varepsilon \le \zeta < \omega_1,{\frak b}
\subseteq {\frak a}$ and ${\frak b} = {\frak a} \cap 
{\frak b}_{\aleph_{\zeta +1}}$ mod $J_{< \aleph_{\zeta +1}}[{\frak a}]$
\ub{then} rk$'_{\bar {\frak c}}({\frak b}) \ge \varepsilon$.
\ermn
[Why?  For $\varepsilon = 0$ this is trivial and also for $\varepsilon$   
limit.  If $\varepsilon = \xi + 1$ we can find ${\frak d} \in \aleph
_{\beta_\varepsilon} \cap \text{ Reg} \backslash \aleph_{\beta_\xi}$ such
that $\aleph_{\zeta +1} \in \text{ pcf}_{J^{\text{bd}}_{\frak d}}$ so \wilog
\, $\aleph_{\zeta +1} = \text{max pcf}({\frak d})$.  By \cite[Ch.I]{Sh:g},
\sciteu{1.12} we have that the set ${\frak d}'$ of $\theta \in {\frak d}$
such that ${\frak b}_\theta \cap {\frak a} \subseteq {\frak b} \text{ mod }
J_{< \theta}[{\frak a}]$ is $= {\frak d}$ mod $J_{< \aleph_{\zeta +1}}
[{\frak d}]$ hence is not bounded in ${\frak d}$, hence is not  bounded in
$\aleph_{\beta_\varepsilon}$.  But ${\frak d}' \subseteq {\frak d} \subseteq
\aleph_{\beta_\varepsilon} \cap \text{ Reg} \backslash \aleph_{\beta_
\varepsilon}$, hence by the induction hypothesis $\theta \in {\frak d}'
\Rightarrow \text{ rk}'_{\bar {\frak c}}({\frak b}_\theta \cap {\frak b}) \ge
\xi$, but of course ${\frak b} \backslash {\frak b}_\theta \ne \emptyset$.]

Now apply \scite{2.6}. \hfill$\square_{\scite{2.7}}$ 
\enddemo
\bigskip

\proclaim{\stag{2.8} Proclaim}  If $|{\frak a}| < \text{ min}({\frak a}),
\lambda = \sup({\frak a})$ is singular and pcf$_{J^{\text{bd}}_{\frak a}}
({\frak a})$ contains an interval of \text{\rm Reg\/} of cardinality
$|{\frak a}|^+$ \ub{then} IND$(\lambda)$.
\endproclaim
\bigskip

\demo{Proof}  Similar to the proof of \scite{2.7}.
\enddemo
\bigskip

\demo{\stag{2.9} Discussion}  We can also prove e.g.: if $\lambda =
\text{ tcf}(\dsize \prod_{\varepsilon < \kappa} \lambda_\varepsilon/[\kappa]
^{< \aleph_0})$, satisfies $\lambda > \lambda_\varepsilon =
\text{ cf}(\lambda_\varepsilon) > \kappa$, then for every algebra $M$ on
$\dsize \sum_{\varepsilon < \kappa} \lambda_\varepsilon$ with $< \min\{\lambda
_\varepsilon:\varepsilon < \kappa\}$ functions there are $\alpha_\varepsilon
< \lambda_\varepsilon (\varepsilon < \kappa)$ such that: for finite $u
\subseteq \kappa$ we have $\{\zeta:\alpha_\zeta \in c \ell_M(\{\alpha
_\varepsilon:\varepsilon \in u\})\}$ is finite (and more).  Not clear how
interesting is this statement and where it leads.
\enddemo
\bigskip

\proclaim{\stag{2.10} Claim}  Assume ${\frak a} \subseteq {\frak c} \subseteq
\text{ pcf}({\frak a})$ and $|{\frak c}| < \min({\frak a})$ (or ${\frak c} \in
J_*[\text{pcf}({\frak a})]$, see \cite[Ch.VIII,\S3]{Sh:g} and \cite{Sh:E11}),
so pcf$({\frak a}) = \text{ pcf}({\frak c})$.  Then

$$
\text{rk}'_{\langle J_{< \theta}[{\frak a}]:\theta \in \text{pcf}({\frak a})
\rangle}({\frak a}) =
\text{rk}'_{\langle J_{< \theta}[{\frak c}]:\theta \in \text{pcf}({\frak a}')
\rangle}({\frak c}).
$$
\endproclaim
\bigskip

\demo{Proof}  Let $\bar{\frak b} = \langle{\frak b}_\theta[{\frak a}]:\theta
\in \text{ pcf}({\frak a}) \rangle$ be a generating sequence for ${\frak a}$,
hence we know that letting ${\frak b}_\theta[{\frak c}] =: {\frak c} \cap
\text{ pcf}({\frak b}_\theta[{\frak a}])$, we have: $\bar {\frak b}' = \langle
{\frak b}_\theta[{\frak c}]:\theta \in \text{ pcf}({\frak a}) \rangle$ is a
generating sequence for ${\frak c}$.  Without loss of generality ${\frak b}
_\theta[{\frak c}] \cap {\frak a} = {\frak b}_\theta[{\frak a}]$ and
${\frak b}_{\text{max pcf}({\frak a})}[{\frak a}] = {\frak a}$ hence
${\frak b}_{\text{max pcf}({\frak a})}[{\frak c}] = {\frak a}$.  So we can
prove easily:

$$
{\frak d} \subseteq {\frak a} \Rightarrow
\text{ rk}'_{\langle J_{< \theta}[{\frak a}]:\theta \in \text{pcf}({\frak a})
\rangle}({\frak d}) \le \text{ rk}'_{\langle J_{< \theta}[{\frak c}]:
\theta \in \text{pcf}({\frak a})\rangle}({\frak d})
$$
\mn
(as 
$J_{< \theta}[{\frak a}] = J_{< \theta}[{\frak c}] \restriction {\frak a}$).

For the other direction we prove
\mr
\item "{$(*)$}"  if $n < \omega,\{\theta_1,\dotsc,\theta_n\} \subseteq
\text{ pcf}({\frak a})$ then
\endroster

$$
\text{rk}'_{\langle J_{< \theta}[{\frak c}]:\theta \in \text{pcf}({\frak a})
\rangle}( \dbca^n_{\ell =1} {\frak b}_{\theta_\ell}[{\frak c}]) \le
\text{ rk}'_{\langle J_{< \theta}[{\frak c}]:\theta \in \text{pcf}({\frak a})
\rangle}(\dbca^n_{\ell = 1} {\frak b}_{\theta_\ell}[{\frak a}]).
$$

${{}}$ \hfill$\square_{\scite{2.10}}$
\enddemo
\newpage

\head {\S3 existence of free sets implies restrictions on pcf}
\endhead  \resetall 
\bigskip

\definition{\stag{3.1} Definition}  Suppose $\bar{\Cal I} = \langle (\kappa_n,
I_n):n < \omega \rangle$ is such that: $I_n$ is an ideal on $\kappa_n$. \nl
1) We define $J_n = J^{\bar{\Cal I}}_n$ an ideal on $\dsize \prod_{\ell <n}
\kappa_\ell$:
\mr
\item "{$(*)_0$}"  $J_0 =$ the empty ideal on $\{<>\}$
\sn
\item "{$(*)_1$}"  $J_{n+1} = \{A \subseteq \dsize \prod_{\ell \le n}
\kappa_\ell:\{\alpha < \kappa_n:\{\eta \in \dsize \prod_{\ell < n} 
\kappa_\ell:\eta \char 94 \langle \alpha \rangle \in A\} \notin J_n\} \in
I_n\}$ 
\ermn
we let 
$\bar J^{\bar{\Cal I}} = \langle J^{\bar{\Cal I}}_n:n < \omega \rangle$. \nl
2) We say $\langle J_n:n < \omega \rangle$ is a candidate (for
$\bar{\Cal I}$) if in $(*)_1$ we weaken ``$J_{n+1} = \ldots$" to
``$J_n \subseteq \ldots$".  [So there may be many candidates.]
\enddefinition
\bigskip

\demo{\stag{3.2} Fact}  1) In definition \scite{3.1}(1) above, each $J_n$ is 
an ideal on $\dsize \prod_{\ell < n} \kappa_\ell$. \nl
2) If $I_0,\dotsc,I_{n-1}$ are $\sigma$-complete \ub{then} so is $J_n$. \nl
3) If $\langle J'_n:n < \omega \rangle$ is a candidate for $\bar{\Cal I}$,
then $\dsize \bigwedge_{n < \omega} J'_n \subseteq J^{\bar{\Cal I}}_n$.
\enddemo
\bigskip

\proclaim{\stag{3.3} Claim}  For $\bar{\Cal I},\bar J^{\bar{\Cal I}}$ as in
\scite{3.1}(1), $A \in J_n$ \ub{iff} for some functions $f_0,\dotsc,f_{n-1}$
we have:

$$
\text{Dom}(f_\ell) = \dsize \prod^{n-1}_{m= \ell +1} \kappa_m \text{ and Rang}
(f_\ell) \subseteq I_\ell,
$$
\mn
and $A \subseteq \dbcu_{\ell < n} A^n_\ell(f_\ell)$ where

$$
A^n_\ell(f_\ell) = \{\eta \in \dsize \prod_{m<n} \kappa_m:\eta(m) \in f_\ell
(\eta \restriction (m,n))\}.
$$
\endproclaim
\bigskip

\demo{Proof}  By induction on $n$.
\enddemo
\bigskip

\proclaim{\stag{3.4} Theorem}  Let $\bar{\Cal I} = \langle (\kappa_n,I_n):
n < \omega \rangle$ be as in \scite{3.1}, $\kappa = \dsize \sum_{n < \omega}
\kappa_n \le \mu_0 = \text{cf}(\mu_0) < \mu_1 < \mu = \text{ cf}(\mu)$.
\ub{Then} $\otimes_1 \Rightarrow \otimes_2$ where
\mr
\item "{$\otimes_1$}"  for each $n$ there is $\langle \lambda^n_i:i <
\kappa_n \rangle$ such that: \nl
$\lambda^n_i \in [\mu_0,\mu_1) \cap \text{ Reg}$ and 
$\dsize \prod_{i < \kappa_n} \lambda^n_i/I_n$ is $\mu$-directed
\sn
\item "{$\otimes_2$}"  for some $\bar{\Cal I}$-candidate, $\bar J = \langle
J_n:n < \omega \rangle$ (except for clause $(\delta),J_n = J^{\bar{\Cal I}}_n$
is O.K.) we have:
\sn
\item "{$\otimes^2_{\bar J}$}"  there are $\bar \lambda^n = \langle
\lambda_\eta:\eta \in \dsize \prod_{\ell < n} \kappa_\ell \rangle$ for
$n \in (0,\omega)$ such that for each $n$:
{\roster
\itemitem{ $(\alpha)$ }  $(\dsize \prod_{\eta \in \dsize \prod_{\ell < n}
\kappa_\ell} \lambda_\eta)/J^{\bar{\Cal I}}$ has true cofinality $\mu$
\sn
\itemitem{ $(\beta)$ }  $\mu_0 \le \lambda_\eta = \text{ cf}(\lambda_\eta)
< \mu_1$ (note that by clause $(\alpha)$ we have \nl
$\{\eta \in \dsize \prod_{\ell < n} \kappa_\ell:\lambda_\eta = 
\mu_0\} \in J_n$)
\sn
\itemitem{ $(\gamma)$ }  if $0 < n < \omega,\alpha < \kappa_n$ and $\eta \in
\dsize \prod_{\ell < n} \kappa_\ell$ \ub{then} $\lambda_\eta > \mu_0
\Rightarrow \lambda_\eta > \lambda_{\eta \char 94 \langle \alpha \rangle}$ so
$\{\eta \in \dsize \prod_{\ell < n} \kappa_\ell:\lambda_{\eta \char 94
\langle \alpha \rangle} \nleq \lambda_\eta\} \in J_n$ hence
\sn
\itemitem{ $(\gamma)'$ }  $\{\eta \in \dsize \prod_{\ell \le n} \kappa_\ell:
\lambda_\eta \nleq \lambda_{\eta \restriction n}\} \in J_{n+1}$\sn
\sn
\itemitem{ $(\delta)$ }  $J_n = \{A \subseteq \dsize \prod_{\ell < n}
\kappa_\ell:\text{max pcf}\{\lambda_\eta:\eta \in A\} < \mu\}$.
\endroster}
\endroster
\endproclaim
\bn
\ub{Question}:  Can we prepare the ground to \scite{3.8} with
IND$^+$ instead IND?
\bigskip

\demo{Proof}  We choose $\bar \lambda^n$ by induction on $n$.  For $n=1$ apply
\cite[Ch.II,1.5A]{Sh:g} to the sequence $\langle \lambda^1_i:i < \kappa_0
\rangle$, the ideal $\{A \subseteq \kappa_1:\text{max pcf}\{\lambda^1_i:i <
\kappa_0\} < \mu_0\}$ and the cardinal $\mu$ and get $\langle \lambda_{<i>}:
i < \kappa_0 \rangle$.  For $n+1$ for each $i < \lambda_n$ we apply
\cite[Ch.II,1.5A]{Sh:g} to $\langle \lambda_n:\eta \in \dsize \prod_{\ell < n}
\kappa_\ell \rangle$, the ideal $\{A \subseteq \dsize \prod_{\ell < n}
\kappa_\ell:\text{max pcf}\{\lambda_\eta:\eta \in A\} < \mu\}$ and the
cardinal $\lambda^n_i$ and we get $\langle \lambda_{\eta \char 94 <i>}:\eta
\in \dsize \prod_{\ell < n} \kappa_\ell \rangle$.
\hfill$\square_{\scite{3.4}}$
\enddemo
\bigskip

\proclaim{\stag{3.5} Claim}  In \scite{3.4}, from $\otimes_2$ we can deduce
\mr
\item "{$\otimes^3$}"  there are functions $f_{\ell,n}:\dsize \prod^n_{m =
\ell +1} \kappa_m \rightarrow I_\ell$ (for $\ell < N < \omega$) such that
for every $\eta \in \dsize \prod_{m < \omega} \kappa_m$ for some $\ell < n <
\omega$ we have $\eta(\ell) \in f{\ell,n}(\eta \restriction (\ell,n))$.
\endroster
\endproclaim
\bigskip

\demo{Proof}  Otherwise $\langle \lambda_{\eta \restriction n}:n < \omega
\rangle$ is a strictly decreasing sequence of cardinals. \nl
${{}}$  \hfill$\square_{\scite{3.5}}$
\enddemo
\bigskip

\definition{\stag{3.6} Definition}  1) IND$(\langle J_\varepsilon:\varepsilon
< \varepsilon^* \rangle)$ (note that $|\text{Dom}(J_\varepsilon)|$ is not
necessarily increasing with $\varepsilon$) means that eac $J_\varepsilon$ is
an ideal on Dom$(J_\varepsilon)$, say $\kappa_\varepsilon$ and
\mr
\item "{$(*)$}"  \ub{for every} sequence $\langle f_{\varepsilon,u}:
\varepsilon < \varepsilon^*,u \subseteq \varepsilon^* \backslash (\varepsilon
+1)$ finite$\rangle$ such that $f_{\varepsilon,u}$ a function from
$\dsize \prod_{\zeta \in u} \kappa_\zeta$ to $J_\varepsilon$ \ub{there is} an
increasing sequence $\varepsilon_0 < \varepsilon_1 < \ldots < \varepsilon_n
< \ldots < \varepsilon^*$ (for $n < \omega$) and $\alpha_\ell \in
\kappa_{\varepsilon_\ell}$ (for $\ell < \omega$) such that:
\sn
\item "{$(**)$}"  for $\ell < n < \omega$ we have $\alpha_\ell \notin
f_{\varepsilon_\ell,u}(\langle \alpha_{\varepsilon_{\ell +1}},\dotsc,
\alpha_{\varepsilon_n} \rangle)$ for $u = \{ \alpha_{\varepsilon_{\ell +1}},
\dotsc,\alpha_{\varepsilon_n}\}$.
\ermn
2) IND$^+(\langle J_\varepsilon:\varepsilon < \varepsilon^* \rangle)$ means
$J_\varepsilon$ is an ideal (on Dom$(J_\varepsilon)$ which is say
$\kappa_\varepsilon$) such that:
\mr
\item "{$(*)$}"  \ub{for every} sequence $\langle f_{\varepsilon,u}:
\varepsilon < \varepsilon^*,u \subseteq \varepsilon^* \backslash (\varepsilon
+1)$ finite$\rangle$ such that \nl
$f_{\varepsilon,u}:\dsize \prod_{\zeta \in u} 
\kappa_\zeta \rightarrow J_\varepsilon$
\ub{there is} $\langle \alpha_\varepsilon:\varepsilon < \varepsilon^* \rangle
\in \Pi \kappa_\varepsilon$ such that
\sn
\item "{$(**)$}"  for $\varepsilon < \varepsilon^*,u \subseteq \varepsilon^*
\backslash (\varepsilon +1)$ finite we have: $\alpha_\varepsilon \notin
f_{\varepsilon,u}(\ldots,\alpha_\zeta,\ldots)_{\zeta \in u}$.
\ermn
3) Let function $\langle J_\varepsilon:\varepsilon < \varepsilon^* \rangle$
be the set of $\bar f$ as in $(*)$ of part (1), i.e. $\bar f = 
\langle_{\varepsilon,u}:\varepsilon < \varepsilon^*$ and $u \subseteq
\varepsilon^* \backslash (\varepsilon +1)$ is finite$\rangle$ where
$f_{\varepsilon,u} \in \dsize \prod_{\zeta \in u} \kappa_\zeta$.  We say for
$\bar f \in \text{ function}\langle J_\varepsilon:\varepsilon < \varepsilon^*
\rangle$, that $\bar \varepsilon$ is candidate if $\bar \varepsilon$ is an
increasing sequence of length $\omega$ of ordinals $< \varepsilon^*$.  In 
this case we say that $\bar \alpha$ is $(\bar f,\bar \varepsilon)$-free if
$\bar \alpha \in \dsize \prod_{n,\omega} \text{ Dom}(J_{\varepsilon_n})$ and
the statement $(**)$ of part (1) holds.
\enddefinition
\bigskip

\demo{\stag{3.7} Observation}  1) If IND$(\bar J)$ where $\bar J = \langle
J_\varepsilon:\varepsilon < \varepsilon^* \rangle$ is as in \scite{3.6}, each
$J_\varepsilon$ is $(|\varepsilon^*|^{\aleph_0})^+$-complete, \ub{then} for
some $\varepsilon_0 < \varepsilon_1 < \ldots < \varepsilon^*$ we have
IND$^+(\langle J_{\varepsilon_n}:n < \omega \rangle)$. \nl
2) If IND$(\bar J)$ where $\bar J = \langle J_\varepsilon:\varepsilon <
\varepsilon^* \rangle$ as in \scite{3.6}, each $J_\varepsilon$ is
cov$(|\varepsilon^*|,\mu,\aleph_1,2)^+$-complete \ub{then} for some infinite
$u \in {\Cal S}_{< \mu}(\varepsilon^*)$ we have IND$(\bar J \restriction u)$.
\enddemo
\bn
Before we prove \scite{3.7}
\demo{\stag{3.8} Conjecture}  if IND$(\langle J_\varepsilon:\varepsilon <
\varepsilon^* \rangle),\varepsilon^* < \omega_1$ and each $J_\varepsilon$ is
$\aleph_1$-complete \ub{then} for some c.c.c. forcing $P$ we have:

$$
\Vdash_P ``\text{ for some } \varepsilon_0 < \ldots < \varepsilon_n <
\varepsilon_{n+1} < \ldots < \varepsilon^*, \text{ IND}(\langle
J_{\varepsilon_n}:n < \omega \rangle)".
$$
\enddemo
\bigskip

\remark{\stag{3.9} Remark}  In the proof of \scite{3.7}(2) it is enough to
demand on ${\Cal P}$:
\mr
\item "{$(*)$}"  if $\varepsilon_0 < \varepsilon_1 < \ldots < \varepsilon_n
< \ldots < \varepsilon^*$ (for $n < \omega$) then for some $b \in {\Cal P}$,
\nl
$(\exists^\infty n)\varepsilon_n \in b$
\ermn
this seems to weaken cov$(\ldots)$ but does not.
\endremark
\bigskip

\demo{Proof of 3.7}  1) Similar to the proof of part (2), as
cov$(\lambda,\aleph_1,\aleph_2,2) \le |\varepsilon^*|^{\aleph_0}$. \nl
2)  Let $\mu =: \text{ cov}|\varepsilon^*,\mu,\aleph_1,2)$ and let
${\Cal P} \subseteq [\varepsilon^*]^{< \mu}$ be of cardinality $\mu$
exemplifying its definition i.e. $(\forall a)[a \subseteq \varepsilon^* \and
|a| = \aleph_0 \Rightarrow (\exists b \in {\Cal P})[a \subseteq b]]$.

If for some $b \in {\Cal P}$, IND$(\bar J \restriction u)$ hold then we are
done.  Otherwise for each $b \in {\Cal P}$, we can find $\bar f^b = \langle
f^b_{\varepsilon,u}:\varepsilon \in b \text{ and } b \subseteq u \backslash
(\varepsilon +1)$ is finite$\rangle \in \text{ function}(J \restriction u)$
such that for no $\bar \varepsilon = \langle \varepsilon_n:n < \omega \rangle$
strictly increasing sequence of ordinals from $b$ and $\alpha_n \in
\text{ Dom}(J_{\varepsilon_n})$ (for $n < \omega$) do we have $n < \omega \and
u \subseteq \{\varepsilon_{n+1},\varepsilon_{n+2},\ldots\} \Rightarrow
\alpha_n \notin f^b_{\varepsilon,u}(\{\alpha_m:m \in u\})$.  Let us define for
$\varepsilon < \varepsilon^*$ and $u \subseteq \varepsilon^* \backslash
(\varepsilon +1)$ finite a function $f_{u,\varepsilon}$ from $\dsize \prod
_{\zeta \in u} \text{ Dom}(J_\zeta)$ to $J_\varepsilon$ by:
\mr
\item "{$(*)$}"  if $\alpha_\zeta \in \text{ Dom}(J_\zeta)$ for $\zeta \in u$
then $f_{\zeta,u}(\ldots,\alpha_\zeta,\ldots)_{\zeta \in u} =$ \nl
$\bigcup \{f^b_{\zeta,u}:\zeta \in b \text{ and } u \subseteq b 
\text{ and } b \in {\Cal P}\}$.
\ermn
(As each $J_\zeta$ is $\mu^+$-complete (by assumption) Rang$(f_{\zeta,u})
\subseteq J_\zeta$.  As IND$(\langle J_\varepsilon:\varepsilon < \varepsilon^*
\rangle)$ necessarily there is a strictly increasing $\langle \varepsilon_n:
n < \omega \rangle,\varepsilon_n < \varepsilon^*$, and $\alpha_n \in
\text{Dom}(J_{\varepsilon_n})$ (for $n < \omega$) such that:
\mr
\item "{$(**)$}"  $n < \omega,u \subset \{\varepsilon_{n+1},\varepsilon_{n+2},
\ldots\}$ finite, implies $\alpha_n \notin f_{\varepsilon_n,u}(\ldots,
\alpha_m,\ldots)_{\varepsilon_m \in u}$.
\ermn
By the choice of ${\Cal P}$ for some $b \in {\Cal P}$ we have 
$\{\varepsilon_n:n < \omega\} \subseteq b$, but then $\langle (\varepsilon_n,
\alpha_n):n < \omega \rangle$ contradict the choice of $\bar f^b = \langle
f^b_{\zeta,u}:\zeta \in b,u \subseteq b \backslash (\zeta +1)$ finite$\rangle$.\hfill$\square_{\scite{3.7}}$
\enddemo
\bigskip

\demo{\stag{3.10} Conclusion}  1) If IND$^+(\langle I_n:n < \omega \rangle)$
and Dom$(I_n) = \kappa_n$ \ub{then}
\mr
\item "{$(a)$}"  the conclusion of \scite{3.5} (i.e. $\otimes^3$ there) fails,
hence $\otimes^2$ of \scite{3.4} fails hence $\otimes^1$ of \scite{3.4} fails
\sn
\item "{$(b)$}"  if $\lambda > \dsize \sum_{n < \omega} \kappa_n,\kappa_n
< \text{ cf}(\kappa_{n+1})$, \ub{then} for every $n$ large enough for no
$\lambda_i \in (\dsize \sum_{n < \omega} \kappa_n,\lambda) \cap \text{ Reg}$
(for $i < \kappa_n$) is $\dsize \prod_{i < \kappa_n} \lambda_i/I_n
\lambda$-directed.
\ermn
2) If we weaken the assumption to IND$(\langle I_n:n < \omega \rangle)$ then
in (b) we have just for arbitrarily large $n < \omega$.
\enddemo
\bigskip

\demo{Proof}  1) ? 
\mr
\item "{$(a)$}"  straight
\sn
\item "{$(b)$}"  our problem is to get $\mu_1$, which is not serious.
\ermn
2) Similarly.  \hfill$\square_{\scite{3.10}}$
\enddemo
\bigskip

\demo{\stag{3.11} Conclusion}  1) Assume $\langle \kappa_\varepsilon:
\varepsilon < \delta \rangle$ is strictly increasing, $|\delta| \le \sigma <
\kappa_0,\kappa = \dsize \sum_{i < \delta} \kappa_i$ and IND$(\kappa,\sigma)$.
If $\lambda > \kappa$ \ub{then} for every large enough $\varepsilon < \delta$,
there are no $\lambda_\alpha \in (\kappa,\lambda) \cap \text{ Reg}$ for
$\alpha < \kappa_\varepsilon$ such that $\dsize \prod_{\alpha < \kappa_
\varepsilon} \lambda_\alpha/[\kappa_\varepsilon]^{\le \sigma}$ is
$\lambda$-directed recalling $[\kappa_\varepsilon]^{\le \sigma} = \{a
\subseteq \kappa_\varepsilon:|a| \le \sigma\}$. \nl
2) If IND$(\kappa)$, cf$(\kappa) = \aleph_0 = \sigma,\kappa = \dsize \sum
_{n < \omega} \kappa_n,\kappa_n < \kappa_{n+1}$ \ub{then} the conclusion of
1) holds. \nl
3) If IND$(\kappa,\sigma),\delta = \omega,\kappa_\varepsilon = \kappa$,
\ub{then} the conclusion of (1) holds. \nl
4) If $\varepsilon < \lambda_0,\langle \lambda_\varepsilon:\varepsilon <
\varepsilon^* \rangle$ is a strictly increasing sequence of regular
cardinals, $J_\varepsilon$ is an ideal on $\lambda_\varepsilon$ extending
$J^{\text{bd}}_{\lambda_\varepsilon}$ and we have IND$(\langle J_\varepsilon:
\varepsilon < \varepsilon^* \rangle)$ \ub{then} for some $\varepsilon' \le
\varepsilon^*$, cf$(\varepsilon') = \aleph_0$ and not only IND$(\langle
J_\varepsilon:\varepsilon' \rangle)$ but we can even demand $\varepsilon' =
\dbcu_{n < \omega} \varepsilon_n$. \nl
5) IND$(\langle \lambda_n:n < \omega \rangle)$ iff IND$^+(\langle[\lambda_n]
^{\le \lambda_{n-1}}:n < \omega \rangle)$ (stipulating $\lambda_{-1} =
\aleph_0$). \nl
6) If IND$^0(\lambda,\sigma)$ then 
IND$^+(\langle \lambda:n < \omega \rangle)$.
\enddemo
\bigskip

\demo{Proof}  Check.
\enddemo
\bn
\ub{\stag{3.12} Discussion}:  Let $\bar J = \langle J_\varepsilon:\varepsilon
< \varepsilon^* \rangle$ and assume IND$(\bar J)$. \nl
1) Note that if $P$ is a $\theta$-c.c. forcing notion, each $J_\varepsilon$
is $\theta$-complete then for any $\bar f = \langle f_{\varepsilon,u}:
\varepsilon < \varepsilon^*,u \in [\varepsilon^* \backslash (\varepsilon +1)]
^{< \aleph_0} \rangle \in V^P$ as in Definition \scite{3.6} we can find
$\bar f' = \langle f'_{\varepsilon,u}:\varepsilon < \varepsilon^*,u \in
[\varepsilon^* \backslash (\varepsilon +1)]^{< \aleph_0} \rangle \in V$ such
that for every $\bar \alpha \in \dsize \prod_{\zeta \in u} \text{ Dom}
(J_\zeta)$ we have $f_{\varepsilon,u}(\bar \alpha) \subseteq
f'_{\varepsilon,u}(\bar \lambda)$, so we can consider only $\bar f \in V$.
For each such $\bar f$ let $A_{\bar f} = \{v:$ for some strictly increasing
sequence $\bar \varepsilon = \langle \varepsilon_n:n < \omega \rangle$ of
ordinals $< \varepsilon^*$ and $\bar \alpha \in \dsize \prod_{n < \omega}
\text{ Dom}(J_{\varepsilon_n})$ the conclusion $(**)$ of Definition
\scite{3.6} holds and $v = \{\varepsilon_\ell:\ell < m\}$ for some $m <
\omega\}$.

For $\bar g,\bar f \in \text{ function}(\bar J)$ where $\bar J = \langle
J_\varepsilon:\varepsilon < \varepsilon^* \rangle$, let $\bar g \le \bar f$
iff for every $\varepsilon < \varepsilon^*$ and $u \in [\varepsilon^*
\backslash (\varepsilon +1)]^{< \aleph_0}$ we have $\bar \alpha \in \dsize
\prod_{\varepsilon \in u} \text{ Dom}(J_\varepsilon) \Rightarrow
g_{\varepsilon,u}(\bar \alpha) \subseteq f_{\varepsilon,u}(\bar \alpha)$.
Clearly $\bar g \le \bar f \Rightarrow A_{\bar f} \subseteq A_g$. \nl
2) In $V$ we can define a filter $D$ on $[\dbcu_{\varepsilon < \varepsilon^*}
\text{ Dom}(\bar J_\varepsilon)]^{< \aleph_0}$:

$$
A \in D \text{ \ub{iff} for some } \bar f \in \text{ function}(\langle
J_\varepsilon:\varepsilon < \varepsilon^* \rangle) \text{ we have } A_f
\subseteq A.
$$
\mn
Now $D \subseteq {\Cal P}([\dbcu_{\varepsilon < \varepsilon^*}
\text{ Dom}(J_\varepsilon)]^{< \aleph_0})$ and $D$ is upward closed trivially.
Also $D$ is closed under intersection of countable many members if each
$J_n$ is $\aleph_1$-closed (similarly $\sigma$-closed) because if $A_n \in
D$ let $\bar f^n \in \text{ function}(\bar J)$ be such that $A_{\bar f^n}
\subseteq A$.  Now for some $g \in \text{ function}(J)[n < \omega \Rightarrow
\bar f^n \subseteq \bar g]$, hence $A_g \subseteq A_{\bar f^n}$, so
$A_{\bar g} \subseteq A_n$ for $n < \omega$ and obviously $A_{\bar g} \in D$.
Lastly $\emptyset \notin D$ as IND$(\bar J)$ holds.
\bigskip

\proclaim{\stag{3.13} Claim}  Suppose for $\alpha < \alpha^*,\bar I^\alpha =
\langle I^\alpha_n:n < \omega \rangle,\kappa = \sup\{|\text{Dom}(I^\alpha_n)|:
\alpha < \alpha^*\}$, IND$^+(\langle I^\alpha_n:n < \omega \rangle)$, and:
\mr
\item "{$(*)$}"  if $\alpha < \alpha^*,f_n:\text{Dom}(I^\alpha_n) \rightarrow
\text{ Ord}$ \ub{then} for some $n(*) < \omega,\beta < \alpha^*$ and ordinal
$\gamma$
\endroster
$$
I^\beta_{1+n} \cong I^\alpha_{n(*)+1+n} \restriction \{x \in \text{ Dom}
(I^\alpha_{n(*)+1+n}):f_{1+n}(x) > \gamma\}
$$

$$
I^\beta_0 = I^\alpha_{n(*)} \restriction \{x \in \text{ Dom}(I^\alpha_{n(*)}):
f_0(x) < \gamma\}.
$$
\mn
\ub{Then} for no $\lambda > \kappa$ and $\alpha < \alpha^*$ do we have for
every $n < \omega$:

$$
x \in \text{ Dom}(I^\alpha_n) \Rightarrow \lambda^n_x \in (\kappa,\lambda)
\cap \text{ Reg and } \dsize \prod_{x \in \text{ Dom}(I^\alpha_n)}
\lambda^n_x/I^\alpha_n \text{ is } \lambda^+\text{-directed}.
$$
\endproclaim
\bigskip

\demo{Proof}  No new point, like the proof of \scite{2.6}.
\enddemo
\bigskip

\remark{Remark}  1) This claim is used in the proof of \scite{5.4}. \nl
2) If in \scite{3.13}, $\alpha^* = 1$ and $I^\alpha_{n+1}$ is
$\lambda_n$-complete, $\lambda_n > |\text{Dom}(I^\alpha_\ell))|$ for
$\ell < n$ then $(*)$ there holds.
\endremark
\bn
\ub{Question}:  If IND$(\lambda,\sigma)$, cf$(\lambda) = \aleph_0$ do we have
IND$(\langle J^{\text{bd}}_{\lambda_n}:n < \omega \rangle)$ for some
$\lambda_n = \text{ cf}(\lambda_n) < \lambda,\sigma < \lambda_n$?
\bigskip

\proclaim{\stag{3.14} Claim}  In \scite{2.6}(1) we can deduce even
IND$(\langle J^{\text{bd}}_\theta:\theta \in {\frak a} \rangle)$.
\endproclaim
\bigskip

\demo{Proof}  Reread the proof of \scite{2.6}: let $f_{u,\varepsilon}$ be as
in Definition \scite{3.6}, so without loss of generality 
Rang$(f_{\varepsilon,u}) \subseteq
\kappa_\varepsilon$ (as $\alpha = \{\beta:\beta < \alpha\})$ and $M = 
(\kappa,\kappa_\varepsilon,f_{u,\varepsilon})_{\varepsilon,u}$.  Now repeat
the proof of \scite{2.6}.  \hfil$\square_{\scite{3.14}}$
\enddemo
\bn
Next we improve the ideals from ``bounded" to ``nonstationary"
\proclaim{\stag{3.15} Claim}  1) Assume $\lambda > \text{ cf}(\lambda),
|{\frak a}| < \text{ min}({\frak a})$ and $\lambda = \sup({\frak a})$ and \nl
rk$'_{\langle J_{< \theta}[{\frak a}]:\theta \in \text{pcf}({\frak a})
\rangle}({\frak a}) \ge |{\frak a}|^+,\sigma^* \in (|{\frak a}|,\text{min}
({\frak a})) \cap \text{ Reg}$, and

$$
I_\theta =: \{S:S \subseteq \theta \text{ and } \{\delta \in S:\text{cf}
(\delta) = \sigma^*\} \text{ is not stationary}\}
$$
\mn
\ub{then} IND$(\langle I_\theta:\theta \in {\frak a} \rangle)$. \nl
2) Moreover for any sequence $\bar H = \langle H_\theta:\theta \in {\frak a}
\rangle$, where ${\frak a} \subseteq \lambda,H_\theta$ a function from
$[\lambda]^{< \aleph_0}$ to $I_\theta$ and ${\frak c},{\frak a} \subseteq
{\frak c} \subseteq \text{ pcf}({\frak a}),|{\frak c}| < \text{ min}
({\frak a})$ we can find $\bar \alpha = \langle \alpha_\tau:\tau \in {\frak a}
\rangle \in \Pi {\frak a}$ such that defining for ${\frak b} \subseteq
{\frak a}$

$$
c \ell_{\bar H,\bar \alpha}({\frak b}) = \{\tau \in {\frak a},\alpha_\tau
\in H_\tau(\bar \alpha \restriction {\frak e}) \text{ for some finite }
{\frak e} \subseteq {\frak b}\}
$$
\mn
we have
\mr
\item "{$(*)$}"  $\theta \in {\frak c} \and {\frak b} \in J_{< \theta}
[{\frak a}] \Rightarrow c \ell_{\bar H,\bar \alpha}({\frak b}) \in
J_{< \theta}[{\frak a}]$.
\endroster
\endproclaim
\bigskip

\demo{Proof}  1) We can prove it from part 2) exactly as in the proof of
\scite{2.6}(1). \nl
2) Let $\bar{\frak b} = \langle {\frak b}_\theta[{\frak a}]:\theta \in
\text{ pcf}({\frak a}) \rangle$ be a generating sequence for ${\frak a}$,
without loss of generality max pcf$({\frak a}) \in {\frak c}$ and

$$
[\theta \in {\frak c} \backslash \text{ max pcf}({\frak a}) \Rightarrow
\min(\text{pcf}({\frak a}) \backslash \theta^+) \in {\frak c}].
$$
\enddemo
\bn
Before we continue, recall that we know: \nl
\ub{\stag{3.16} Fact}  If $|{\frak a}| < \min({\frak a}),{\frak a} \subseteq
\text{ Reg and } \langle {\frak b}_\sigma[{\frak a}]:\sigma \in \text{ pcf}
({\frak a}) \rangle$ a generating sequence for ${\frak a}$ then
\mr
\item "{$(*)_{\frak a}$}"  there is $\langle \bar f^\theta:\theta \in
\text{ pcf}({\frak a}) \rangle$ such that
{\roster
\itemitem{ $(a)$ }  $\bar f^\theta = \langle f^\theta_\alpha:\alpha < \theta
\rangle$ is $<_{J_{< \theta}[{\frak a}]}$-increasing
\sn
\itemitem{ $(b)$ }  $\bar f^\theta$ is cofinal in $(\Pi {\frak a},
<_{J_\theta[{\frak a}]})$ where $J_\theta[{\frak a}] =: J_{< \theta}
[{\frak a}] + ({\frak a} \backslash {\frak b}_\theta[{\frak a}])$
\sn
\itemitem{ $(c)$ }  if $\sigma \in \theta \cap {\frak c},{\frak a} < \theta$
and ${\frak b} = {\frak b}_\sigma[{\frak a}]$ \ub{then} for some $n <
\omega$ and $\alpha_\ell < \theta_\ell \le \sigma$ (for $\ell \le n$) we have
$f^\theta_\delta \restriction {\frak b} = (\max\{f^{\theta_\ell}_{\alpha_\ell}
:\ell < n\}) \restriction {\frak b}$ (the max is pointwise)
\sn
\itemitem{ $(d)$ }  if $\delta < \theta \in \text{ pcf}({\frak a})$, 
cf$(\delta) \in (|{\frak a}|,\min({\frak a}))$ \ub{then} for every $\tau \in
{\frak a}$
$$
f^\theta_\delta(\tau) = \min\{\cup\{f^\theta_\alpha(\tau):\alpha \in C\}:
C \text{ a club of } \delta\}
$$
\ub{provided that} the function defined satisfies condition (c) above
\endroster}
(exist by \cite[Ch.VIII,\S1]{Sh:g} or \cite[6.7x]{Sh:430} we choose by induction
on $\theta$).
\ermn
Let

$$
\align
S^{\text{gd}}_\theta =: \bigl\{ \delta < \theta:&\text{ cf}(\delta) \in
(|{\frak a}|,\min({\frak a})),f^\theta_\delta \text{ is a }
<_{J^{\text{bd}}_{< \theta}[{\frak a}]} \text{-eub of }
  \bar f^\theta \restriction \delta \text{ and} \\
  &\, \{\tau \in {\frak a}:\text{cf}(f^\theta_\alpha(\tau)) = \text{ cf}
(\delta)\} = {\frak b}_\theta[{\frak a}] \text{ mod } J_{< \theta}
[{\frak a}] \bigr\}
\endalign
$$
\mn
(alternatively use simultaneous witnesses for $I[\theta]$ as in
\cite[\S1]{Sh:420}. \nl
Note:
\mr
\item "{$(*)_{\frak a}$}" $\,\,\, (e) \quad$ if $E_\tau$ is a club of $\tau$
for $\tau \in {\frak a}$ and $\theta \in \text{ pcf}({\frak a})$ \ub{then}
for some club $E$ 

$\qquad$ of $\theta$:
$$
\delta \in E \cap S^{\text{gd}}_\theta \Rightarrow \{\tau \in {\frak a}:
f^\theta_\delta(\tau) \in E_\tau\} = {\frak a} \text{ mod } J_\theta
[{\frak a}].
$$
\ermn
Let $\bar H = \langle H_\tau:\tau \in {\frak a} \rangle$ be as in the claim,
so $H_\tau$ is a function from $[\lambda]^{< \aleph_0}$ to $I_\tau$.  Now
for every $\theta \in \text{ pcf}({\frak a})$ and $\alpha < \theta$ and
$\tau \in {\frak a}$ we define $A^{\theta,\tau}_\alpha = A^\tau(f^\theta
_\alpha) \in I_\tau$ as

$$
\bigcup\{H_\tau(u):u \subseteq \text{ Rang}(f^\theta_\alpha)
\text{ is finite}\}
$$
\mn
(as $I_\tau$ is $\tau$-complete, $\tau > |{\frak a}| = |\{u:u \subseteq
\text{ Rang}(f^\theta_\alpha) \text{ finite}\}|$, really $A^\theta_\alpha
\in I_\tau$).

Now for each $\alpha < \theta \in \text{ pcf}({\frak a})$ and $\sigma \in
\text{ pcf}({\frak a})$ by $(*)_{\frak a}(e)$ applied with $\sigma,\langle
A^{\theta,\tau}_\alpha:\tau \in {\frak a} \rangle$ here standing for $\theta,
\langle E_\tau:\tau \in {\frak a} \rangle$ there, we get a club
$C_{\theta,\sigma,\alpha}$ of $\sigma$.

For each $\sigma \in \text{ pcf}({\frak a})$ let $C_\sigma =
\dbca_{\theta \in {\frak c} \cap \sigma,\alpha < \theta} C_{\theta,\sigma,
\alpha}$ (note: $|{\frak c} \cap \sigma| < \sigma$), so $C_\sigma$ is a
club of $\sigma$.  Lastly let $\sigma^* = \text{ cf}(\sigma^*) \in
(|{\frak a}|,\min({\frak a}))$ and $\langle N_i:i \le \sigma^* \rangle$ be
as in \cite[Ch.VIII,1.2,1.4]{Sh:g}, so
\mr
\item "{$\oplus_1$}"  $N_i \prec ({\Cal H}(\chi),\in,<^*_\chi)$ is
increasing continuous, for $i \le \sigma^*,\|N_i\| = \sigma^*$, \nl
$\langle N_j:j \le i \rangle \in N_{i+1}$.
\ermn
Let $\alpha_\theta =: \sup(N_{\sigma^*} \cap \theta)$ for $\theta \in
{\frak c}$; so $\alpha_\theta \in S^{\text{gd}}_\theta$ for every $\theta \in
{\frak c}$, and we shall show that $\bar \alpha = \langle \alpha_\tau:\tau
\in {\frak a} \rangle$ is as required.

For each $\sigma,\theta \in {\frak c},\sigma < \theta$ we have: $\alpha
_\theta \in C_\theta$ hence $\alpha_\theta \in C_{\theta,\sigma,
\alpha_\sigma}$ hence

$$
\{\tau \in {\frak a}:f^\sigma_{\alpha_\sigma}(\tau) \in \tau \backslash
A^{\sigma,\tau}_{\alpha_\sigma}\} = {\frak b}_\tau[{\frak a}] \text{ mod }
J_\sigma[{\frak a}].
$$
\mn
The rest should be clear.  \hfill$\square_{\scite{3.15}}$
\bn
We next point out another connection; if the rank is small and
$|\text{pcf}({\frak a})|$ is large, then we have a case of ``$\Pi {\frak d}/
{\Cal S}_{\le \lambda}({\frak d})$ has large true cofinality".
\proclaim{\stag{3.17} Claim}  If $\zeta > \text{\rm rk\/}'({\frak a},\langle
J_{< \theta}[{\frak a}]:\theta \in \text{ pcf}({\frak a}) \rangle)$ and
$\bar \lambda = \langle \lambda_\varepsilon:\varepsilon \le \zeta \rangle$ 
is strictly increasing and $|\text{pcf}({\frak a})| \ge \lambda_\zeta$
\ub{then} for some $\varepsilon < \zeta$ and ${\frak c} \subseteq {\frak a}$
we have $|\text{pcf}({\frak c})| \ge \lambda_{\varepsilon +1}$ and
${\frak d} \in J_*[\text{pcf}({\frak c})\ \Rightarrow \Pi {\frak d}/{\Cal S}
_{\le \lambda_\varepsilon}({\frak d})$ has the true cofinality which is
max pcf$({\frak c})$.
\endproclaim
\bigskip

\demo{Proof}  If not, prove by induction on $\varepsilon \le \zeta$ that
\mr
\item "{$(*)$}"  if ${\frak c} \subseteq {\frak a},\varepsilon =
\text{ rk}'({\frak c},\langle J_{< \theta}[{\frak a}]:\theta \in
\text{ pcf}({\frak a}) \rangle)$ \ub{then} $|\text{pcf}({\frak c})| <
\lambda_{\varepsilon + 1}$.
\ermn
Let $J_0 = \{{\frak d}:{\frak d} \subseteq {\frak c},\varepsilon >
\text{ rk}'({\frak d},\langle J_{< \theta}[{\frak a}]:\theta \in \text{ pcf}
({\frak a}) \rangle)\}$.  By the induction hypothesis $[{\frak d} \in J_0
\Rightarrow |\text{pcf}({\frak d}) \le \lambda_\varepsilon]$.  Let $J$ be
the ideal on ${\frak c}$ which $J_0$ generates.  We have: $[{\frak d} \in J
\Rightarrow |\text{pcf}({\frak d})| \le \lambda_\varepsilon]$, so by the
assumption toward contradiction ${\frak c} \notin J$.  Let $\theta = 
\text{max pcf}({\frak c})$.  So by the definition of rk$':J_{< \theta}
[{\frak c}] \subseteq J$.  Hence (see \cite[Ch.VIII,\S3]{Sh:g} or
\cite{Sh:E11} for ${\frak d} \subseteq \text{ pcf}({\frak c}),{\frak d} \in
J_*[\text{pcf}({\frak a})]$ we have $\Pi {\frak d}/J^*_{< \theta}
[\text{pcf}({\frak c})]$ has the true cofinality $\theta$ and $J_*[\text{pcf}
({\frak a})]$ is generated by $\{\text{pcf}({\frak b}):{\frak b} \in
J_{< \theta}[{\frak c}]\}$.  So the conclusion holds.
\hfill$\square_{\scite{3.17}}$
\enddemo
\newpage

\head {\S4 Sticks and BA's} \endhead  \resetall
\bigskip

\proclaim{\stag{4.1} Lemma}  Assume $\theta \le \mu < \lambda \le \lambda^*,J$
an ideal on $\theta$ and assume
\mr
\item "{$\otimes^J_{\theta,\mu,\lambda,\lambda^*}$}"  if 
$n < \omega,{\frak a}_i
\in [\text{Rang } \cap \lambda^+ \backslash \mu^+]^n$ for $i < \theta$
\ub{then} \nl
$\{a \in J:\text{max pcf}(\dbcu_{i \in a} {\frak a}_i) \le
\lambda^*\}$ is generated by $\le \mu$ sets. \footnote{see \scite{4.2}(3)}
\ermn
\ub{Then} there is a set $H$ such that
\mr
\item "{$(a)$}"  $H$ a set of partial functions from $\theta$ to
$[\lambda]^{\le \mu}$
\sn
\item "{$(b)$}"  $|H| \le \lambda^*$
\sn
\item "{$(c)$}"  for every function $g:\theta \rightarrow \lambda$ we can
find $h$ and $\bar{\frak a} = \langle {\frak a}_i:i < \theta \rangle$ such
that
{\roster
\itemitem{ $(i)$ }  ${\frak a}_i$ is a finite set of regular cardinals from
$(\theta,\lambda]$
\sn
\itemitem{ $(ii)$ }  $h$ is a function from $\theta$ to $[\lambda]^{\le \mu}$
such that $\dsize \bigwedge_i g(i) \in h(i)$
\sn
\itemitem{ $(iii)$ }  for any $n < \omega$ and $a \subseteq \theta$: \nl
\ub{if} $(\forall i \in a)[|{\frak a}_i| \le n]$ and max pcf$(\dbcu_{i \in a}
{\frak a}_i) \le \lambda^*$ \ub{then} for some $b$ satisfying $a \subseteq b
\subseteq \theta$ we have $h \restriction b \in H$.
\endroster}
\endroster
\endproclaim
\bigskip

\demo{Proof}  Like \cite[\S2]{Sh:430}.
\enddemo
\bigskip

\remark{\stag{4.2} Remark}  1) But we can then change the bound (in clause
(c)(ii)) to $h(i) \in [\lambda]^{< \mu}$.  Then $\otimes_{\theta,\mu,\lambda,
\lambda^*}$ is changed to
\mr
\item "{$\otimes'_{\theta,\mu,\lambda,\lambda^*}$}"  if 
$n < \omega,{\frak a}_i \in [\text{Reg } \cap \lambda^+]^n$ for $i < \theta$
then \nl
$\{a \in J:\text{for some } \mu_0 < \mu,\text{max pcf}(\dbcu_{i \in a}
{\frak a}_i \backslash \mu^+_0) \le \lambda^*\}$ \nl
is generated by $< \mu$ sets.
\ermn
2) We can weaken $\otimes$ to
\mr
\item "{$\otimes^J_{\theta,\mu,\lambda,\lambda^*}$}"  if 
$n < \omega,{\frak a}_i \in [\text{Reg } \cap \lambda^+]^n$ for $i < \theta$
then \nl
$\{a \in J:\text{max pcf}(\dbcu_{i < a} {\frak a}_i) \le \lambda^*\}$ \nl
is generated (as an ideal) by some ${\Cal P} \subseteq J$ such that
$$
\kappa \in \dbcu_{i < \theta} {\frak a}_i \Rightarrow \kappa <
\{a \in {\Cal P}:\kappa \in \dsize \bigvee_{i \in a} {\frak a}_i\}.
$$
\ermn
3) Instead $\otimes_{\theta,\mu,\lambda,\lambda^*}$ we can define a game
\mr
\item "{$\otimes_{\theta,\mu,\lambda,\lambda^*}[{\frak D}]$}"  First player
has no winning strategy in the game defined below
\sn
\item "{$GM'_{\theta,\mu,\lambda,\lambda^*}[{\frak D}]$}"  The play lasts
$\omega$ moves, in the $n$-th move: \nl
first player chooses $\bar \lambda^n = \langle \lambda^n_i:i \in A_n \rangle,
A_0 = \theta, \dsize \bigwedge_{m < n} A_m \subseteq A_n,\dsize \bigwedge
_{m < n} \dsize \bigwedge_{i \in A_n} \lambda^n_i < \lambda^m_i,\lambda^n_i=
\text{ cf}(\lambda^n_i) \in (\mu,\lambda]$ and \nl
second player chooses an ideal $J_n$ on $A_n,J_n \subseteq \{a \subseteq A_n:
\text{max pcf}\{\lambda^n_i:i \in a\} \le \lambda^*\},J_n$ generated by
$\le \theta$ sets. \nl
In the end (clearly $\dbca_{n < \omega} A_n = \emptyset$) they produce the
ideal $J$, the one generated by $\{a \subseteq \theta:\text{for some } n,a
\subseteq A_n \backslash A_{n+1} \text{ and } a \in J_n\}$. \nl
Second player wins if $J \in {\frak D}$.
\endroster
\endremark
\bigskip

\definition{\stag{4.3} Definition}  Assume $\bar J = \langle J_\ell:\ell < 3
\rangle$, where $J_0 \subseteq J_1 \subseteq J_2 \subseteq {\Cal P}(\theta)$,
each $J_\ell$ is downward closed (usually is an ideal); we let $J^+_\ell =:
{\Cal P}(\theta) \backslash J_\ell$. \nl

$$
\align
1) \qquad \text{ dcf}_{\bar J}(\lambda < \mu) = \min\biggl\{ |{\Cal F}:&
{\Cal F} \text{ is a family of functions each with domain from} \\
  &J^+_1 \text{ and range included in } [\lambda]^{< \mu} \text{ such that:}\\
  &(*)_{\Cal F} \quad \text{ for every } b \in J^+_2 \text{ and } f \in
{}^b \lambda \text{ for some } a \in J^+_1 \\
  &\quad \qquad \,\,\text{ and } 
g \in {\Cal F} \cap {}^a \lambda \text{ we have }
(\forall^{J_0} i \in a)(i \in b \and g(i) \in f(i)) \\
  &\quad \qquad \,\, \text{ i.e. } 
\{i:i \in a \text{ and } i \notin b \vee g(i)
\notin f(i)\} \in J_0 \biggr\}
\endalign
$$

$$
\align
2) \qquad \text{ ecf}_{\bar J}(\lambda < \mu) = \min \biggl\{ |{\Cal F}:&
{\Cal F} \text{ is a family of functions each with domain from} \\
  &J^+_1 \text{ and range included in } [\lambda]^{< \mu} \text{ such that:}\\
  &(*)_{\Cal F} \quad \text{ for every } b \in J^+_2 \text{ and } f \in
{}^b \lambda \text{ for some } a \in S^+_1 \\
  &\qquad \quad \text{ and } g \in {\Cal F} \cap 
{}^\alpha \lambda \text{ we have} \\
  &\qquad \quad (\exists^{J_0} i \in a)(i \in b \and g(i) \in f(i)) \\
  &\qquad \quad \,\, \text{ i.e. } 
\{i:i \in a,i \in b \text{ and } g(i) \in f(i)\} \notin J_0 \biggr\}
\endalign
$$
\mn
3) xcf$_{\bar J}(\lambda, < \mu^+)$ is written xcf$_{\bar J}(\lambda,\le \mu)$
or xcf$_{\bar J}(\lambda,\mu)$.  Also xcf$_{\bar J}(\lambda,1)$ is written
xcf$_{\bar J}(\lambda)$.  Also xcf$_{< [\lambda]^{< \sigma}, 
[\lambda]^{< \theta},[\lambda]^{< \chi}}(\lambda, < \mu)$ is written
xcf$_{\chi,\theta,\sigma}(\lambda,< \mu)$, etc.
\enddefinition
\bigskip

\definition{\stag{4.4} Definition}  $\stick_{\lambda,\mu,\theta} =: \min
\{|{\Cal P}|:{\Cal P} \subseteq [\lambda]^\theta$ is such that for every
$A \in [\lambda]^\mu$ for some $a \in {\Cal P}$ we have $a \subseteq A\}$.
If $\mu = \lambda$ we may omit it.  Let $\stick_\lambda$ mean
$\stick_{\lambda,\aleph_1}$ and $\stick = \stick_{\aleph_1}$.
\enddefinition
\bigskip

\proclaim{\stag{4.5} Claim}  1) $\stick_{\lambda,\mu,\theta} = \text{ dcf}
_{\langle\{\emptyset\},[\lambda]^{< \theta},[\lambda]^{< \mu} \rangle}
(\lambda)$ when $\theta \le \mu \le \lambda$. \nl
2) ecf$_{\mu,\theta,\sigma}(\lambda,\theta) = \min\{|{\Cal P}|:{\Cal P}
\subseteq [\lambda]^\theta$ and for every $A \in [\lambda]^\mu$ for some $a
\in {\Cal P}$ we have $|A \cap a| \ge \sigma\}$. \nl
Now we can phrase the analog of \scite{4.1} for dcf.
\endproclaim
\bigskip

\proclaim{\stag{4.6} Lemma}  Assume $\theta < \mu < \lambda \le \lambda^*,J$
an ideal on $\theta$ and assume
\mr
\item "{$\otimes^J_{\theta,\mu,\lambda,\lambda^*}$}"  if
$n < \omega,{\frak a}_i \in [\text{Rang } \cap \lambda^+ \backslash \mu^+]^n$
for $i < \theta$ then \nl
$\{a \in J:\text{max pcf } \dbcu_{i \in a} {\frak a}_i \le \lambda^*\}$ is
generated by $\le \mu$ sets. \footnote{see \scite{4.2}(3)}
\ermn
\ub{Then} there is a set $H$ such that
\mr
\item "{$(a)$}"  $H$ a st of partial functions from $\theta$ to
$[\lambda]^{\le \mu}$
\sn
\item "{$(b)$}"  $|H| \le \lambda^*$
\sn
\item "{$(c)$}"  for every function $g:\theta \rightarrow \lambda$ we can 
find $h$ and $\bar{\frak a} = \langle {\frak a}_i:i < \theta \rangle$ such
that
{\roster
\itemitem{ $(i)$ }  ${\frak a}_i$ is a finite set of regular cardinals from
$(\mu,\lambda]$
\sn
\itemitem{ $(ii)$ }  $h$ is a function from $\theta$ to $[\lambda]^{\le \mu}$
such that $\dsize \bigwedge_i g(i) \in h(i)$
\sn
\itemitem{ $(iii)$ }  for any $n < \omega$ and $a \subseteq \theta$: \nl
\ub{if} $(\forall i \in a)[{\frak a}_i| \le n]$ and max pcf$(\dbcu_{i \in a}
{\frak a}_i) \le \lambda^*$ \ub{then} for some $b,a \subseteq b \subseteq
\theta$ we have $h \restriction b \in H$.
\endroster}
\endroster
\endproclaim
\bigskip

\proclaim{\stag{4.7} Claim}  Assume:
\mr
\item "{$(*)_0$}"  $\aleph_0 < \aleph_{\alpha(*)} \le \stick$ and:
$\aleph_{\alpha(*)} < \aleph_{\omega_2}$ or at least
$$
\text{cov}(\aleph_{\alpha(*)},\aleph_2,\aleph_2,\aleph_1) \le \stick
$$
\sn
\item "{$(*)_1$}"  ${\frak a} \subseteq \text{ Reg} \cap \stick \backslash
\aleph_{\alpha(*)+1} \and |{\frak a}| \le \aleph_0 \Rightarrow |\text{pcf}
({\frak a})| \le \aleph_{\alpha(*)}$
\sn
\item "{$(*)_2$}"  if $\lambda_i \in 
(\aleph_1,\stick) \cap \text{ Reg}$ for
$i < \omega_1$ \ub{then} for some $a \in [\omega_1]^{\aleph_1}$, for every
$b \in [a]^{\aleph_0}$ we have max pcf$(\{\lambda_i:i \in b\}) \le \stick$
\ub{then} $\stick_{\stick} = \stick$.
\endroster
\endproclaim
\bigskip

\remark{\stag{4.8} Remark}  1) This means that the conclusion holds except
when some dubious statements on pcf, ones which have high consistency
strength (or are inconsistent) and $\stick$ is somewhat large. \nl
2) There are obvious monotonicity properties and ecf$_{\bar J}(\lambda,< \mu)
\le \text{ dcf}_{\bar J}(\lambda,< \mu)$.
\endremark
\bigskip

\demo{Proof}  Let $\theta = \aleph_1,\mu = \aleph_{\alpha(*)},\lambda =
\stick,\lambda^* = \stick,J = [\omega_1]^{\le \aleph_0}$.  Apply 4.2,
more exactly the variant \scite{4.2}(2).  The assumption $\otimes^J_{\theta,
\mu,\lambda,\lambda} = \otimes^J_{\aleph_1,\aleph_{\alpha(*)},\stick,
\stick}$ holds by $(*)_1$.  So let $H \subseteq \{h:h:\aleph_1 \rightarrow
[\stick]^\mu\},|H| = \lambda^* = \lambda$ be as in the conclusion there.
Let $\chi = \beth^+_7,{\frak B} \prec ({\Cal H}(\chi),\in,<^*_\chi),
|{\frak B}| = \lambda,\lambda + 1 \subseteq {\frak B},H \in {\frak B}$.  We
want to show ${\Cal P} = {\frak B} \cap [\lambda]^{\aleph_0}$ exemplifies
$\stick_\lambda = \lambda$.  So let $g:\theta \rightarrow \lambda$, so
there are $h,\left< \langle \lambda^n_i:n < n_i \rangle:i < \omega_1 \right>$,
as there.  Let $n^*$ be such that $B_0 =: \{i < \omega_1:n_i = n^*\}$ is
uncountable.  By using $(*)_2 \, n$ times we can find an uncountable
$B \subseteq B_0 (\subseteq \omega_1)$ such that
\mr
\item "{$(*)$}"  $a \subseteq B \and a \in J \Rightarrow \text{ max pcf}
\{\lambda^n_i:n < n^*,i \in a\} \le \lambda$.
\ermn
So for every $a \in [B]^{\aleph_0}$, for some $b \in [\omega_1]^{\aleph_0}$
we have $a \subseteq b$ and $h \restriction b \in H \subseteq {\frak B}$.

Let for a set $b \in H(\chi)$ of ordinals, $f_b$ be the $<^*_\chi$-first one-
to-one function from $|b|$ onto $b$, let $g'(i) = f^{-1}_{h(i)}(g(i))$, so
$g'$ is a function from $\theta = \aleph_1$ to $\mu = \aleph_{\alpha(*)}$
(as $|h(i)| \le \mu$).  Now cov$(\aleph_{\alpha(*),\aleph_2,\aleph_2,
\aleph_1}) \le \stick$ so cov$(\aleph_{\alpha(*),\aleph_2,\aleph_2,
\aleph_1}) \le \lambda$ (the only property of $\alpha(*)$ we use) so there
is ${\Cal P}' \subseteq [\aleph_{\alpha(*)}]^{\le \aleph_1},|{\Cal P}'| \le
\lambda$ exemplifying this and \wilog \, ${\Cal P}' \subseteq {\frak B} \and
{\Cal P}' \in {\frak B}$.  So the set $\{g'(i):i \in B\}$ is included in a
countable union of members of ${\Cal P}'$, so for some $Y \in {\Cal P}'$ (so
$Y \in [\aleph_{\alpha(*)}]^{\le \aleph_1},Y \in {\frak B}$) we have $B^* =:
\{i \in B:g'(i) \in Y\}$ is uncountable.

Define $h'$:

$$
\text{Dom}(h') = \omega_1,h'(i) = \{\alpha \in h(i):f_{h(i)}(\alpha) \in Y\}.
$$
\mn
So $h'$ is a function from $\omega_1$ to $[\lambda]^{\le \aleph_1}$ (as $|Y|
\le \aleph_1$) and $i \in B^* \Rightarrow g(i) \in h'(i)$; remember $a \in
[B^*]^{\le \aleph_0} \Rightarrow (\exists b)[a \subseteq b \subseteq \omega_1
\and h \restriction b \in {\frak B}]$.

Let $Z =: \{(i,f^{-1}_{h'(i)}(g(i))):i \in B\}$, it is a subset of $\omega_1
\times \omega_1$ of cardinality $\lambda$, but $\stick = \lambda,\lambda
+1 \subseteq {\frak B}$, so for some infinite $z \in {\frak B},z \subseteq Z$.
Let $z_0 = \{i < \omega_1:\dsize \bigvee_j (i,j) \in z\}$, so $z_0 \in
{\frak B},z_0 \in [\omega_1]^{\aleph_0}$ and even $z_0 \subseteq B^*$, hence
$h' \restriction z_0 \in {\frak B}$.  So as $h' \restriction z_0 \in 
{\frak B}$ and $\{(i,f^{-1}_{h'(i)}(g(i))):i \in z_0\} \in {\frak B}$ also
$g \restriction z_0 \in {\frak B}$, so Rang$(g \restriction z_0) \in B$, so
Rang$(g \restriction z_0) \in {\Cal P}$ and we are done.
\hfill$\square_{\scite{4.7}}$
\enddemo
\bigskip

\definition{\stag{4.9} Definition}  
$$
\align
St^3_{\lambda,\kappa} = \min \bigl\{ |{\Cal P}|:&(a) \quad {\Cal P} \subseteq
[\lambda]^{\aleph_0} \\
  &(b) \quad (\forall A \in [\lambda]^\kappa)(\exists b \in {\Cal P})(b \cap
A \text{ infinite}) \\
  &(c) \quad {\Cal P} \text{ is } AD \text{ which means } a \ne b \in {\Cal P}
\Rightarrow a \cap b \text{ finite} \bigr\}
\endalign
$$
\mn
(main case $\kappa = \aleph_1$).
\enddefinition
\bigskip

\definition{\stag{4.10} Definition}  
$$
\align
St^4_{\lambda,\kappa} = \min \bigl\{ |{\Cal P}|:&(a) \quad {\Cal P} \subseteq
[\lambda]^\omega \\
  &(b) \quad (\forall A \in [\lambda]^\kappa)(\exists b \in {\Cal P})(b \cap
A \text{ infinite}) \\
  &(c) \quad \sup\{\text{otp}(a):a \in {\Cal P}\} < \omega_1 \bigr\}.
\endalign
$$
\enddefinition
\bigskip

\definition{\stag{4.11} Definition}  
$$
\align
St^5_{\lambda,\kappa} = \min \bigl\{ {\Cal P}:&(a) \quad {\Cal P} \subseteq
[\lambda]^\omega \\
  &(b) \quad (\forall A \in [\lambda]^\kappa)(\exists b \in {\Cal P})(b \cap
A \text{ infinite}) \\
  &(c) \quad \text{ the } BA \text{ of subsets of } \lambda \text{ which }
{\Cal P} \text{ and the singletons} \\
  &\qquad \quad \text{ generate is superatomic of rank } 
< \omega_1 \bigr\}.
\endalign
$$
\enddefinition
\bn
\ub{\stag{4.12} Fact}:  dcf$_{\kappa,\aleph_0,\aleph_0} \le St^\ell_{\lambda,
\kappa}$.
\bn
\centerline{$* \qquad * \qquad *$}
\bigskip

\proclaim{\stag{4.13} Claim}  1) Given $\lambda \ge \kappa = \text{ cf}(\kappa)
> \aleph_0$, the following cardinals are equivalent for $k < \omega,k > 0$:
\mr
\item "{$(a)$}"  ecf$_{\kappa,\aleph_0,\aleph_0}(\lambda)$
\sn
\item "{$(b)_k$}"  min$\{|{\Cal F}|:$
{\roster
\itemitem{ $(i)$ }  ${\Cal F}$ is a family of partial functions $f$ from
$\lambda$ to $k$
\sn
\itemitem{ $(ii)$ }  $f \in {\Cal F} \Rightarrow |\text{Dom }f| = \aleph_0$,
\sn
\itemitem{ $(iii)$ }  $f \in {\Cal F},\ell < k \Rightarrow f^{-1}(\{\ell\})$
is infinite
\sn
\itemitem{ $(iv)$ }  if $\langle A_0,\dotsc,A_{k-1} \rangle$ are pairwise
disjoint subsets of $\lambda$ each of cardinality $\kappa$ \ub{then} for some
$f \in {\Cal F}$ we have $\ell < k \Rightarrow f^{-1}(\{\ell\}) \cap A_\ell$
is infinite$\}$
\endroster}
\item "{$(c)_k$}"  like $(b)_k$ replacing (iii), (iv) by
{\roster
\itemitem{ $(iii)^+$ }  if $\langle \alpha_{\varepsilon,\ell}:\varepsilon
< \kappa,\ell < k \rangle$ is a sequence of ordinals, with no repetitions 
then for infinitely many $\varepsilon < \kappa$, for some $\ell < k,
f(\alpha_{\varepsilon,\ell}) = \ell$ \nl
(so $\alpha_{\varepsilon,\ell} \in \text{ Dom}(f)$).
\endroster}
\endroster
\endproclaim
\bigskip

\demo{Proof}  Let $\lambda^b_\kappa,\lambda^c_k$ be the cardinal from
$(b)_k,(c)_k$ respectively and $\lambda^*$ the cardinal from (a).  Clearly 
$\lambda^b_k \le \lambda^b_{k+1},\lambda^c_k \le \lambda^c_{k+1},
\lambda^b_k \le \lambda^c_k,\lambda^c_1 = \lambda^b_1 = \lambda^*$.  So it
suffices to prove $\lambda^c_k \le \lambda^*$, assume ${\Cal P}$ exemplifies
$\lambda^* = \text{ ecf}_{\kappa,\aleph_0,\aleph_0}(\lambda)$ by \scite{4.5}
(2).  If $\lambda^* \ge 2^{\aleph_0}$ let ${\Cal F}^* = \{f:\text{for some }
a \in {\Cal P},f$ is a function from $a$ to $k$ such that $\ell < k
\Rightarrow |f^{-1}(\{\ell\})| = \aleph_0\}$ clearly it exemplifies
$\lambda^c_k \le |{\Cal F}^*| = 2^{\aleph_0} \times \lambda^* = \lambda^*$.
So assume $\lambda^* < 2^{\aleph_0}$, hence $\lambda^* < 2^{\aleph_0}$ and
let $\bar \eta = \langle \eta_i:i < \lambda^* \rangle$ be a sequence of
pairwise distinct members of ${}^\omega 2$.  Let $\bar g = \langle g_\ell:
\ell < k \rangle$ be such that: $g_\ell:\lambda \rightarrow \lambda$ and
$(\forall \alpha_0 \ldots \alpha_{k-1} < \lambda)(\exists \beta < \lambda)
[\dsize \bigwedge_{\ell < k} g_\ell(\beta) = \alpha_\ell]$.  Let ${\Cal P}' =
\{a^{\bar g}:a \in {\Cal P}\}$ where $a^{\bar g} = \{g_\ell(x):\ell < k,x \in
a\}$ and let ${\Cal F} = \{f_{b,h}:b \in {\Cal P}$ and for some $n$ we have
$h \in {}^{{}^n2)}k$ and $\ell < k \Rightarrow |f^{-1}_{b,h}(\{\ell\}) =
\aleph_0\}$ where $f_{b,h}$ is the function with domain $b^{\bar g},
f_{b,h}(i) = h(\eta_i \restriction n)$ where $h \in {}^{({}^n 2)} k$.

Clearly ${\Cal F}$ has the right cardinality and form.  Let us show that it
satisfies the main requirement: let $\langle A_0,\dotsc,A_{k-1} \rangle$ be
a sequence of subsets of $\lambda$ each of cardinality $\kappa$.  Let
$A_\ell = \{\gamma_{\ell,\varepsilon}:\varepsilon < \kappa\}$ (no repetition).
Let $\gamma_\varepsilon < \lambda$ be such that $\dsize \bigwedge_{\ell < k}
g_\ell(\gamma_\varepsilon) = \gamma_{\ell,\varepsilon}$.  For each
$\varepsilon$ for some $n(\varepsilon)$ we have: $\langle
\eta_{\gamma_{\ell,\varepsilon}} \restriction n(\varepsilon):\ell < k \rangle$
is with no repetitions.  As $\kappa = \text{ cf}(\kappa) > \aleph_0$ \wilog
\, for some $\bar \nu = \langle \nu_\ell:\ell < k \rangle$ and $n(*)$ we have
$\varepsilon < \kappa \Rightarrow \eta_{\gamma_{\ell,\varepsilon}}
\restriction n(*) = \nu_\ell,\varepsilon < \kappa \Rightarrow n(\varepsilon)
= n(*)$.  Now by the choice of ${\Cal P}$ for some $a \in {\Cal P},W =
\{\varepsilon:\gamma_\varepsilon \in a\}$ is infinite.  Let $b = a_{\bar g}$,
let $h:{}^{n(*)}2 \rightarrow k$ be such that $h(\nu_\ell) = \ell$, now
$f_{b,h} \in {\Cal F}$ is as required.  \hfill$\square_{\scite{4.13}}$
\enddemo
\bigskip

\proclaim{\stag{4.14} Claim}  Assume $\lambda_{n+1} \ge \text{ ecf}_{\kappa,
\aleph_0,\aleph_0}(\lambda_n,\aleph_0)$ for $n < \omega,\lambda_0 \ge \kappa
= \text{ cf}(\kappa) > \aleph_0$.  \ub{Then} there is a Boolean Algebra $B$ of
cardinality $\dsize \sum_{n < \omega} \lambda_n$ into which the free Boolean
algebra generated by $\dsize \sum_{n < \omega} \lambda_n$ elements can be
embedded \ub{but} such that there is no homomorphism from $B$ onto the free
Boolean Algebra generated by $\kappa$ elements.
\endproclaim
\bigskip

\remark{\stag{4.15} Remark}  1) So we can find quite many pairs $\lambda =
\dsize \sum_{n < \omega} \lambda_n$ and $\kappa$ as required in \scite{4.14}
(or \scite{4.15}).  E.g. any $\lambda = \lambda_n > \beth_\omega$ and large
enough $\kappa = \text{ cf}(\kappa) \le \beth_\omega$ is as required by
\cite{Sh:460}; also $\lambda = \lambda_n = \aleph_{\kappa^{+4}} > \kappa =
\text{ cf}(\kappa) > \aleph_0$ is as required by \cite[\S4]{Sh:400}. \nl
2) On the problem see Fuchino, Shelah, Soukup \cite{FShS:543},
\cite{FShS:544}. \nl
3) From the proof we can strengthen the last phrase in the conclusion to
``no homomorphism from $B$ into $Fr(\kappa)$ with range of cardinality
$\kappa$".  Similarly in \scite{4.16}.
\endremark
\bigskip

\demo{Proof}  Let ${\Cal F}_\alpha = \{f^n_\alpha:\alpha < \lambda_{n+1}\}$
be as required in clause $(b)_2$ of \scite{4.13} (exists as $\lambda_{n+1} \ge
\text{ ecf}_{\kappa,\aleph_0,\aleph_0}(\lambda_n))$.

Remember the variety of Boolean rings has the operation $x \cup y,x \cap y,
x-y$ and constant 0 (but no 1 and no $-x$), so any ideal of a Boolean algebra
is a Boolean ring and if the ideal is maximal, the Boolean algebra is 
definable in the Boolean ring.

Let $B_0$ be the Boolean ring freely generated by $\{x^0_i:i < \lambda_0\}$.
Let $B_1$ be the Boolean ring generated by $B_0 \cup \{x^1_i:i < \lambda_1\}$
freely except:
\mr
\item "{$(a)$}"  the equations which holds in $B_0$
\sn
\item "{$(b)$}"  $x^1_{2 \alpha} \cap x^0_i = 0$ if $f^n_\alpha(i) = 0$
(so $i \in \text{ Dom}(f^n_\alpha))$
\sn
\item "{$(c)$}"  $x^1_{2 \alpha} \cap x^0_i = x^0_i$ if $f^n_\alpha(i) = 1$
(so $i \in \text{ Dom}(f^n_\alpha))$.
\ermn
(Why 
$x^1_{2 \alpha}$ and not $x^1_\alpha$? to embed $Fr(\dsize \sum_{n < \omega}
\lambda_n)$ into the Boolean Algebra we are constructing.)

Similarly let $B_{n+1}$ be the Boolean ring generated by $B_n \cup
\{x^{n+1}_\alpha:\alpha < \lambda_{n+1}\}$ freely except
\mr
\item "{$(a)$}"  the equations which holds in $B_n$
\sn
\item "{$(b)$}" $x^{n+1}_{2 \alpha} \cap x^n_i = 0$ if $f^n_\alpha(i) = 0$
(if $i \in \text{ Dom}(f^n_\alpha))$
\sn
\item "{$(c)$}"  $x^{n+1}_\alpha \cap x^n_i = x^n_i$ if $f^n_\alpha(i) = 1$
(if $i \in \text{ Dom}(f^n_\alpha))$.
\ermn
Now $B_\omega = \dbcu_{n < \omega} B_n$ is a Boolean ring and let $B$ be the
Boolean algebra for which $B_\omega$ is a maximal ideal.  Assume $f$ is a
homomorphism from $B$ onto $Fr(\kappa)$, the Boolean algebra freely generated
say by $\{z_i:i < \kappa\}$.  Now $B$ is generated by $\{x^n_\alpha:n <
\omega,\alpha < \lambda_n\}$.  So as $f$ is onto, for some $n$, for every 
$\zeta < \kappa$ for some $\alpha,f(x^n_{\alpha,\sigma_0}) \notin \langle
z_\varepsilon:\varepsilon < \zeta \rangle_{Fr(\kappa)}$.  By the
$\Delta$-system lemma, we can find a stationary $S \subseteq \kappa$, Boolean
term $\sigma_1 = \sigma_1(x_0,\dotsc,x_{n(*)-1}),m(*) < n(*)$ ordinals
$\varepsilon(0) < \ldots < \varepsilon(m(*)-1) < \min(S)$, and for each
$\zeta \in S$, ordinals $\varepsilon(m(*),\zeta) < \ldots < \varepsilon
(n(*)-1,\zeta)$ all in the interval $[\zeta,\min(S \backslash (\zeta +1))]$
and $\alpha_\zeta < \lambda_n$ such that $f(x^n_{\alpha_\zeta,\sigma}) = 
\sigma_1(z_{\varepsilon(0)},\dotsc,z_{\varepsilon(m(*)-1)},
z_{\varepsilon(m(*),\zeta)},\dotsc,z_{\varepsilon(n(*)-1,\zeta)})$, where all
the $n(*)$ variables are needed in the term $\sigma_1$.

Let $S = \{\zeta(i):i < \kappa\}$ with $\zeta(i)$ increasing in $i$, let $f$
be the function $f(\alpha_{\zeta(2i+1)}) = 1,f(\alpha_{\zeta(2i)})=0$.  So
for some $\alpha < \lambda_{n+1}$ for $\ell = 0,1$ the sets $W_\ell =:
\{\alpha_{\zeta(2i + \ell)}:f^n_\alpha(\alpha_{\zeta(2i + \ell)}) = \ell\}$
are infinite.  So $z^* = f(x^{n+1}_{2 \alpha}) \in Fr(\kappa)$ satisfies
\mr
\item "{$(*)_0$}"  $\alpha_\zeta \in W_0 \Rightarrow Fr(\kappa) \models
z^* \cap f(x^n_{\alpha_\zeta}) = 0$
\sn
\item "{$(*)_1$}"   $\alpha_\zeta \in W_1 \Rightarrow Fr(\kappa) \models
z^* \cap f(x^n_{\alpha_\zeta}) = f(x^n_{\alpha_\zeta})$.
\ermn
But for some finite $u \subseteq \kappa,z^* \in \langle z_\gamma:\gamma \in u
\rangle{Fr(\kappa)}$, so there is $\alpha_{\zeta(2i_0)} \in W_0$, such that
$u$ is disjiont to $\{\varepsilon(m(*),\zeta(2i_0)),\dotsc,\varepsilon(n(*)-1,
\zeta(2i_0 +1))\}$ and there is $\alpha_{\zeta(2i_1+1)} \in W_1$ such that
$u$ is disjoint to $\{\varepsilon(m(*),\zeta(2i_1+1)),\dotsc,\varepsilon(n(*)
-1,\zeta(2i_1 +1))\}$.  For them $(*)_0,(*)_1$ gives a contradiction.  So
$B$ is a Boolean algebra, of cardinality $\le \lambda$ (as it is generated
by $\{x^n_\alpha:\alpha < \lambda_n,n < \omega\})$, we can embed into it the
free Boolean algebra with $\lambda$ generators $\{x^n_{2 \alpha +1}:\alpha <
\lambda_n,n < \omega\})$ with no homomorphism onto the free Boolean algebra
with $\kappa$ generators.  \hfill$\square_{\scite{4.14}}$
\enddemo
\bigskip

\definition{\stag{4.15A} Definition}  Let $B^{fcf}_\mu$ be the Boolean Algebra
of finite and cofinite subsets of $\mu$.
\enddefinition
\bigskip

\proclaim{\stag{4.16} Claim}  Assume $\lambda_{n+1} \ge \text{ ecf}_{\kappa,
\aleph_0,\aleph_0}(\lambda_n,\aleph_0), \kappa = \text{ cf}(\kappa) >
\aleph_0$ and $\lambda = \dsize \sum_n \lambda_n$.  \ub{Then} there is a
Boolean Algebra $B$ of cardinality $\lambda$ into which $B^{fcf}_\lambda$
can be embedded but such that there is no homomorphism from $B$ onto
$B^{fcf}_\kappa$.
\endproclaim
\bigskip

\demo{Proof}  Let

$$
{\Cal P}_n = \{a^n_\alpha:\alpha < \lambda_{n+1}\} \subseteq [\lambda_n]
^{\aleph_0}
$$
\mn
exemplifies $\lambda_{n+1} \ge \text{ ecf}_{\kappa,\aleph_0,\aleph_0}
(\lambda_n,\aleph_0)$ by \scite{4.5}(2).  We define by induction on $n$ a
countable subset $x^n_\alpha$ of $n \times \lambda$ for each $\alpha <
\lambda_n$.  For $n=0$ let $x^0_\alpha = \{(0,\alpha)\}$.  For $n+1$ let

$$
x^{n+1}_{2 \alpha} = \{(n+1,2 \alpha)\} \cup \dbcu_{\beta \in a^n_\alpha}
x^n_\beta \text{ and } x^{2n}_{2 \alpha +1} = \{(n+1,2 \alpha +1)\}.
$$
\mn
Let $B$ be the Boolean Algebra of subsets of $\omega \times \lambda$
generated by $\{x^n_\alpha:\alpha < \lambda_n,n < \omega\}$.

Clearly $|B| \le \dsize \sum_n \lambda_n = \lambda$, also
$\{x^n_{2 \alpha +1}:\alpha < \lambda,n < \omega\}$ generate a subalgebra
isomorphic to $B^{fcf}_\lambda$ hence $|B| \ge \lambda$ (so $|B| = \lambda$)
and $B^{fcf}_\lambda$ can be embedded into $B$.

Lastly, suppose $g$ is a homomorphism from $B$ onto $B^{fcf}_\kappa$.  Let
$z_\zeta =: \{\zeta\} \in B^{fcf}_\kappa$.  For each $\zeta < \kappa$ for
some $n(\zeta) < \omega,\alpha(\zeta) < \lambda_n$ we have
$g(x^{n(\zeta)}_{\alpha(\zeta)}) \notin \langle z_\zeta:\varepsilon < \zeta
\rangle_{B^{fcf}_\kappa}$, so for some stationary $S \subseteq \kappa,
[\zeta \in S \Rightarrow n(\zeta) = n(*)]$.  If $g(x^n_\alpha) \in B^{fcf}
_\kappa$ is infinite then $\{g(x):x \in B \text{ and } g(x) \le g(x^n_\alpha)
\} = \{g(x \cap x^n_\alpha):x \in B\}$ is countable so $g$ is not onto $B$,
a contradiction.  So possibly shrinking $S$ \wilog \, $\langle 
g(x^{n(*)}_{\alpha(\zeta)}):\zeta \in S \rangle$ is a $\Delta$-system of
finite subsets of $\kappa$ with heart called $w$.

For some $\beta < \lambda_{n+1}$ the set $u = \{\zeta \in S:\alpha(\zeta) \in
a^{n+1}_\beta\}$ is infinite, clearly $\zeta \in u \Rightarrow x^n
_{\alpha(\zeta)} \le x^{n+1}_{2 \beta}$ hence $g(x^{n(*)}_{\alpha(\zeta)}) \le
g(a^{n+1}_{2 \beta})$ hence $g(x^{n+1}_{2 \beta})$ is infinite hence it is
co-finite, contradicting an earier statement. \hfill$\square_{\scite{4.16}}$
\enddemo
\bigskip

\definition{\stag{4.18} Definition}  

$$
\align
St^6_{\lambda,\kappa} = \min \bigl\{\lambda^*:&\text{there is } {\Cal P} 
\subseteq [\lambda]^{\aleph_0} \text{ such that} \\
  &(i) \quad {\Cal P} \text{ is } \aleph_2 \text{-free i.e. if } a_i \in
P \text{ for } i < i^* < \aleph_2 \text{ are disjoint} \\
  &\qquad \quad \text{ then for some finite } b_i \subseteq a_i,
\langle a_i \backslash b_i:i < i^* \rangle \\
  &\qquad \quad \text{ are pairwise disjoint } \\
  &(ii) \quad \text{for every } f:\kappa \rightarrow \lambda
\text{ for some } a \in P \\
  &\qquad \qquad \qquad (\exists^\infty \alpha < \kappa)
(f(\alpha) \in a) \bigr\}.
\endalign
$$
\enddefinition
\bigskip

\proclaim{\stag{4.19} Claim}  Assume $\lambda_{n+1} \le 
St^6_{\lambda_{n,\kappa}}$ for $n < \omega,\lambda = \dsize \sum_{n < \omega}
\lambda_n$.  Then there is a Boolean Algebra $B$ as in the previous claim
which is superatomic.
\endproclaim
\bigskip

\demo{Proof}  Like the previous claim.
\enddemo
\newpage

\head {\S5 More on independence} \endhead  \resetall
\bigskip

\proclaim{\stag{5.1} Claim}  If $|\text{\rm pcf\/}({\frak a})| \ge \mu = 
\text{ \rm cf\/}(\mu) > \lambda^+ > |{\frak a}|^+,\sigma < \lambda$ and
$\text{\rm cf\/}(\chi) = \mu,\chi = |\chi \cap \text{\rm pcf\/}({\frak a}|,
\chi$ is $J^{bd}_\lambda$-inaccessible, \ub{then} for any $\chi_0 < \chi$
there is ${\frak b} \subseteq \text{\rm pcf\/}\,\sigma$-complete
$({\frak a}) \cap (\chi_0,\chi) \backslash \lambda,|{\frak b}| < \lambda$,
such that pcf$_{\sigma\text{-complete}}({\frak b}) \cap \chi$ is unbounded
in $\chi$.
\endproclaim
\bigskip

\demo{Proof}  We use \cite{Sh:430}, \scite{6.7}F(5) (localization for
pcf$_{\sigma\text{-complete}}$), \cite{Sh:420},\S1.  
\hfill$\square_{\scite{5.1}}$
\enddemo
\bigskip

\proclaim{\stag{5.2} Claim}  Assume IND$(\langle J^{bd}_{\lambda_\varepsilon}:
\varepsilon < \varepsilon(*) \rangle)$ and $\lambda = \dsize \sum_{\varepsilon
< \varepsilon(*)} \lambda_\varepsilon$.  If $\mu > \lambda$ and $\theta_i \in
\text{ Reg } \cap \mu \backslash \lambda$ for $i < \lambda$, \ub{then} for
some $\varepsilon$ we can find ${\frak c} \subseteq \mu \cap
\text{\rm pcf\/} \{\theta_i:i < \lambda\},|{\frak c}| \le \lambda$ and
${\frak b}_\tau \in J_{\le \tau}[\{\theta_i:i < \lambda\}]$ for $\tau \in
{\frak c}$ such that
\mr
\item "{$(*)$}"  there is no $a \in [\lambda]^{\lambda_\varepsilon}$ such
that $[\tau \in {\frak c} \Rightarrow (\forall^{|{\frak a}|} i \in a)
\theta_i \notin {\frak b}_\tau]$ \nl
$((\forall^{|{\frak a}|} i \in a)$ means: for all but $< |a|$ members $i$
of $a$).
\endroster
\endproclaim
\bigskip

\demo{Proof}  Repeat the proof of \sciteu{3.x}.
\enddemo
\bigskip

\remark{\stag{5.3} Remark}  If we deal with a normal ideal it seems that, for
a given $\chi \in \text{ acc}(c \ell(\text{pcf}({\frak a})))$, we can get long
intervals of $\lambda$ with $\chi$ being $J^{bd}_\lambda$-inaccessible.
\endremark
\bigskip

\proclaim{\stag{5.4} Claim}  1) If IND$(\lambda,\kappa)$, \ub{then} for every
${\frak a} \subseteq \text{ Reg}$, $\lambda \le |{\frak a}| < \min({\frak a})$
for some ${\frak b} \subseteq {\frak a},\kappa < |{\frak b}| \le \lambda,
\Pi({\frak b})/[{\frak b}]^{\le \kappa}$ has true cofinality (so if $\lambda$
is minimal for this $\kappa$, this holds for any $\kappa' < \lambda$). \nl
2) If IND$(\langle J^{bd}_{\lambda_n}:n < \omega \rangle),\lambda_n <
\lambda_{n+1},\lambda_n$ regular, $|{\frak a}| < \lambda_0,|{\frak a}| <
\min({\frak a})$ and

$$
\mu = \underset{n < \omega} {}\to \sup \min\{|{\Cal P}|:{\Cal P} \subseteq
[\lambda_n]^{|{\frak a}|^+} \text{ and } 
(\forall A \in [\lambda_n]^{\lambda_n})
(\exists B \in {\Cal P})[B \subseteq A]\}
$$
\mn
\ub{then} $|\text{pcf}({\frak a}) \le \mu$. \nl
3) Assume $\sigma < \theta < \lambda_n \, (n < \omega)$.  (We can guess
filters: which are $(< \theta)$-based.)  IND$(\langle J_n:n < \omega \rangle),
J_n$ an ideal on $\lambda_n$ and $\mu$ satisfies: $\dsize \bigwedge_n \mu >
\kappa_n$ and
\mr
\item "{$(*)_{\mu,J_n}$}"   there is a set ${\Cal E},|{\Cal D}| \le \mu$, each
member of ${\Cal E}$ is an ideal on some bounded subset of $\kappa_n$ such
that:
{\roster
\itemitem{ $\otimes$ }  \ub{if} $Y \in J^+_n$ (so $Y \subseteq \lambda_n$),
and $I$ is a $(< \theta)$-based $\sigma$-complete ideal on $Y$ generated by
$\le \mu$ sets \ub{then} for some $I' \in {\Cal E}$, we have $(\text{Dom }
I') \cap Y \in I^+$ and Dom$(I') \backslash Y \in I',I' \restriction (Y \cap
\text{ Dom}(I')) \supseteq I \restriction (Y \cap \text{ Dom}(I'))$.
\endroster}
\ermn
\ub{Then} for 
${\frak a} \subseteq \text{ Reg}$, (pcf$_{\sigma\text{-complete}}({\frak a}))
\le \mu$.
\endproclaim
\bigskip

\demo{Proof}  1) Straight. \nl
2) Suppose not.  Let $\langle \tau_{i+1}:i < \mu^+ \rangle$ be the first
$\mu^+$ members of pcf$({\frak a})$ listed in increasing order. Let
$\tau_\delta = \dbcu_{i < \delta} \tau_{1+i}$ for limit $\delta \le \mu^+$.
For each limit $\delta < \mu^+$ for some $n = n_\delta < \omega,\tau_\delta$
is $\{J^{bd}_{\lambda_n}\}$-inaccessible (by \scite{3.13}?).  So for some
$n(*) < \omega,\{\delta < \mu^+:n_\delta = n(*)\}$ is stationary, hence
\mr
\item "{$(*)$}"  for no $\theta_\alpha \in \text{ Reg } \cap \tau_{\mu^+}$
for $\alpha < \lambda_{n(*)}$, do we have
$$
\dsize \prod_{\alpha < \lambda_{n(*)}} \theta_\alpha/J^{bd}_{\lambda_{n(*)}}
\text{ is } \tau_{\mu^+}\text{-directed}.
$$
\ermn
By \cite[\S1]{Sh:420} we can find $\langle C_\alpha:\alpha \in S \rangle,
S \subseteq \mu^+,C_\alpha \subseteq \alpha$, otp$(C_\alpha) \le
\lambda_{n(*)},[\beta \in C_\alpha \Rightarrow \beta \in S \and C_\beta =
C_\alpha \cap \beta]$ and otp$(C_\alpha) = \lambda_{n(*)} \Rightarrow \alpha
= \sup(C_\alpha)$ and $\{\alpha \in S:\text{otp}(C_\alpha) = \lambda_{n(*)}\}$
is stationary.  Now we imitate \cite[\S2]{Sh:400}. \nl
3) Similar to 2).  \hfill$\square_{\scite{5.4}}$
\enddemo
\newpage

\head {\S6 Odds and ends} \endhead  \resetall
\bn
As in \cite[\S6]{Sh:430} this section is dedicated to things I forgot to say.
We repeat and elaborate older things \cite{Sh:430}, \sciteu{6.6D},
\sciteu{6.6E}, \sciteu{6.6F}, \cite{Sh:410}, \scite{3.7}.
\proclaim{\stag{6.1} Claim}  Suppose $D$ is a $\sigma$-complete filter on
$\theta = \text{ cf}(\theta)$ such that $[\alpha < \theta \Rightarrow \theta
\backslash \alpha \in D],\sigma$ is regular $> \kappa^+ + |\alpha|^\kappa$
for $\alpha < \sigma$, and for each $\alpha < \theta,\bar \beta = \langle
\beta^\alpha_\epsilon:\epsilon < \kappa \rangle$ is a sequence of ordinals.
\ub{Then} for every $X \subseteq \theta,X \ne \emptyset$ mod $D$ there is
$\langle \beta^*_\epsilon:\epsilon < \kappa \rangle$ (a sequence of ordinals)
and $w \subseteq \kappa$ such that:
\mr
\item "{$(a)$}"  $\epsilon \in \kappa \backslash w \Rightarrow \sigma \le
\text{ cf}(\beta^*_\epsilon) \le \theta$,
\sn
\item "{$(b)$}"  $B =: \{\alpha \in X:\text{if } \epsilon \in w$ then
$B^\alpha_\epsilon = \beta^*_\epsilon$ and: if $\epsilon \in \kappa
\backslash w$ then $\beta^\alpha_\epsilon$ is $< \beta^*_\epsilon$ but
$> \sup\{\beta^*_\zeta:\zeta < \kappa,\beta^*_\zeta < \beta^*_\epsilon\}\}$
is $\ne \emptyset$ mod $D$
\sn
\item "{$(c)$}"  if $\beta'_\epsilon < \beta^*_\epsilon$ for $\epsilon \in
\kappa \backslash w$ then $\{\alpha \in B:\text{if } \epsilon \in \kappa
\backslash w$ then $B'_\epsilon < \beta^\alpha_\epsilon\} \ne \emptyset$ mod
$D$.
\endroster
\endproclaim
\bigskip

\remark{\stag{6.2} Remark}  1) Of course, we can replace $\kappa$ by any set
of this cardinality. \nl
2) May look at \cite{Sh:620}, \S7, there more is said concerning \scite{6.1}.
\endremark
\bigskip

\demo{Proof}  Let $f_\alpha:\kappa \rightarrow \text{ Ord}$ be $f_\alpha(i)=
\beta^\alpha_i$.

Let $\chi$ be large enough.  We choose by induction on $i < \sigma$, a model
$N_i$ such that:
\medskip

$N_i \prec ({\Cal H}(\chi),\in,<^*_\chi)$;

$\|N_i\| \le 2^\kappa + |i|^\kappa$;

$2^\kappa \subseteq N_0$;

$\kappa,\sigma,\theta \in N_0,\langle f_\alpha:\alpha < \theta \rangle \in
N_0$;

$i < j \Rightarrow N_i \prec N_j$;

$\langle N_j:j \le i \rangle \in N_{i+1}$;

$N_i$ increasing continuous.
\mn

Let $\delta_i =: \min(\cap \{B:B \in N_i \cap D\})$, now $\delta_i$ is well
defined (as $D$ is $\sigma$-complete and $\sigma > 2^\kappa + |i|^\kappa \ge
\|N_i\|$, hence the intersection is in $D$).  Now $\delta_i \ge \sup(\theta
\cap N_i)$.  As $\alpha \in \theta \cap N_i \Rightarrow \theta \backslash
\alpha \in D \cap N_i$ and as $\delta_i \in N_{i+1}$ (as $\{N_i,D,\theta\}
\in N_{i+1})$ clearly $\langle \delta_i:i < \sigma \rangle$ is strictly
increasing continuous.  We define for $i < \sigma$, a function $g_i \in
{}^\kappa \chi$ by

$$
g_i(\zeta) = \min(N_i \cap \chi \backslash f_{\delta_i}(\zeta))
$$
\mn
(it is well defined as $\dbcu_{\alpha < \theta} (f_\alpha(i) + 1) \in N_0
\prec N_i$).  Clearly $E =: \{\alpha < \sigma:N_\alpha \cap \sigma = \alpha\}$
is a club of $\sigma$, and as $(\forall \alpha < \sigma)[|\alpha|^\kappa <
\sigma]$ clearly $\alpha < \beta \in E \and a \subseteq N_\alpha \and
|a| \le \kappa \Rightarrow a \in N_\beta$.  Now $i \in E$, cf$(i) = \kappa^+$
implies $N_i = \dbcu_{j<i} N_j$ and Rang$(g_i) \subseteq \dbcu_{j<i} N_j$
hence $\dsize \bigvee_{j<i}[\text{Rang}(g_i) \subseteq N_j]$; but by the
previous sentence every subset of $N_j$ of cardinality $\le \kappa$ belongs
to $N_i$, hence $g_i \in \dbcu_{j<i} N_j$.  So by Fodor Lemma for some
stationary subset $S$ of $\{i \in E:\text{cf}(i) = \kappa^+\}$ and some
$g^*:\kappa \rightarrow \sigma$ and some $u \subseteq \kappa$ and some $i(*)
< \sigma$ we have: $[i \in S \Rightarrow g_i = g^*],(\forall i \in S)(\forall
\zeta < \kappa)[f_{\delta_i}(\zeta) = g^*(\zeta) \Leftrightarrow \zeta \notin
u]$ and $g^* \in N_{i(*)}$; note $u \in N_0 \subseteq N_{i(*)}$ as $u
\subseteq \kappa$ and we can assume $i(*) < \min(S)$ and $i(*) \in E$.

Let $w =: \kappa \backslash u,\beta^*_i =: g^*(i)$ for $i < \kappa$ now we
show that $w,\langle \beta^*_i:i < \kappa \rangle$ are as required.
\enddemo
\bn
\ub{Clause (b)}:  

The set $B$ is defined from: $\langle B^*_i:i < \kappa \rangle$ and $w$ and
$\bar f = \langle f_\alpha:\alpha < \theta \rangle$. As all of them belong to
$N_{i(*)}$ clearly $B \in N_{i(*)}$, so if $B = \emptyset$ mod $D$ then
$(\theta \backslash B) \in D \cap N_{i(*)}$ hence $\zeta \in S \Rightarrow
\delta_\zeta \in \theta \backslash B \Rightarrow \delta_\zeta \notin B$; but
$\delta_\zeta \in B$ by the definition of $B,g^*,g_\zeta,S$.
\bn
\ub{Clause (a)}:

If $\epsilon \in \kappa \backslash w$, cf$(\beta^*_\epsilon) > \theta$ then
$\gamma^*_\epsilon =: \sup\{f_\alpha(\epsilon):\alpha < \theta,f_\alpha
(\epsilon) < \beta^*_\epsilon\}$ is $< \beta^*_\epsilon$ and it belongs to
$N_{i(*)}$ (as $\epsilon,\langle f_\alpha:\alpha < \lambda \rangle$ and
$\beta^*_\epsilon$ belongs to $N_{i(*)}$) and for any $\zeta \in S$ we get
a contradiction to to $g_\zeta(\epsilon) = \beta^*_\zeta$.

Clearly $\zeta \in \kappa \backslash w \Rightarrow \text{ cf}[g^*(\zeta)] \ge
\sigma$ (otherwise $g^*(\zeta) = \sup(N_i \cap g^*(\zeta))$ (as $N_i \prec
({\Cal H}(\chi),\in,<^*_\chi)$ and $N_{i(*)} \cap \sigma = i(*)$ because
$i(*) \in E$) and easy contradiction.
\bn
\ub{Clause (c)}:  

If there is $\bar \beta' = \langle \beta'_\epsilon \in u \rangle$ 
contradicting clause (c), then there is such a sequence defined from $B,
\langle f_\alpha:\alpha < \theta \rangle,d,w,\langle \beta^*_i:i < \kappa
\rangle$, just use the $<^*_\chi$-first one, hence \wilog \, $\bar \beta' \in
N_{i(*)}$, so for any $\zeta \in S$ we get a contradiction.  
\hfill$\square_{\scite{6.1}}$
\bigskip

\demo{\stag{6.3} Observation}  If $|{\frak a}| < \min({\frak a}),H \subseteq
\Pi {\frak a},|H| = \theta = \text{ cf}(\theta) \notin \text{ pcf}({\frak a})$
and also $\theta > \sup(\theta \cap \text{ pcf}({\frak a}))$ \ub{then} for
some $g \in \Pi {\frak a}$, the set $H_g =: \{f \in H:g < g\}$ has
cardinality $\theta$; in face $H$ is the union of $\le \sup (\theta \cap
\text{ pcf}({\frak a}))$ sets of the form $H_g$. \nl
[Why?  This is as $\Pi {\frak a}/J_{< \theta}[{\frak a}]$ is 
min(pcf${\frak a}) \backslash \theta$)-directed and the ideal
$J_{< \theta}[{\frak a}]$ is generated by $< \theta$ sets.

In details, let $\langle {\frak b}_\sigma[{\frak a}]:\sigma \in
\text{ pcf}({\frak a}) \rangle$ be a generating sequence for ${\frak a}$
(exists by \cite[Ch.VIII]{Sh:g}, \scite{2.6}).

For $\sigma \in \text{ pcf}({\frak a})$ let 
$f^\sigma_\alpha \in \Pi({\frak a})$ for $\alpha < \sigma$ be such that
$\langle f^\sigma_\alpha:\alpha < \sigma \rangle$ is
$<_{J_\theta[{\frak a}]}$-increasing and cofinal and moreover $\{f^\sigma
_\alpha \restriction {\frak b}_\sigma[{\frak a}]:\alpha < \sigma\}$ is cofinal
in $\Pi {\frak b}_\sigma[{\frak a}]$ (where $J_\sigma[{\frak a}] =
J_{< \sigma}[{\frak a}] + {\frak b}_\sigma[{\frak a}])$, exists as $(\Pi
{\frak b}_\sigma,[{\frak a}],<)$ has cofinality $\sigma$ by \cite[Ch.II]{Sh:g},
\scite{3.1}.

Now as $\Pi {\frak a}/J_{< \theta}[{\frak a}]$ is min($\text{pcf}({\frak a})
\backslash \theta^+)$-directed (as $J_{< \theta}[{\frak a}] = 
J_{< \theta^+}[{\frak a}] = J_{< \text{min(pcf}({\frak a}) \backslash
\theta^+)}[{\frak a}]$ and by \cite[Ch.I]{Sh:g}, \scite{1.5}) there is $g \in
\Pi {\frak a}$ such that: $h \in H \Rightarrow h > g \text{ mod }
J_{< \theta}[{\frak a}]$; hence for each $h \in F$ for some finite $\Theta(h)
\subseteq \theta \cap \text{ pcf}({\frak a})$ we have

$$
\{\sigma \in {\frak a}:h(\sigma) \ge g(\sigma)\} \subseteq \bigcup
\{{\frak b}_\sigma[{\frak a}]:\sigma \in \Theta(h)\}.
$$
\mn
Also for every $\sigma \in \text{ pcf}({\frak a})$ we can find $\alpha_\sigma
(h)$ (for $\sigma \in \text{ pcf}({\frak a}),h \in H)$ such that $h
\restriction {\frak b}_\sigma[{\frak a}] < f^\sigma_{\alpha_\sigma(h)}
\restriction {\frak b}_\sigma[{\frak a}]$.  So $h < \max(\{g,f^\sigma
_{\alpha_\sigma(h)}:\sigma \in \Theta(h)\})$.  Let

$$
\align
G =: \bigl\{ \max\{g,f^{\sigma_1}_{\alpha_1},\dotsc,f^{\sigma_n}
_{\alpha_n}\}:&n < \omega,\{\sigma_1,\dotsc,\sigma_n\} \subseteq \theta \cap
\text{ pcf}({\frak a}) \\
  &\text{and } \alpha_1 < \sigma_1,\dotsc,\alpha_n < \sigma_n \bigr\}
\endalign
$$
\mn
it has cardinality $\le \aleph_0 + \sup(\theta \cap \text{ pcf}({\frak a})
< \theta$ and $\theta = \text{ cf}(\theta) = |H|$ and $\forall h \in H
\exists g' \in G(h < g)$. \nl
So for some $g^* \in G$ the set $\{h \in H:h < g^*\}$ has cardinality
$\theta$ as required.]  \hfill$\square_{\scite{6.3}}$
\enddemo
\bn
We comment to \cite{Sh:410}, \scite{3.7} (which solve a problem from Gerlits
Hajnal Szentmiklossy \cite{GHS}).
\proclaim{\stag{6.4} Claim}  1) Suppose $\gamma^*$ and $i^* = i(*)$ are
ordinals and $\bar \chi = \langle \chi_i:i < i^* \rangle$ is a sequence of
infinite cardinals; so of course we can find $n,\bar w = \langle w_\ell:
\ell < n \rangle,\bar \kappa = \langle \kappa_\ell:\ell < n \rangle,
\bar \sigma = \langle \sigma_\ell:\ell < n \rangle$ such that $\bar w$ is a
partition of $i^*,|w_\ell| = \kappa_\ell,\dsize \bigwedge_{i \in w_\ell}
\chi_i \le \sigma_\ell$, and

$$
(\forall \chi < \sigma_\ell)(\exists^{\kappa_\ell} i)(i \in w_\ell \and \chi
< \chi_i \le \sigma_\ell)
$$
\mn
and $\bar \kappa$ is strictly increasing, $\bar \sigma$ is strictly
decreasing, in fact $\bar w, \bar \kappa, \bar \sigma$ are unique.  \ub{Then}
the following are equivalent
\mr
\item "{$(A)_{\bar \chi,\gamma^*}$}"  there are $f_\alpha \in \dsize \prod
_{i < i(*)} \chi_i$ for $\alpha < \gamma^*$ satisfying $\alpha < \beta
\Rightarrow (\exists i < i^*)[f_\alpha(i) > f_\beta(i)]$
\sn
\item "{$(B)_{\bar \chi,\gamma^*}$}"  for some $\ell < n$ we have:
$2^{\kappa_\ell} \ge |\gamma^*|$ \ub{or} for every regular $\mu_1 \in
(\gamma^* +1) \backslash \sigma^+_\ell$ for some singular $\lambda^* \le
\sigma_\ell$ we have
{\roster
\itemitem{ $(*)$ }  cf$(\lambda^*) \le \kappa_\ell,\lambda^* > 
2^{\kappa_\ell}$ and pp$^+(\lambda^*) > \mu_1$
\endroster}
\ermn
2) If
\mr
\item "{$\otimes$}"  $(\forall \mu_1)(\mu_1 = \text{ cf}(\mu_1) \le
|\gamma^*| \rightarrow (\exists \mu_2)[\mu_2 = \text{ cf}(\mu_2) \and \mu_1
\le \mu_2 \le |\gamma^*| \and (\forall \alpha < \mu_2)|\alpha|^{\aleph_0}
< \mu_2])$,
\ermn
\ub{then} in part (1), $(B)_{\bar \chi,\gamma^*}$ the demand $(*)$ on
$\lambda^*$ can be replaced by
\mr
\item "{$(*)'$}"  cf$(\lambda^*) \le \kappa_\ell,\lambda^* > 2^{\kappa_\ell}$
and $(\forall \mu < \lambda^*)(\mu^{\kappa_\ell} < \lambda^*)$.
\ermn
Now we call it $(B)'_{\bar \chi,\gamma^*}$ (so if $\otimes$ then
$(A)_{\bar \chi,\gamma^*} \Leftrightarrow (B)_{\bar \chi,\gamma^*}
\Leftrightarrow (B)'_{\bar \chi,\gamma^*}$). \nl
3) If $\otimes$ from part (2) holds, \ub{then} also $(A)_{\bar \chi,\gamma^*}
\Leftrightarrow (B)''_{\bar \chi,\gamma^*} \Leftrightarrow (B)^+_{\bar \chi,
\gamma^*}$ where
\mr
\item "{$(B)''_{\bar chi,\gamma^*}$}"  $|\gamma^*| \le
\underset{\ell < n} {}\to \max(\sigma_\ell)^{\kappa_\ell}$
\sn
\item "{$(B)^*$}"  for some $\ell < n$ we have: $2^{\kappa_\ell} \ge
|\gamma^*|$ \ub{or} for every regular $\mu_1 \in (\gamma^* + 1) \backslash
\sigma^+_\ell$ for some singular $\lambda^*$ we have
{\roster
\itemitem{ $(*)^+$ }  cf$(\lambda^*) \le \kappa_\ell,\lambda^* >
2^{\kappa_\ell}$ and $(\forall \mu < \lambda^*)(\mu^{\kappa_\ell} 
< \lambda^*)$ and pp$^+(\lambda^*) > \mu_1$.
\endroster}
\ermn
4) In part (2), (3) instead $\otimes$ we may let $\lambda_0 = \max\{
2^{\kappa_\ell}:\ell < n\},\lambda_1 = |\gamma^*|$, demand
\mr
\item "{$\oplus_{\lambda_0,\lambda_1}$}"  if $\lambda_0 < \mu \le \lambda_1$,
cf$(\mu) = \aleph_0$ and $(\forall \lambda < \mu)(|\lambda|^{\aleph_0} <
\mu)$ then pp$(\mu) =^+ \mu^{\aleph_0}$.
\endroster
\endproclaim
\bigskip

\remark{Remark}  1) On $\oplus_{\lambda_0,\lambda_1}$ from \scite{6.4}(4),
see \cite{Sh:430},\S1. \nl
2) Note that we could in \scite{6.4}(1) demand

$$
(\forall \chi < \sigma_\ell)(\exists^{\kappa_\ell}i)(i \in {}^w \ell \and \chi
\le \chi_1 < \sigma_\ell)
$$
\mn
and can allow $\chi_i$ (hence $\sigma_i$) to be any ordinal, and even let 
$\kappa_i$ be ordinal so the demand is $(\forall \alpha < \sigma_\ell)
[\kappa_\ell = \text{ otp}\{i \in w_\ell:\alpha \le \chi_i < \sigma_\ell\}]$.
This causes no serious change.
\endremark
\bigskip

\demo{Proof}  1) 
$(B)_{\bar \chi,\gamma^*} \Rightarrow (A)_{\bar \chi,\gamma^*}$.

Let $\ell < n$ exemplifies $(B)_{\bar \chi,\gamma^*}$ so there are $f'_\alpha
\in {}^{\kappa_\ell)}(\sigma_\ell)$ for $\alpha < \gamma^*$ such that
$\alpha < \beta < \gamma^* \Rightarrow (\exists j < \kappa_\ell)[f'_\alpha(j)
> g'_\beta(j)]$. \nl
[Why?  Easy by cases.]

Let $h:w_\ell \rightarrow \kappa_\ell$ be such that:

$$
(\forall j < \kappa_\ell)(\forall \sigma < \sigma_\ell)(\exists^{\kappa_\ell}
i)[i \in w_\ell \and \sigma < \chi_i \le \sigma_\ell \and h(i) = j].
$$
\mn
Let $f_\alpha \in \dsize \prod_{i < i(*)} \chi_i$ be: $f_\alpha(i) =
f'_\alpha(h(i))$ if $i \in w_\ell,f'_\alpha(h(i)) < \chi_i$ and $f_\alpha(i)
=0$ otherwise; so if $\alpha < \beta < \gamma^*$ then for some $j < \kappa
_\ell,f'_\alpha(j) < f'_\beta(j) < \sigma_\ell$ so for some $i \in w_\ell$ we
have: $h(i) = j \and \chi_i > f'_\beta(j)$.  So $f_\alpha(i) = f'_\alpha(j)
< f'_\beta(j) < f_\beta(i)$ as required.
\bn
$(A)_{\bar \chi,\gamma^*} \Rightarrow (B)_{\bar \chi,\gamma^*}$

Assume this fails, note that clearly $\gamma_1 < \gamma_2 \and
(A)_{\bar \chi,\gamma_2} \Rightarrow (A)_{\bar \chi,\gamma_1}$ and $\gamma_1
< \gamma_2 \and (B)_{\bar \chi,\gamma_2} \Rightarrow 
(B)_{\bar \chi,\gamma_1}$.Without loss of generality 
$\gamma^*$ is minimal (for our $\bar chi$) for which this fails; so $\gamma^*$
is minimal such that $(B)_{\bar \chi,\gamma^*}$ fails inspecting
$(B)_{\bar \chi,\gamma^*}$ as $n$ is finite, clearly $\gamma^*$ is a regular
cardinal, call it $\theta$.  By renaming we can assume that $i^*$ is a
cardinal.  Now let $\langle f_\alpha:\alpha < \theta \rangle$ exemplifies
$(A)_{\bar \chi,\theta}$; now apply \scite{6.1} above with $i^*,\langle
f_\alpha:\alpha < \theta \rangle,D^{cb}_\theta$ the filter of cobounded
subsets of $\theta,\underset{\ell < n} {}\to \max(2^{\kappa_\ell})^+ =
(2^{|i^*|})^+$ here standing for $\kappa,\left< \langle \beta^\alpha_i:i <
\kappa \rangle:\alpha < \theta \right>,D,\sigma$ there.

So we get $w,\langle \beta^*_i:i < i^* \rangle,B$ as there, and let
${\frak a} = \{\text{cf}(\beta^*_i):i \in i^* \backslash w\}$, so ${\frak a}
\subseteq \text{ Reg } \cap (\underset{\ell < n} {}\to \max \sigma_\ell)
\backslash (2^{|i^*|})^+$ (see clause (a) of \scite{6.1}) and $|{\frak a}| \le
|i^*|$.  Now if $\theta \le \text{ max pcf}({\frak a})$ then for some $\ell,
\theta \le \text{ max pcf}\{\text{cf}(\beta^*_i):i \in w_\ell \backslash w\}$,
and so $(B)_{\bar \chi,\theta}$ does not fail, contradicting an earlier
assumption.  So $\theta > \text{ max pcf}({\frak a})$, so there is a cofinal
$H \subseteq \dsize \prod_{i \in (i^* \backslash w)} \beta^*_i$ of cofinality
$< \theta$, so there are $h_\alpha \in H$ such that $f_\alpha \restriction
(i^* \backslash w) < h_\alpha$ but $|B| = \theta > |H|$ (by the choice of
$D^{cb}_\theta$ as $D$) so for some $h^* \in H$ the set $B_1 = \{\alpha \in B:
h_\alpha = h^*\}$ is unbounded in $\theta$.  By clause (c) of the conclusion
of \scite{6.1} for some $\alpha \in B$ we have $i \in i^* \backslash w
\Rightarrow h^*(i) < f_\alpha(i)$.  Choose $\beta \in B_1 \backslash (\alpha
+1)$, so $\alpha < \beta$ are in $B$, hence $f_\alpha \restriction w =
f_\beta \restriction w$ and $i \in i^* \backslash w \Rightarrow f_\beta(i) <
h^*(i) < f_\alpha(i)$, so $\dsize \bigwedge_i f_\beta(i) \le f_\alpha(i)$,
a contradiction to $\langle f_\gamma:\gamma < \theta \rangle$ exemplifies
$(A)_{\bar \chi,\theta}$. \nl
2), 3)  Clearly $(B)^+_{\bar \chi,\gamma^*} \Rightarrow 
(B)_{\bar \chi,\gamma^*} \Rightarrow (B)''_{\bar \chi,\gamma^*}$ by checking.
So we should just prove:
\bn
$(B)''_{\bar \chi,\gamma^*} \Rightarrow (B)'_{\bar \chi,\gamma^*}$

We can assume $|\gamma^*| > 2^{\kappa_\ell}$ for $\ell < n$ (otherwise the
conclusion is trivial).  We know by $(B)''_{\bar \chi,\gamma^*}$ that for
some $\ell,(\sigma_\ell)^{\kappa_\ell} \ge |\gamma^*|$.  Let us check now
$(B)'_{\bar \chi,\gamma^*}$ so let a regular $\mu_1 \in (\gamma^* +1)
\backslash \sigma^+_\ell$ be given.  So $(\sigma_\ell)^{\kappa_\ell} \ge
|\gamma^*| \ge \mu_1$ hence $\lambda =: \min\{\lambda:\lambda^{\kappa_\ell}
\ge \mu_1\}$ is $\le \sigma_\ell$ but is $> 2^{\kappa_\ell}$, hence it is
singular, cf$(\lambda) \le \kappa_\ell$ and $(\forall \alpha < \lambda)
(|\alpha|^{\kappa_\ell} < \lambda)$, i.e. as required in $(*)'$.
\bn
$(B)'_{\bar \chi,\gamma^*} \Rightarrow (B)^+_{\bar \chi,\gamma^*}$

Again we can assume $|\gamma^*| > 2^{\kappa_\ell}$ for $\ell < n$.  Let us
check $(B)_{\bar \chi,\gamma^*}$, so let a regular $\mu_1 \in (\gamma^* +1)
\backslash \sigma^+_\ell$ be given.  As we are assuming $\otimes$ from
\scite{6.4}(2) there is $\mu_2 \in \text{ Reg } \cap (\gamma^* +1) 
\backslash \mu_1$ such that $(\forall \alpha < \mu_2)(|\alpha|^{\aleph_0} <
\mu_2)$.  Apply $(B)'_{\chi,\gamma^*}$ for $\mu_2$ and get $\lambda^*_2$ as in
$(*)$ for $\mu_2$ instead of $\mu_1$.  Clearly $(\lambda^*_2)^{\kappa_\ell}
\ge \mu_2$ and let $\lambda^* = \min\{\lambda:\lambda^{\kappa_\ell} \ge
\mu_2\}$, so $\lambda^* \le \lambda^*_2 \le \sigma_\ell,\lambda^* >
2^{\kappa_\ell}$ and $(\forall \mu < \lambda^*)[\mu^{\kappa_\ell} 
< \lambda^*]$, and clearly cf$(\lambda^*) \le \kappa_\ell$, so
$(\lambda^*)^{\kappa_\ell} = (\lambda^*)^{\text{cf}(\lambda)}$.

By the choice of $\mu_2$ necessarily cf$(\lambda^*) > \aleph_0$ (otherwise
$\mu_2 \le (\lambda^*)^{\kappa_\ell} = (\lambda^*)^{\aleph_0} < \mu_2)$.
By \cite{Sh:g}, (see \scite{6.5} below), pp$(\lambda^*) =^+ (\lambda^*)
^{\text{cf}(\lambda^*)} = (\lambda^*)^\kappa$ as required in $(*)^+$. \nl
4) The only place we use was in choosing $\mu_2$ in the proof of
$(B)'_{\bar \chi,\gamma^*} \Rightarrow (B)^+_{\bar \chi,\gamma^*}$ and the
use of its property is to show cf$(\lambda^*) > \aleph_0$ (to be able to use
\cite[Ch.VIII]{Sh:g},\S1) but we can use instead 
$\oplus_{\lambda_0,\lambda_1}$.   \hfill$\square_{\scite{6.4}}$
\enddemo
\bn
Remember that by \cite{Sh:g}:
\demo{\stag{6.5} Observation}  If $\mu$ is singular, cf$(\mu) > \aleph_0$
and $\alpha < \mu \Rightarrow |\alpha|^{\text{cf}(\mu)} < \mu$ then
$\mu^{\text{cf}(\mu)} = \text{ pp}(\mu)$.
\enddemo
\bigskip

\demo{Proof}  By \cite{Sh:430}, \scite{3.5}, \cite[CH.VIII]{Sh:g}, \scite{1.8},
\cite[Ch.II]{Sh:g}, \sciteu{5.6}.  \hfill$\square_{\scite{6.5}}$
\enddemo
\bigskip

\definition{\stag{6.6} Definition}  1) For $F \subseteq {}^\delta\text{Ord}$,
we say $F$ is free$^\ell$ when we can find $\zeta_f < \delta$ for $f \in F$
such that:
\mr
\item "{$(a)$}"  if $\ell=1$ then
$$
f \ne g \in F,\zeta = \max\{\zeta_f,\zeta_g\} \Rightarrow f \restriction \zeta
\ne g \restriction \zeta
$$
\sn
\item "{$(b)$}"  if $\ell=2$ then
$$
f \ne g \in F,\delta > \zeta \ge \max\{\zeta_f,\zeta_g\} \Rightarrow 
f(\zeta) \ne g (\zeta)
$$
\sn
\item "{$(c)$}"  if $\ell=3$ then
$$
f \ne g \in F,\delta > \zeta > \varepsilon = \max\{\zeta_f,\zeta_g\},
f(\varepsilon) \le g(\varepsilon) \Rightarrow f(\zeta) \le g(\zeta)
$$
\sn
\item "{$(d)$}"  if $\ell=4$ then for $f,g \in F,f \le g$ (i.e. \nl
$\zeta < \delta \Rightarrow f(\zeta) \le g(\zeta))$ or $g \le f$
\sn
\item "{$(e)$}"  if $\ell=5$ then for some $\zeta$ and $h$ we have
$$
f \in F \Rightarrow f \restriction \zeta = h
$$

$$
\{f(\zeta):f \in F\} \text{ is with no repetition}.
$$
\ermn
2)  Let free$^{\ell,m}$ means free$^\ell$ and free$^m$.  For $J$ an ideal on
$\delta$ we write $J$-free$^\ell$ if $\zeta_f < \delta$ is replaced by
$s_f \in J$.
\enddefinition
\bigskip

\definition{\stag{6.7} Definition}  For $F \subseteq {}^\delta\text{Ord}$: \nl
1) We say $F$ is $\mu$-free$^x$ if every $F' \in [F]^\mu$ is free$^x$. \nl
2) We say $F$ is $(\mu,\kappa)$-free$^x$ if every $F' \in [F]^\mu$ there is 
$F'' \in [F']^\kappa$ which is free$^x$. \nl
\enddefinition
\bn
\ub{\stag{6.8} Fact}:  1) ``$F$ is free$^\ell$ implies $F$ is free$^m$" when
$(\ell,m)$ is one of (2,1), (4,3), (5,1). \nl
2)  Similarly for $\mu$-free$^x$ and $(\mu,\kappa)$-free$^x$.
\bigskip

\demo{Proof}  Straight. 

On \scite{6.9} see \cite[Ch.II]{Sh:g},\S1, \cite[Ch.II]{Sh:g},\S3,
\cite{Sh:282}, \cite[Ch.II]{Sh:g}, \scite{4.10}, Shelah Zapletal
\cite{ShZa:561}.
\enddemo
\bigskip

\proclaim{\stag{6.9} Claim}  1) If $J^{bd}_\delta \subseteq J,J$ an ideal on
$\delta,\dsize \prod_{i < \delta} \lambda_i/J$ is $\lambda^+$-directed,
$\langle \lambda_i:i < \delta \rangle$ an increasing sequence of regulars
$> \delta$ with limit $\mu,\mu < \lambda = \text{ cf}(\lambda)$, \ub{then}
there are regulars $\lambda'_i < \lambda_i$ with tlim$(\lambda'_i) = \mu$
and $f_\alpha \in \Pi \lambda'_i$ for $\alpha < \lambda$ such that $\langle
f_\alpha:\alpha < \lambda \rangle$ is $<_J$-increasing cofinal in
$\dsize \prod_{i < \delta} \lambda'_i$ and $\{f_\alpha:\alpha < \lambda\}$
is $\mu^+-J$-free. \nl
2) Assume ${\frak a}$ is a set of regular $> |{\frak a}|$ with no last
element, $J$ an ideal on ${\frak a}$ extending $J^{bd}_{\frak a}$ and
${\frak c} = \{\theta \in {\frak a}:\theta > \text{ max pcf}(\theta \cap
{\frak a})\}$ and $\lambda = \text{ max pcf}({\frak a})$.  \ub{Then} there is
$\langle f_\alpha:\alpha < \lambda \rangle$ cofinal in $\Pi {\frak a},
<_J$-increasing, such that:
\mr
\item "{$(*)$}"  if $\theta \in {\frak c}$ then $\{f_\alpha \restriction
(\theta \cap {\frak a}):\alpha < \lambda\}$ has cardinality $< \theta$.
\ermn
3) If $\mu > \mu_0 \ge \kappa \ge \text{ cf}(\mu),\lambda = \mu^+$ (or just
$\mu < \lambda = \text{ cf}(\lambda) < \text{ pp}^+_\kappa(\mu)$ \ub{then}
for some ${\frak a} \subseteq (\mu_0,\mu) \cap \text{ Reg},|{\frak a}| \le
\kappa,[\theta \in {\frak a} \Rightarrow \text{ max pcf}(\theta \cap 
{\frak a}) < \theta]$ and $\lambda = \text{ max pcf}({\frak a})$ (if
$[\alpha < \mu \Rightarrow |\alpha|^\kappa < \mu]$ we can have $\mu = \sup
({\frak a})$, otp$({\frak a}) = \text{ cf}(\mu)$ (so part (2) is not empty)).
\nl
4) If $J$ is an ideal on ${\frak a},\langle f_\alpha:\alpha < \lambda \rangle$
is $<_J$-increasing $<_J$-cofinal in $\Pi {\frak a}$ and $J$ is generated by
$< \text{ min}({\frak a})$ sets (as an ideal) \ub{then} every $A \in 
[\lambda]^\lambda$ for some ${\frak d} \in J$ we have: for every $g \in \Pi
{\frak a}$ for $\lambda$ ordinals $\alpha \in A,g \restriction ({\frak a}
\backslash {\frak d}) < f_\alpha \restriction ({\frak a} \backslash 
{\frak d})$.  Hence $\bar f$ is $(\lambda,\min({\frak a})-J$-free
$^{\{2,3\}}$. \nl
5) Assuming $|{\frak a}| < \min({\frak a}),\lambda = 
\text{ tcf}(\Pi{\frak a},<_{J^{bd}_{\frak a}}),{\frak c} = \{\theta \in
{\frak a}:\text{max pcf}(\theta \cap {\frak a})\}$ is unbounded in ${\frak a}$
and $\langle f_\alpha:\alpha < \lambda \rangle$ as in part 2.  \ub{Then} not
only for each $\theta \in {\frak a}$ is $\{f_\alpha:\alpha < \lambda\},
(\lambda,\theta)$-free$^{\{2,3\}}$, but for any $A \in [\lambda]^\lambda$ for
every large enough $\theta \in {\frak c}$ there is $B \in [A]^\theta$ such
that $\langle f_\alpha \restriction (\theta \cap {\frak a}):\alpha \in B
\rangle$ is constant and $\langle f_\alpha \restriction ({\frak a} \backslash
\theta):\alpha \in B \rangle$ is strictly increasing. \nl
6) Assume $\lambda = \text{ tcf}(\Pi{\frak a},<_J),J^{bd}_{\frak a} \subseteq
J,\lambda > \mu = \sup({\frak a})$ and $\langle f_\alpha:\alpha \in \lambda
\rangle$ is $<_J$-increasing and $<_J$-cofinal
\mr
\item "{$(a)$}"  if $\kappa \le \mu^+$ and \footnote{note: if $\lambda$ is
a successor of regular $> \kappa^+$ then this holds} $\{\delta < \lambda:
\text{cf}(\delta) < \kappa\} \in I[\lambda]$, \ub{then} for some $A \in
[\lambda]^\lambda$ we have $\{f_\alpha:\alpha \in A\}$ is $\mu$-free$^1$,
also we can find $\bar f = \langle f_\alpha:\alpha < \lambda \rangle$ as
above with $A=\lambda$, such that $\bar f$ is ${}^b$ continuous (see
\cite[Ch.I,\S3]{Sh:g})
\sn
\item "{$(b)$}"  if $\kappa = \text{ cf}(\kappa) < \mu^+$ and $\{\delta <
\lambda:\text{cf}(\delta) < \kappa\}$ we get similar results for
$(\kappa,\kappa)$-free$^1$.
\ermn
7) Assume $\langle \mu_i:i \le \kappa \rangle$ is an increasing continuous
sequence of singulars $> \kappa,\kappa = \text{cf}(\kappa) > \aleph_0,
\theta \in \text{ Reg } \cap \mu_0 \backslash \kappa^+$ and
\footnote{note: if $\lambda$ is a successor of a regular $> \kappa^+$ then
this holds} $\delta < \mu^+_i:\text{cf}(\delta) = \theta\} \in I[\mu^+_i]$
for $i \le \kappa$ and pp$_\kappa(\mu_i) < \mu^+_{i+1}$ for each $i \le
\kappa$.  \ub{Then} $<_{J^{bd}_C}$-increasing and cofinal and is
$(\theta,\theta)$-free$^1$.  If we demand in addition $\{\delta < \mu^+_i:
\text{cf}(\delta) \le \theta\} \in I[\mu^+_i]$ then $F$ is $\theta^+$-free
$^1$.
\endproclaim
\bigskip

\demo{Proof}  1) By \cite[Ch.II]{Sh:g},\S1. \nl
2) By \cite[Ch.II]{Sh:g}, \scite{3.5}. \nl
3) By \cite[Ch.II]{Sh:g},\S3. \nl
4) Can prove as in \cite{Sh:282}.  Or as in \scite{6.1}, as below.

Let $A \in [\lambda]^\lambda$ and let $\{{\frak d}_\zeta:\zeta < \zeta^*\}
\subseteq J$ be a family generating $J$ closed under finite union such that
$\zeta^* < \lambda_0$.  We shall prove that for some $\zeta < \zeta^*$ we
have $(\forall g \in \Pi{\frak a})(\exists^\lambda \alpha \in A)(g
\restriction ({\frak a} \backslash {\frak d}_\zeta) < f_\alpha \restriction
({\frak a} \backslash {\frak d}_\zeta))$.  If not, for each $\zeta < \zeta^*$
some $g_\zeta \in \Pi{\frak a}$ and $\alpha_\zeta < \lambda$ exemplifies it,
i.e. $\alpha \in A \backslash \alpha_\zeta \Rightarrow \neg(g \restriction
({\frak a} \backslash {\frak d}_\zeta) < f_\alpha \restriction ({\frak a}
\backslash {\frak d}_\zeta))$.  So defined by the function $g \in 
\Pi{\frak a},g(\theta) =: \sup[\{g_\zeta(\theta):\zeta < \zeta^*\} \cup
\{f_{\alpha^*}(\theta)+1\}]$ is well defined and let $\alpha^* = \sup
\{\alpha_\zeta:\zeta < \zeta^*\} < \lambda$.  Now $\{f_\alpha:\alpha <
\lambda\}$ is $<_J$-increasing and $<_J$-cofinal in $\Pi{\frak a}$ so for some
$\alpha \in A \backslash (\alpha^* +1),g <_J f_\alpha$, so for some $\zeta <
\zeta^*$ we have $\{\theta \in {\frak a}:\neg(g(\theta) < f_\alpha(\theta))
\} \subseteq {\frak d}_\zeta$, so for ${\frak d}_\zeta$, the pair $(g_\zeta,
\alpha_\zeta)$ is not as required, a contradiction. \nl
5), 6, 7)  Left to the reader.  \hfill$\square_{\scite{6.9}}$ 
\enddemo
\bigskip

\definition{\stag{6.10} Definition}  We say $(\Pi {\frak a},<_J)$ is
$X$-free$^y$ if there is a $<_J$-increasing $<_J$-cofinal $\langle f_\alpha:
\alpha < \lambda \rangle$ from $\Pi {\frak a}$ which is $y$-free$^x$.
\enddefinition
\bn
\ub{\stag{6.11} Fact}:  If $y \in \{1,2,3,\{2,3\}\},x \in \{\mu-J,(\mu,\theta)
-J\}$ and $(\Pi{\frak a},<_J)$ is $x$-free$^y$ and $\langle f_\alpha:\alpha <
\lambda \rangle$ is $<_J$-increasing $<_J$-cofinal in $\Pi{\frak a}$ \ub{then}
for some $A \in [\lambda]^\lambda$ we have $\{f_\alpha:\alpha \in A\}$ is
$x$-free$^y$.
\bigskip

\demo{Proof}  Straight.
\enddemo
\bn
\ub{\stag{6.12} Question}:  Let $\mu > \kappa \ge \text{ cf}(\mu)$.  For how
many ${\frak a} \subseteq \text{ Reg } \cap \mu,\mu = \sup({\frak a}),
J^{bd}_{\frak a} \subseteq J,J$ an ideal on ${\frak a},\lambda =
\text{ tcf}(\Pi{\frak a},<_J)$, is $(\Pi{\frak a},<_J)$ not $\mu$-free?
\bn
The proof of \cite[Ch.IX]{Sh:g}, \scite{3.5} has a gap (in the reference to
\cite[Ch.IX]{Sh:g}, \sciteu{3.3A}).  What we know is only
\proclaim{\stag{6.13} Lemma}  1) Assume cf$(\lambda) \in [\sigma,\theta),
\lambda > \theta > \text{ cf}(\lambda) \ge \sigma = \text{ cf}(\sigma) >
\aleph_0$.  \ub{Then} cov$(\lambda,\lambda,\theta,\sigma) =^+ \sup
\{\text{pcf}_{\Gamma(\theta,\sigma),J^{bd}_{\frak a}}({\frak a}):{\frak a}
\subseteq \text{ Reg } \cap \lambda,\lambda = \sup({\frak a}),|{\frak a}|
< \min({\frak a})\}$ where

$$
\align
\text{pcf}_{\Gamma(\theta,\sigma),J}({\frak a}) = \bigl\{ \text{tcf}
(\Pi{\frak a},<_I):&I \text{ an ideal on } {\frak a} \text{ extending } J,
\text{tcf}(\Pi{\frak a},<_I) \text{ well defined}, \\
  &I \text{ is } \sigma \text{-complete and for some } {\frak b} \in I,
|{\frak a} \backslash {\frak b}| < \theta \bigr\}
\endalign
$$
\mn
(the $=^+$ means that if the left side is regular then the supremum in the
right side is obtained). \nl
2) If in addition $(\forall \mu < \lambda)(\text{cov}(\mu,\theta,\theta,
\sigma) < \lambda),(\theta,\sigma$ regular) \ub{then}
$$
\text{cov}(\lambda,\lambda,\theta,\sigma) =^+ \text{ pp}_{\Gamma(\theta,
\sigma)}(\lambda).
$$
\mn
3)  So for ${\Cal Y} = {\Cal Y}_\mu,Eq = Eq_\mu$ be as in \relax
\cite[\S3- \S5]{Sh:420}, cf$(\mu) = \sigma > \aleph_0$ for simplicity, $\mu > \theta >
\text{ cf}(\mu)$. \nl
If ${\frak a} \subseteq \text{ Reg } \backslash \mu,|{\frak a}| < \mu,
|{\frak a}| < \min({\frak a}),\lambda$ inaccessible, $J = J^{bd}_{\frak a},
\lambda = \text{ sup pcf}_{\Gamma(\theta,\sigma),J^{bd}}({\frak a})$, then
we can find $e \in Eq,\bar \lambda = \langle \lambda_x:x \in {\Cal Y}_\mu/e
\rangle$ and $D \in FIL({\Cal Y}_\mu)$ such that:
\mr
\item "{{}}"  $\lambda = \text{ tcf}(\Pi \bar \lambda/D)$
\sn
\item "{{}}"  $\text{lim}_D(\bar \lambda) = \mu$
\sn
\item "{{}}"  $\lambda_x = \text{ cf}(\lambda_x)$.
\endroster
\endproclaim
\bigskip

\demo{Proof}  As in \cite{Sh:410},\S1 (replacing normal by $\sigma$-complete)
or make the following changes in the proof of \cite[Ch.IX]{Sh:g}, \scite{3.5}:
$\|N^x_k\| = \mu_k,\mu_k + 1 \subseteq N^x_k$.
\enddemo
\bigskip

\proclaim{\stag{6.14} Claim}  1) Assume
\mr
\widestnumber\item{$(iii)$}
\item "{$(i)$}"  $\lambda$ is inaccessible
\sn
\item "{$(ii)$}"  $|{\frak a}|^+ < \min({\frak a}),\mu =: \sup({\frak a}) <
\lambda$
\sn
\item "{$(iii)$}"  $R \subseteq \lambda \cap \text{ Reg } \backslash \mu,
|R| = \lambda$,
\sn
\item "{$(iv)$}"  for $\tau \in R,{\frak b}_\tau \subseteq {\frak a},\sup
({\frak b}_\tau) = \mu,J_\tau$ an ideal on ${\frak b}_\theta$ including
$J^{bd}_{{\frak b}_\theta},\tau = \text{ tcf}(\Pi {\frak b}_\theta,
<_{J \theta})$.
\ermn
\ub{Then} for some $\langle \lambda_\theta:\theta \in {\frak a} \rangle$ we
have
\mr
\item "{$(a)$}"  $\lambda_\theta = \text{ cf}(\lambda_\theta) < \theta,
\lambda_\theta > |{\frak a}|$
\sn
\item "{$(b)$}"  lim$_{J_\tau} \lambda_\theta = \mu$
\sn
\item "{$(c)$}"  $\lambda = \text{ tcf}(\dsize \prod_{\theta \in {\frak a}}
\lambda_\theta,<_{J_{< \lambda}[{\frak a}]})$
\sn
\item "{$(d)$}"  $\{\min \text{ pcf}_{J_\tau}(\dsize \prod_{\theta \in
{\frak b}_\tau} \lambda_\theta,<_{J_\tau}): \tau \in R\}$ is unbounded in
$\lambda$.
\ermn
So ${\frak a}' =: \{\lambda_\theta:\theta \in {\frak a}\}$ satisfies (ii),
(iii), (iv) and
\mr
\widestnumber\item{$(i)^-$}
\item "{$(i)^-$}"  $|{\frak a}'|^+ < \min({\frak a}')$ instead
\sn
\item "{$(i)$}"  $|{\frak a}|^+ < min({\frak a})$, immaterial and
\sn
\item "{$(v)$}"  max pcf$({\frak a}') = \lambda$.
\ermn
2) Assume $(i)^-,(ii),(iii),(iv),(v)$ are satisfied by ${\frak a},R,
{\frak b}_\tau,J_\tau$.  Then for some $f_\tau \in \Pi {\frak b}_\tau$ for
$\tau \in R$ we have:
\mr
\item "{$(*)$}"  for every $g \in \Pi {\frak a}$ for some $\tau,g
\restriction {\frak b}_\tau <_{J_\tau} f_\tau$.
\endroster
\endproclaim
\bigskip

\demo{Proof}  1) By \cite[Ch.II]{Sh:g}, \sciteu{1.5A} we can find $\langle
\lambda_\theta:\theta \in {\frak a} \rangle$ satisfying (a), (b), (c).  If
(d) fails choose for each $\tau \in R,{\frak b}'_\tau \in J^+_\tau$ such
that $\chi_\tau = \text{ tcf}(\dsize \prod_{\theta \in {\frak b}'_\tau}
\lambda_\theta,<_{J_{\tau \restriction {\frak b}'_\tau}})$ is well defined
and equal to min pcf$_{J_\tau}(\dsize \prod_{\theta \in {\frak b}_\tau}
\lambda_\theta,<_{J_\tau})$.  So $\{\chi_\tau:\tau \in R\}$ is unbounded in
$\lambda$ hence for some $\chi,R' = \{\tau \in R:\chi_\tau = \chi\}$ is
unbounded in $\lambda$.  Hence we can find $\zeta^* < |{\frak a}|^+$ and
$\tau_\zeta \in R'$ for $\zeta < \zeta^*$ such that $\lambda \le
\text{ max pcf}\{\tau_\zeta:\zeta < \zeta^*\}$ (choose by induction on
$\zeta < |{\frak a}|^+,\tau_\zeta \in R',\tau_\zeta > \text{ max pcf}
\{\tau_\varepsilon:\varepsilon < \zeta\}$ and use localization).  Let 
$D_\zeta$ be an ultrafilter on ${\frak a},b'_{\tau_\zeta} \in D_\tau,
J_{\tau_\zeta} \cap D = \emptyset$, so tcf$(\Pi{\frak a},<_{D_\zeta}) =
\tau_\zeta$, and $E$ is an ultrafilter on $\{\tau_\zeta:\zeta < \zeta^*\}$
such that tcf$(\dsize \prod_\zeta \tau_\zeta,<_E) \ge \lambda$.  Let
$D = \{{\frak c} \subseteq {\frak a}:\{\zeta:{\frak c} \in D_\zeta\} \in
E\}$, so tcf$(\Pi {\frak a},<_D) \ge \lambda$ hence $D \cap J_{< \lambda}
[{\frak a}] = \emptyset$.  Apply this to $\langle \lambda_\theta:\theta \in 
{\frak a} \rangle$, by \cite[Ch.I]{Sh:g}, \sciteu{1.10} - \sciteu{1.11} we get
a contradiction. \nl
3) Let $\langle f^*_\alpha:\alpha < \lambda \rangle$ be 
$<_{J_{< \lambda}[{\frak a}]}$-increasing and cofinal.  Let $f_\alpha =
f^*_\alpha \restriction {\frak b}_\alpha$.   \hfill$\square_{\scite{6.14}}$
\enddemo
\bn
Concerning \cite{Sh:430},\scite{3.1} we comment
\proclaim{\stag{6.15} Theorem}  1) Assume $\lambda > \theta > \aleph_0$ are
regular and
\mr
\item "{$(*)_{\theta,\kappa}$}"  if ${\frak a} \subseteq \text{ Reg } \cap
\lambda \backslash \theta$ and $|{\frak a}| < \theta$ \ub{then} there are
$\zeta^* < \theta$ and ${\frak b}_\zeta \in J_{< \lambda}[{\frak a}]$ for
$\zeta < \zeta^*$ such that for every ${\frak b} \subseteq {\frak a},
|{\frak b}| < \kappa$ for some $\zeta < \zeta^*,{\frak b} \subseteq {\frak b}
_\zeta$.
\ermn
\ub{Then} the following conditions are equivalent:
\mr
\item "{$(A) = (A)_{\lambda,\theta,\kappa}$}"    for every $\mu < \lambda$
we have cov$(\mu,\theta,\kappa,2) < \lambda$
\sn
\item "{$(B) = (B)_{\lambda,\theta,\kappa}$}"  if $\mu < \lambda$ and
$a_\alpha \in {\Cal S}_{< \kappa}(\mu)$ for $\alpha < \lambda$ then for some
$W \subseteq \lambda$ of cardinality $\lambda$ we have $|\dbcu_{\alpha \in W}
a_\alpha| < \theta$
\sn
\item "{$(C) = (C)_{\lambda,\theta,\kappa}$}"  if $a_\alpha$ is a set of
cardinality $< \kappa$ for $\alpha < \lambda$ and $W_0 \subseteq \{\delta <
\lambda:\text{cf}(\delta) \ge \kappa\}$ \ub{then} for some stationary $W
\subseteq W_0$ and set $b$ of cardinality $< \theta$ we have 
$\langle a_\alpha \backslash b:\alpha \in W \rangle$ is a sequence of
pairwise disjoint sets.
\ermn
2) If $\lambda > \theta_1 > \theta_2 \ge \kappa > \aleph_0$ where $\lambda,
\kappa$ are regular, $(A)_{\lambda,\theta,\kappa} \Leftrightarrow
(B)_{\lambda,\theta_1,\kappa}$ and cov$(\theta_1,\theta_2,\kappa,2) <
\lambda$ \ub{then} $(A)_{\lambda,\theta_2,\kappa} \Leftrightarrow
(B)_{\lambda,\theta_2,\kappa}$ (so if for some $\theta_1,
(*)_{\theta_1,\kappa},\lambda > \theta_1 > \theta_2 = 
\text{ cf}(\theta) \ge \kappa$ and
cov$(\theta_1,\theta,\kappa,2) < \lambda$ then the conclusions holds).
\endproclaim
\bigskip

\demo{Proof}  1) Read the proof of \cite{Sh:430}, \scite{3.1} (which was
written in a way appropriate to this generalization), but defining the
$M_n, \langle N^n_\zeta:\zeta < \theta \rangle$, we omit clause (d), that is,
$N^n_\zeta \in {\frak A}_{\delta(*)}$ and instead demand
\mr
\item "{$(d)'$}"  for each $n$ we can find ${\Cal P}_n \subseteq [\theta]
^{< \theta}$ such that $(\forall a \in [\theta]^{< \kappa})(\exists b \in
{\Cal P})(a \subseteq b)$ and $\langle N^n_b:b \in {\Cal P}_n \rangle$ such
that $N^n_b \prec {\frak B},N^n_\zeta = \bigcup\{N^n_b:b \in {\Cal P}_n,
b \subseteq \zeta\}$ and $b_1 \subseteq b_2 \Rightarrow N^n_{b_1} \prec
N^n_{b_2}$.
\ermn
2) Left to the reader.
\enddemo
\bn
\ub{\stag{6.16} Question}:  Can we in \cite{Sh:430}, \scite{4.2}(1) weaken
clause $(\beta)$ in the conclusion to ``$\lambda_x > \mu_0$ for $D$-almost
all $x \in {\Cal Y}/e$" then we can weaken the hypothesis \cite{Sh:420},
\sciteu{6.1C} (was stated in \cite{Sh:430}, earlier version clear).

\newpage
    
REFERENCES.  
\bibliographystyle{lit-plain}
\bibliography{lista,listb,listx,listf,liste}

\shlhetal
\enddocument

\bye

%% file: mathdefs.tex
\expandafter\ifx\csname mathdefs.tex\endcsname\relax
  \expandafter\gdef\csname mathdefs.tex\endcsname{}
\else \message{Hey!  Apparently you were trying to
  \string\input{mathdefs.tex} twice.   This does not make sense.} 
\errmessage{Please edit your file (probably \jobname.tex) and remove
any duplicate ``\string\input'' lines} \fi




\catcode`\X=12\catcode`\@=11

\def\n@wcount{\alloc@0\count\countdef\insc@unt}
\def\n@wwrite{\alloc@7\write\chardef\sixt@@n}
\def\n@wread{\alloc@6\read\chardef\sixt@@n}
\def\r@s@t{\relax}\def\v@idline{\par}\def\@mputate#1/{#1}
\def\l@c@l#1X{\firstpart.#1}\def\gl@b@l#1X{#1}\def\t@d@l#1X{{}}

\def\crossrefs#1{\ifx\all#1\let\tr@ce=\all\else\def\tr@ce{#1,}\fi
   \n@wwrite\cit@tionsout\openout\cit@tionsout=\jobname.cit 
   \write\cit@tionsout{\tr@ce}\expandafter\setfl@gs\tr@ce,}
\def\setfl@gs#1,{\def\@{#1}\ifx\@\empty\let\next=\relax
   \else\let\next=\setfl@gs\expandafter\xdef
   \csname#1tr@cetrue\endcsname{}\fi\next}
\def\m@ketag#1#2{\expandafter\n@wcount\csname#2tagno\endcsname
     \csname#2tagno\endcsname=0\let\tail=\all\xdef\all{\tail#2,}
   \ifx#1\l@c@l\let\tail=\r@s@t\xdef\r@s@t{\csname#2tagno\endcsname=0\tail}\fi
   \expandafter\gdef\csname#2cite\endcsname##1{\expandafter
     \ifx\csname#2tag##1\endcsname\relax?\else\csname#2tag##1\endcsname\fi
     \expandafter\ifx\csname#2tr@cetrue\endcsname\relax\else
     \write\cit@tionsout{#2tag ##1 cited on page \folio.}\fi}
   \expandafter\gdef\csname#2page\endcsname##1{\expandafter
     \ifx\csname#2page##1\endcsname\relax?\else\csname#2page##1\endcsname\fi
     \expandafter\ifx\csname#2tr@cetrue\endcsname\relax\else
     \write\cit@tionsout{#2tag ##1 cited on page \folio.}\fi}
   \expandafter\gdef\csname#2tag\endcsname##1{\expandafter
      \ifx\csname#2check##1\endcsname\relax
      \expandafter\xdef\csname#2check##1\endcsname{}%
      \else\immediate\write16{Warning: #2tag ##1 used more than once.}\fi
      \multit@g{#1}{#2}##1/X%
      \write\t@gsout{#2tag ##1 assigned number \csname#2tag##1\endcsname\space
      on page \number\count0.}%
   \csname#2tag##1\endcsname}}
\def\multit@g#1#2#3/#4X{\def\t@mp{#4}\ifx\t@mp\empty%
      \global\advance\csname#2tagno\endcsname by 1 
      \expandafter\xdef\csname#2tag#3\endcsname
      {#1\number\csname#2tagno\endcsnameX}%
   \else\expandafter\ifx\csname#2last#3\endcsname\relax
      \expandafter\n@wcount\csname#2last#3\endcsname
      \global\advance\csname#2tagno\endcsname by 1 
      \expandafter\xdef\csname#2tag#3\endcsname
      {#1\number\csname#2tagno\endcsnameX}
      \write\t@gsout{#2tag #3 assigned number \csname#2tag#3\endcsname\space
      on page \number\count0.}\fi
   \global\advance\csname#2last#3\endcsname by 1
   \def\t@mp{\expandafter\xdef\csname#2tag#3/}%
   \expandafter\t@mp\@mputate#4\endcsname
   {\csname#2tag#3\endcsname\lastpart{\csname#2last#3\endcsname}}\fi}
\def\t@gs#1{\def\all{}\m@ketag#1e\m@ketag#1s\m@ketag\t@d@l p
   \m@ketag\gl@b@l r \n@wread\t@gsin
   \openin\t@gsin=\jobname.tgs \re@der \closein\t@gsin
   \n@wwrite\t@gsout\openout\t@gsout=\jobname.tgs }
\outer\def\localtags{\t@gs\l@c@l}
\outer\def\globaltags{\t@gs\gl@b@l}
\outer\def\newlocaltag#1{\m@ketag\l@c@l{#1}}
\outer\def\newglobaltag#1{\m@ketag\gl@b@l{#1}}

\newif\ifpr@ 
\def\m@kecs #1tag #2 assigned number #3 on page #4.%
   {\expandafter\gdef\csname#1tag#2\endcsname{#3}
   \expandafter\gdef\csname#1page#2\endcsname{#4}
   \ifpr@\expandafter\xdef\csname#1check#2\endcsname{}\fi}
\def\re@der{\ifeof\t@gsin\let\next=\relax\else
   \read\t@gsin to\t@gline\ifx\t@gline\v@idline\else
   \expandafter\m@kecs \t@gline\fi\let \next=\re@der\fi\next}
\def\pretags#1{\pr@true\pret@gs#1,,}
\def\pret@gs#1,{\def\@{#1}\ifx\@\empty\let\n@xtfile=\relax
   \else\let\n@xtfile=\pret@gs \openin\t@gsin=#1.tgs \message{#1} \re@der 
   \closein\t@gsin\fi \n@xtfile}

\newcount\sectno\sectno=0\newcount\subsectno\subsectno=0
\newif\ifultr@local \def\ultralocal{\ultr@localtrue}
\def\firstpart{\number\sectno}
\def\lastpart#1{\ifcase#1 \or a\or b\or c\or d\or e\or f\or g\or h\or 
   i\or k\or l\or m\or n\or o\or p\or q\or r\or s\or t\or u\or v\or w\or 
   x\or y\or z \fi}

\def\resetall{\global\advance\sectno by 1\subsectno=0
   \gdef\firstpart{\number\sectno}\r@s@t}
\def\resetsub{\global\advance\subsectno by 1
   \gdef\firstpart{\number\sectno.\number\subsectno}\r@s@t}
\def\newsection#1\par{\resetall\vskip0pt plus.3\vsize\penalty-250
   \vskip0pt plus-.3\vsize\bigskip\bigskip
   \message{#1}\leftline{\bf#1}\nobreak\bigskip}
\def\subsection#1\par{\ifultr@local\resetsub\fi
   \vskip0pt plus.2\vsize\penalty-250\vskip0pt plus-.2\vsize
   \bigskip\smallskip\message{#1}\leftline{\bf#1}\nobreak\medskip}

\def\t@gsoff#1,{\def\@{#1}\ifx\@\empty\let\next=\relax\else\let\next=\t@gsoff
   \def\@@{p}\ifx\@\@@\else
   \expandafter\gdef\csname#1cite\endcsname##1{\zeigen{##1}}
   \expandafter\gdef\csname#1page\endcsname##1{?}
   \expandafter\gdef\csname#1tag\endcsname##1{\zeigen{##1}}\fi\fi\next}
\def\verbatimtags{\ifx\all\relax\else\expandafter\t@gsoff\all,\fi}
\def\zeigen#1{\hbox{$\langle$}#1\hbox{$\rangle$}}

\def\(#1){\edef\dot@g{\ifmmode\ifinner(\hbox{\noexpand\etag{#1}})
   \else\noexpand\eqno(\hbox{\noexpand\etag{#1}})\fi
   \else(\noexpand\ecite{#1})\fi}\dot@g}

\newif\ifbr@ck
\def\eat#1{}
\def\[#1]{\br@cktrue[\br@cket#1'X]}
\def\br@cket#1'#2X{\def\temp{#2}\ifx\temp\empty\let\next\eat
   \else\let\next\br@cket\fi
   \ifbr@ck\br@ckfalse\br@ck@t#1,X\else\br@cktrue#1\fi\next#2X}
\def\br@ck@t#1,#2X{\def\temp{#2}\ifx\temp\empty\let\neext\eat
   \else\let\neext\br@ck@t\def\temp{,}\fi
   \def\teemp{#1}\ifx\teemp\empty\else\rcite{#1}\fi\temp\neext#2X}
\def\resetbr@cket{\gdef\[##1]{[\rtag{##1}]}}
\def\references{\resetbr@cket\newsection References\par}

\newtoks\symb@ls\newtoks\s@mb@ls\newtoks\p@gelist\n@wcount\ftn@mber
    \ftn@mber=1\newif\ifftn@mbers\ftn@mbersfalse\newif\ifbyp@ge\byp@gefalse
\def\defm@rk{\ifftn@mbers\n@mberm@rk\else\symb@lm@rk\fi}
\def\n@mberm@rk{\xdef\m@rk{{\the\ftn@mber}}%
    \global\advance\ftn@mber by 1 }
\def\rot@te#1{\let\temp=#1\global#1=\expandafter\r@t@te\the\temp,X}
\def\r@t@te#1,#2X{{#2#1}\xdef\m@rk{{#1}}}
\def\b@@st#1{{$^{#1}$}}\def\str@p#1{#1}
\def\symb@lm@rk{\ifbyp@ge\rot@te\p@gelist\ifnum\expandafter\str@p\m@rk=1 
    \s@mb@ls=\symb@ls\fi\write\f@nsout{\number\count0}\fi \rot@te\s@mb@ls}
\def\byp@ge{\byp@getrue\n@wwrite\f@nsin\openin\f@nsin=\jobname.fns 
    \n@wcount\currentp@ge\currentp@ge=0\p@gelist={0}
    \re@dfns\closein\f@nsin\rot@te\p@gelist
    \n@wread\f@nsout\openout\f@nsout=\jobname.fns }
\def\m@kelist#1X#2{{#1,#2}}
\def\re@dfns{\ifeof\f@nsin\let\next=\relax\else\read\f@nsin to \f@nline
    \ifx\f@nline\v@idline\else\let\t@mplist=\p@gelist
    \ifnum\currentp@ge=\f@nline
    \global\p@gelist=\expandafter\m@kelist\the\t@mplistX0
    \else\currentp@ge=\f@nline
    \global\p@gelist=\expandafter\m@kelist\the\t@mplistX1\fi\fi
    \let\next=\re@dfns\fi\next}
\def\symbols#1{\symb@ls={#1}\s@mb@ls=\symb@ls} 
\def\bigsymbol{\textstyle}
\symbols{\bigsymbol\ast,\dagger,\ddagger,\sharp,\flat,\natural,\star}
\def\ftnumbers{\ftn@mberstrue} \def\ftsymbols{\ftn@mbersfalse}
\def\paginal{\byp@ge} \def\resetftnumbers{\ftn@mber=1}
\def\ftnote#1{\defm@rk\expandafter\expandafter\expandafter\footnote
    \expandafter\b@@st\m@rk{#1}}

\long\def\jump#1\endjump{}
\def\ssum{\mathop{\lower .1em\hbox{$\textstyle\Sigma$}}\nolimits}

\def\qed{\nobreak\kern 1em \vrule height .5em width .5em depth 0em}
\def\newneq{\hbox{\rlap{\hbox to 1\wd9{\hss$=$\hss}}\raise .1em 
   \hbox to 1\wd9{\hss$\scriptscriptstyle/$\hss}}}
\def\subsetne{\setbox9 = \hbox{$\subset$}\mathrel{\hbox{\rlap
   {\lower .4em \newneq}\raise .13em \hbox{$\subset$}}}}
\def\supsetne{\setbox9 = \hbox{$\subset$}\mathrel{\hbox{\rlap
   {\lower .4em \newneq}\raise .13em \hbox{$\supset$}}}}

\def\vbar{\mathchoice{\vrule height6.3ptdepth-.5ptwidth.8pt\kern-.8pt}
   {\vrule height6.3ptdepth-.5ptwidth.8pt\kern-.8pt}
   {\vrule height4.1ptdepth-.35ptwidth.6pt\kern-.6pt}
   {\vrule height3.1ptdepth-.25ptwidth.5pt\kern-.5pt}}
\def\f@dge{\mathchoice{}{}{\mkern.5mu}{\mkern.8mu}}
\def\b@c#1#2{{\rm \mkern#2mu\vbar\mkern-#2mu#1}}
\def\b@b#1{{\rm I\mkern-3.5mu #1}}
\def\b@a#1#2{{\rm #1\mkern-#2mu\f@dge #1}}
\def\bb#1{{\count4=`#1 \advance\count4by-64 \ifcase\count4\or\b@a A{11.5}\or
   \b@b B\or\b@c C{5}\or\b@b D\or\b@b E\or\b@b F \or\b@c G{5}\or\b@b H\or
   \b@b I\or\b@c J{3}\or\b@b K\or\b@b L \or\b@b M\or\b@b N\or\b@c O{5} \or
   \b@b P\or\b@c Q{5}\or\b@b R\or\b@a S{8}\or\b@a T{10.5}\or\b@c U{5}\or
   \b@a V{12}\or\b@a W{16.5}\or\b@a X{11}\or\b@a Y{11.7}\or\b@a Z{7.5}\fi}}

\catcode`\X=11 \catcode`\@=12

%% file: citeadd.tex

\expandafter\ifx\csname citeadd.tex\endcsname\relax
\expandafter\gdef\csname citeadd.tex\endcsname{}
\else \message{Hey!  Apparently you were trying to
\string\input{citeadd.tex} twice.   This does not make sense.} 
\errmessage{Please edit your file (probably \jobname.tex) and remove
any duplicate ``\string\input'' lines} \fi

\def\sciteu{\sciteerror{undefined}}

\def\sciteerror#1#2{{\mathortextbf{\scite{#2}}}\complainaboutcitation{#1}{#2}}
\def\mathortextbf#1{\hbox{\bf #1}}
\def\complainaboutcitation#1#2{%
\vadjust{\line{\llap{---$\!\!>$ }\qquad scite$\{$#2$\}$ #1\hfil}}}

%% file: alice2jlem.tex

\expandafter\ifx\csname alice2jlem.tex\endcsname\relax
  \expandafter\gdef\csname alice2jlem.tex\endcsname{}
\else \message{Hey!  Apparently you were trying to
\string\input{alice2jlem.tex}  twice.   This does not make sense.}
\errmessage{Please edit your file (probably \jobname.tex) and remove
any duplicate ``\string\input'' lines} \fi

\input bib4plain

\def\widestnumber#1#2{}

\def\rm{\fam0 \tenrm}

\def\fakesubhead#1\endsubhead{\bigskip\noindent{\bf#1}\par}

\input rsfs

%% file: bib4plain.tex
\expandafter\ifx\csname bib4plain.tex\endcsname\relax
  \expandafter\gdef\csname bib4plain.tex\endcsname{}
\else \message{Hey!  Apparently you were trying to \string\input
  bib4plain.tex twice.   This does not make sense.}
\errmessage{Please edit your file (probably \jobname.tex) and remove
any duplicate ``\string\input'' lines} \fi

\def\renewcommand{\newcommand}	       
\edef\cite{\the\catcode`@}%
\catcode`@ = 11
\let\@oldatcatcode = \cite
\chardef\@letter = 11
\chardef\@other = 12
%
%
%
%
\def\@innerdef#1#2{\edef#1{\expandafter\noexpand\csname #2\endcsname}}%
%
%
\@innerdef\@innernewcount{newcount}%
\@innerdef\@innernewdimen{newdimen}%
\@innerdef\@innernewif{newif}%
\@innerdef\@innernewwrite{newwrite}%
%
%
%
\def\@gobble#1{}%
%
%
%
\ifx\inputlineno\@undefined
   \let\@linenumber = \empty 
\else
   \def\@linenumber{\the\inputlineno:\space}%
\fi
%
%
%
\def\@futurenonspacelet#1{\def\cs{#1}%
   \afterassignment\@stepone\let\@nexttoken=
}%
\begingroup 
\def\\{\global\let\@stoken= }%
\\ 
\endgroup
\def\@stepone{\expandafter\futurelet\cs\@steptwo}%
\def\@steptwo{\expandafter\ifx\cs\@stoken\let\@@next=\@stepthree
   \else\let\@@next=\@nexttoken\fi \@@next}%
\def\@stepthree{\afterassignment\@stepone\let\@@next= }%
%
%
%
\def\@getoptionalarg#1{%
   \let\@optionaltemp = #1%
   \let\@optionalnext = \relax
   \@futurenonspacelet\@optionalnext\@bracketcheck
}%
%
%
\def\@bracketcheck{%
   \ifx [\@optionalnext
      \expandafter\@@getoptionalarg
   \else
      \let\@optionalarg = \empty
      \expandafter\@optionaltemp
   \fi
}%
\def\@@getoptionalarg[#1]{%
   \def\@optionalarg{#1}%
   \@optionaltemp
}%
%
%
%
\def\@nnil{\@nil}%
\def\@fornoop#1\@@#2#3{}%
\def\@for#1:=#2\do#3{%
   \edef\@fortmp{#2}%
   \ifx\@fortmp\empty \else
      \expandafter\@forloop#2,\@nil,\@nil\@@#1{#3}%
   \fi
}%
\def\@forloop#1,#2,#3\@@#4#5{\def#4{#1}\ifx #4\@nnil \else
       #5\def#4{#2}\ifx #4\@nnil \else#5\@iforloop #3\@@#4{#5}\fi\fi
}%
\def\@iforloop#1,#2\@@#3#4{\def#3{#1}\ifx #3\@nnil
       \let\@nextwhile=\@fornoop \else
      #4\relax\let\@nextwhile=\@iforloop\fi\@nextwhile#2\@@#3{#4}%
}%
%
%
%
\@innernewif\if@fileexists
\def\@testfileexistence{\@getoptionalarg\@finishtestfileexistence}%
\def\@finishtestfileexistence#1{%
   \begingroup
      \def\extension{#1}%
      \immediate\openin0 =
         \ifx\@optionalarg\empty\jobname\else\@optionalarg\fi
         \ifx\extension\empty \else .#1\fi
         \space
      \ifeof 0
         \global\@fileexistsfalse
      \else
         \global\@fileexiststrue
      \fi
      \immediate\closein0
   \endgroup
}%
%
%
%
%
\def\bibliographystyle#1{%
   \@readauxfile
   \@writeaux{\string\bibstyle{#1}}%
}%
\let\bibstyle = \@gobble
%
%
\let\bblfilebasename = \jobname
\def\bibliography#1{%
   \@readauxfile
   \@writeaux{\string\bibdata{#1}}%
   \@testfileexistence[\bblfilebasename]{bbl}%
   \if@fileexists
      \nobreak
      \@readbblfile
   \fi
}%
\let\bibdata = \@gobble
%
%
\def\nocite#1{%
   \@readauxfile
   \@writeaux{\string\citation{#1}}%
}%
\@innernewif\if@notfirstcitation
%
%
\def\cite{\@getoptionalarg\@cite}%
%
%
\def\@cite#1{%
   \let\@citenotetext = \@optionalarg
   \printcitestart
   \nocite{#1}%
   \@notfirstcitationfalse
   \@for \@citation :=#1\do
   {%
      \expandafter\@onecitation\@citation\@@
   }%
   \ifx\empty\@citenotetext\else
      \printcitenote{\@citenotetext}%
   \fi
   \printcitefinish
}%
\def\@onecitation#1\@@{%
   \if@notfirstcitation
      \printbetweencitations
   \fi
   \expandafter \ifx \csname\@citelabel{#1}\endcsname \relax
      \if@citewarning
         \message{\@linenumber Undefined citation `#1'.}%
      \fi
      \expandafter\gdef\csname\@citelabel{#1}\endcsname{%
\strut
\vadjust{\vskip-\dp\strutbox
\vbox to 0pt{\vss\parindent0cm \leftskip=\hsize 
\advance\leftskip3mm
\advance\hsize 4cm\strut\openup-4pt 
\rightskip 0cm plus 1cm minus 0.5cm ?  #1 ?\strut}}
         {\tt
            \escapechar = -1
            \nobreak\hskip0pt
            \expandafter\string\csname#1\endcsname
            \nobreak\hskip0pt
         }%
      }%
   \fi
   \csname\@citelabel{#1}\endcsname
   \@notfirstcitationtrue
}%
%
%
\def\@citelabel#1{b@#1}%
%
%
\def\@citedef#1#2{\expandafter\gdef\csname\@citelabel{#1}\endcsname{#2}}%
%
%
%
\def\@readbblfile{%
   \ifx\@itemnum\@undefined
      \@innernewcount\@itemnum
   \fi
   \begingroup
      \def\begin##1##2{%
         \setbox0 = \hbox{\biblabelcontents{##2}}%
         \biblabelwidth = \wd0
      }%
      \def\end##1{}
      %
      %
      \@itemnum = 0
      \def\bibitem{\@getoptionalarg\@bibitem}%
      \def\@bibitem{%
         \ifx\@optionalarg\empty
            \expandafter\@numberedbibitem
         \else
            \expandafter\@alphabibitem
         \fi
      }%
      \def\@alphabibitem##1{%
         \expandafter \xdef\csname\@citelabel{##1}\endcsname {\@optionalarg}%
         \ifx\biblabelprecontents\@undefined
            \let\biblabelprecontents = \relax
         \fi
         \ifx\biblabelpostcontents\@undefined
            \let\biblabelpostcontents = \hss
         \fi
         \@finishbibitem{##1}%
      }%
      \def\@numberedbibitem##1{%
         \advance\@itemnum by 1
         \expandafter \xdef\csname\@citelabel{##1}\endcsname{\number\@itemnum}%
         \ifx\biblabelprecontents\@undefined
            \let\biblabelprecontents = \hss
         \fi
         \ifx\biblabelpostcontents\@undefined
            \let\biblabelpostcontents = \relax
         \fi
         \@finishbibitem{##1}%
      }%
      \def\@finishbibitem##1{%
         \biblabelprint{\csname\@citelabel{##1}\endcsname}%
         \@writeaux{\string\@citedef{##1}{\csname\@citelabel{##1}\endcsname}}%
         \ignorespaces
      }%
      %
      %
      \let\em = \bblem
      \let\newblock = \bblnewblock
      \let\sc = \bblsc
      \frenchspacing
      \clubpenalty = 4000 \widowpenalty = 4000
      \tolerance = 10000 \hfuzz = .5pt
      \everypar = {\hangindent = \biblabelwidth
                      \advance\hangindent by \biblabelextraspace}%
      \bblrm
      \parskip = 1.5ex plus .5ex minus .5ex
      \biblabelextraspace = .5em
      \bblhook
      \input \bblfilebasename.bbl
   \endgroup
}%
%
%
\@innernewdimen\biblabelwidth
\@innernewdimen\biblabelextraspace
%
%
%
\def\biblabelprint#1{%
   \noindent
   \hbox to \biblabelwidth{%
      \biblabelprecontents
      \biblabelcontents{#1}%
      \biblabelpostcontents
   }%
   \kern\biblabelextraspace
}%
%
%
%
\def\biblabelcontents#1{{\bblrm [#1]}}%
%
%
\def\bblrm{\rm}%
%
%
\def\bblem{\it}%
%
%
\def\bblsc{\ifx\@scfont\@undefined
              \font\@scfont = cmcsc10
           \fi
           \@scfont
}%
%
%
\def\bblnewblock{\hskip .11em plus .33em minus .07em }%
%
%
\let\bblhook = \empty
%
%
%
\def\printcitestart{[}
\def\printcitefinish{]}
\def\printbetweencitations{, }
\def\printcitenote#1{, #1}
%
%
%
\let\citation = \@gobble
%
%
%
\@innernewcount\@numparams
%
%
\def\newcommand#1{%
   \def\@commandname{#1}%
   \@getoptionalarg\@continuenewcommand
}%
%
%
\def\@continuenewcommand{%
   \@numparams = \ifx\@optionalarg\empty 0\else\@optionalarg \fi \relax
   \@newcommand
}%
%
%
\def\@newcommand#1{%
   \def\@startdef{\expandafter\edef\@commandname}%
   \ifnum\@numparams=0
      \let\@paramdef = \empty
   \else
      \ifnum\@numparams>9
         \errmessage{\the\@numparams\space is too many parameters}%
      \else
         \ifnum\@numparams<0
            \errmessage{\the\@numparams\space is too few parameters}%
         \else
            \edef\@paramdef{%
               \ifcase\@numparams
                  \empty  No arguments.
               \or ####1%
               \or ####1####2%
               \or ####1####2####3%
               \or ####1####2####3####4%
               \or ####1####2####3####4####5%
               \or ####1####2####3####4####5####6%
               \or ####1####2####3####4####5####6####7%
               \or ####1####2####3####4####5####6####7####8%
               \or ####1####2####3####4####5####6####7####8####9%
               \fi
            }%
         \fi
      \fi
   \fi
   \expandafter\@startdef\@paramdef{#1}%
}%
%
%
%
%
\def\@readauxfile{%
   \if@auxfiledone \else 
      \global\@auxfiledonetrue
      \@testfileexistence{aux}%
      \if@fileexists
         \begingroup
            \endlinechar = -1
            \catcode`@ = 11
            \input \jobname.aux
         \endgroup
      \else
         \message{\@undefinedmessage}%
         \global\@citewarningfalse
      \fi
      \immediate\openout\@auxfile = \jobname.aux
   \fi
}%
%
%
\newif\if@auxfiledone
\ifx\noauxfile\@undefined \else \@auxfiledonetrue\fi
%
%
%
%
\@innernewwrite\@auxfile
\def\@writeaux#1{\ifx\noauxfile\@undefined \write\@auxfile{#1}\fi}%
%
%
%
\ifx\@undefinedmessage\@undefined
   \def\@undefinedmessage{No .aux file; I won't give you warnings about
                          undefined citations.}%
\fi
%
%
\@innernewif\if@citewarning
\ifx\noauxfile\@undefined \@citewarningtrue\fi
%
%
%
\catcode`@ = \@oldatcatcode

%% file: rsfs.tex

%
%
%

%

\font\textrsfs=rsfs10
\font\scriptrsfs=rsfs7
\font\scriptscriptrsfs=rsfs5

\newfam\rsfsfam
\textfont\rsfsfam=\textrsfs
\scriptfont\rsfsfam=\scriptrsfs
\scriptscriptfont\rsfsfam=\scriptscriptrsfs

\edef\oldcatcodeofat{\the\catcode`\@}
\catcode`\@11

\def\Cal@@#1{\noaccents@ \fam \rsfsfam #1}

\catcode`\@\oldcatcodeofat